\documentclass[review,12pt]{elsarticle}
\biboptions{authoryear,round}

\usepackage{hyperref}

\usepackage{amsmath,amssymb}
\usepackage{booktabs}

\usepackage{algorithm}
\usepackage{algpseudocode}
\usepackage{multirow}
\usepackage{graphicx}
\usepackage{subcaption}
\usepackage{makecell}
\usepackage{threeparttable}
\usepackage{array}
\usepackage{float}
\usepackage{setspace}

\journal{Transport Policy}
\begin{document}

\begin{frontmatter}

\title{Assessing Autonomous Mobility-on-Demand Services and the Impacts of Operational Strategies: A Case Study of Chengdu, China}

\author[scu,nhh]{Youkai Wu}
\author[scu]{Zhaoxia Guo}
\author[scu]{Qi Liu\corref{cor1}}
\author[scu,nhh]{Stein W. Wallace}

\cortext[cor1]{Corresponding author}

\address[scu]{Business School, Sichuan University, Chengdu, 610065, China}
\address[nhh]{NHH Norwegian School of Economics, Bergen, Norway}

\begin{abstract}
The Autonomous Mobility-on-Demand (AMoD) service is emerging as a potential alternative to on-demand urban mobility, but its operational performance relative to traditional street-hailing services and the effectiveness of related operational strategies remain unclear. This study presents a simulation framework integrating a graph theory-based trip-vehicle matching mechanism and uses historical street-hailing operations data to simulate AMoD services in Chengdu, China. The operational performance of these two urban mobility modes is evaluated using three key performance indicators: average passenger waiting time (APWT), average deadheading mileage (ADM), and average deadheading energy consumption (ADEC). We further evaluate the impacts of four operational strategies on simulated AMoD performance: vehicle repositioning, fleet size management, geofencing, and request rejection. Simulation results indicate that, under the same historical trip demand, fleet-size constraints, and road network as the observed street-hailing system, the simulated AMoD service is estimated to have lower values of APWT, ADM, and ADEC by 73.3\% -- 83.4\%, 75.0\%, and 74.0\%, respectively, reflecting the potential operational gains associated with centralized dispatch in simulation settings. These differences are most pronounced during early-morning low-demand hours and in remote areas such as airports. Operational strategy analysis indicates that for AMoD services, (1) vehicle repositioning shortens APWT at the cost of increased ADM and ADEC; (2) fleet expansion yields diminishing marginal returns; (3) geofencing improves the performance of city-centre services but increases unserved requests and worsens service performance for other zones, resulting in poorer overall system performance; and (4) rejecting service requests destined for high-attraction but low-demand areas yields significant improvements in all performance indicators among accepted trips.
\end{abstract}


\end{frontmatter}

\section{Introduction}

Autonomous Mobility-on-Demand (AMoD) services, such as those operated by Waymo and Apollo Go, are an emerging mode of on-demand urban transportation. Unlike traditional on-demand mobility service modes such as street-hailing and ride-hailing services which are provided by human drivers and affected by drivers’ preferences and cruising behaviours~\citep{Rayle2016,Aguilera-Garcia2022,NahmiasBiran2019TRR}, AMoD refers to the on-demand mobility system in which a centrally managed fleet of autonomous vehicles (AVs) is dispatched to serve passenger trip requests~\citep{Oh2020_AMOD_Singapore,siddiq2022ride}. In recent years, AMoD services have expanded to a growing number of urban areas. Since the first AMoD service was offered in Singapore in 2016, AMoD services have expanded to around 30 major cities around the world~\citep{Waymo2025YIR,Yang2024Robotaxi}. As an example, Apollo Go delivered 2.2 million fully driverless rides in Q2 of 2025, representing a 148\% year-over-year increase~\citep{Baidu2025HKEXQ2}. In light of the rapid growth of AMoD services, assessing the performance of AMoD services and their comparison with street-hailing services is important for assessing the economic feasibility and strategy implications of AMoD services and further developing targeted operational strategies.

Most existing comparative studies between AMoD and traditional taxi services focus on surveys about passenger travel preferences and acceptance or cost analysis of the two modes. Some studies compare the operational performance of these two modes by explicitly modelling traditional taxi drivers’ behavioural chains and conducting experiments using simulated demand data~\citep{NahmiasBiran2019TRR}; others analyse two modes in the context of ride-sharing~\citep{Lokhandwala2018TRC}. Limited work has compared the operational performance of AMoD and street-hailing services under historical and identical trips and supply conditions and systematically evaluated the potential operational strategies and their impacts on AMoD services. Given these gaps, this study focuses on two major questions for AMoD platforms:
\begin{enumerate}
  \item How does the operational performance of the AMoD system differ from that of traditional street-hailing taxi services under the same historical trips, fleet size and road network?
  \item How do operational strategies that adjust repositioning, fleet size, geofencing, and request rejection affect the performance of AMoD services?
\end{enumerate}

To answer these questions, this study conducts a simulation-based evaluation, taking Chengdu, a megacity in China, as a case city. First, we develop a simulation framework grounded in real street-hailing taxi operations data to evaluate and compare AMoD and street-hailing services. Second, we formulate the minimization of the trip pickup time as the objective function and employ theoretically optimal and efficient graph theory-based trip-vehicle matching mechanisms to simulate AMoD services within rolling time windows. We further test four operational strategies: vehicle repositioning, fleet size management, geofencing, and request rejection, using this framework to analyse their impacts on AMoD performance.

This study contributes to the literature in the following ways: i) it compares the performance of AMoD and street-hailing services by using the same historical passenger trips, fleet-size constraints, and urban road network, ii) it identifies how AMoD performance varies across different spatial and temporal conditions, and iii) it simulates and evaluates the impacts of four representative operational strategies (i.e., vehicle repositioning, fleet size management, geofencing, and request rejection) on the performance of AMoD services.

The remainder of the paper is organized as follows. Section~\ref{sec:lit_review} reviews the relevant literature. Section~\ref{sec:method} describes the methodology and data used in this study. Section~\ref{sec:results} compares results of AMoD and street-hailing services based on spatiotemporal analyses. Section~\ref{sec:discussion} discusses the effects of different operational strategies. Section~\ref{sec:conclusion} summarizes the main conclusions and outlines future research directions.

\section{Literature review}\label{sec:lit_review}
The emergence of AMoD has prompted extensive debate and attracted increasing attention from researchers and policymakers~\citep{Zhou2023RTE,Pande2023JEPR}. This section reviews previous studies on (i) the comparison of AMoD and street-hailing services, and (ii) operational strategies affecting AMoD services.

\subsection{Comparison of AMoD and street-hailing services}
Numerous studies have evaluated and compared the performance of on-demand mobility services based on different analytical methods, which can be mainly divided into survey-based analyses, empirical data analyses, and simulation model analyses. Some researchers adopted survey-based analyses to study passenger choice patterns between AMoD and other mobility service modes across on-demand mobility services~\citep{Asgari2018TRR,Pakusch2020TBS,Pemberton2025BIT,Mo2021TRC}. However, research based on surveys suffers from sampling and response bias. Due to the subjective nature of passenger perceptions, these studies provided qualitative insights and needed to be extended to quantitative analysis.

To obtain quantitative insights into on-demand mobility services, a number of researchers used empirical data analyses to compare different service modes based on real-world operations data.~\cite{Bosch2018TP} developed a bottom-up cost-accounting model to evaluate cost-competitiveness under different settings, concluding that AMoD services are cheaper in lower-density demand settings.~\cite{Freedman2018IE} used a Markov chain–based simulation and cost-effectiveness framework and found that AV deployment is more likely to be cost-effective in taxi-like, high-utilization service.~\cite{Wadud2021TRA} computed a generalized total cost of ownership using household travel data, and then analysed how costs and socio-economic factors dominate choice patterns. Related empirical studies have also compared service quality between app-based ride-hailing and traditional taxis.~\cite{Brown2021Transportation} examined ride-hailing and taxi trips in Los Angeles and found that ride-hailing users waited substantially less than taxi users for comparable origin–destination pairs. Although app-based ride-hailing is not equivalent to AMoD, this evidence suggests that centralized dispatching can lead to substantial differences in passenger waiting time relative to decentralized taxi services.

To evaluate the performance of AMoD services, some researchers adopted simulation model analyses using the same dataset to evaluate the performance of different on-demand mobility services.~\cite{Vazifeh2018} and \cite{Zhan2016TITS} modelled the operation of AVs based on graph theory using data from street-hailing taxis in New York and addressed the minimum fleet problem in on-demand mobility. They found respectively that AMoD needs 30\% fewer vehicles to serve most orders of street-hailing services and achieves at most 87\% better performance on deadheading mileages than street-hailing.~\cite{FengKongWang2021} constructed a queuing network on a closed circular road and simulated the efficiency of point-to-point on‑demand matching system under varying demand levels and travel distances. They found that under medium demand and long travel distances, their point-to-point on‑demand system performed worse than traditional street‑hailing services.~\cite{NahmiasBiran2019TRR} proposed an event-based, agent-based simulation framework in SimMobility that jointly models decentralized traditional MoD driver decision-making and centralized AMoD fleet management. Using a Singapore case study, they demonstrated how switching from traditional MoD to AMoD changes service performance and required fleet size. These studies compared different on-demand mobility services using the same datasets. However, they either adopted simplified or synthetic spatial environments that did not fully reflect real-world road network constraints or did not evaluate different service modes in the same supply and demand contexts.

\subsection{Operational strategies affecting AMoD services}

Previous studies have investigated a range of operational strategies and regulatory measures affecting AMoD services from both the supply and demand sides.

Supply-side research has primarily examined the operational behaviour of platforms and drivers, including fleet sizing, vehicle repositioning, ride-sharing, and service-area restrictions. Studies showed that a reduced fleet size can serve a similar number of trips under centralized dispatch compared with street-hailing~\citep{NahmiasBiran2019TRR,Vazifeh2018}. Platform-level interventions, such as limiting autonomous vehicle service areas, have also been modelled to explore equilibrium solutions in mixed or monopolistic markets~\citep{siddiq2022ride} or tested on the accessibility and network performance~\citep{Nahmias-Biran2021Transportation}. However, studies on geofencing focused on restricting the operational region of AMoD services to the city centre, rather than examining citywide service-area partitioning strategies. Another important research question in the operational design and deployment of AMoD systems concerns fleet sizing, dynamic ride-sharing, shared-ride operations, and integration with public transit. \cite{FagnantKockelman2018Transportation} investigated dynamic ride-sharing and fleet-sizing decisions for a system of shared autonomous vehicles in Austin, USA. Their results showed that dynamic ride-sharing could reduce passenger service times and travel costs while limiting the additional vehicle miles travelled associated with shared autonomous vehicle operations. Using an agent-based stochastic dynamic simulation framework, \cite{HylandMahmassani2020TRA} analysed the operational benefits and challenges of shared-ride AMoD services and demonstrated that ride-sharing strategies, fleet size, and other operational parameters can substantially affect system efficiency. In addition to stand-alone AMoD services, \cite{Cai2025TIMOD} investigated the incorporation of mobility-on-demand services into public transit in suburban areas and compared their cost-effectiveness with fixed-route bus services, private driving, and commercial ride-hailing. It highlights the potential role of mobility-on-demand services as a complement to public transit, particularly in suburban and low-density areas.

Demand-side research has mainly examined how pricing can be used to guide, regulate, stimulate, or restrict passenger travel demand using pricing as the tool.~\cite{SimoniKockelmanGurumurthyBischoff2019} used agent-based transport simulator MATSim with multiple future AV scenarios to quantify impacts of pricing on congestion and traveller welfare, finding that all pricing schemes reduce congestion but welfare outcomes differ materially.~\cite{SalazarLanzettiRossiSchifferPavone2020} formulated a public-transit network flow model that maximizes social welfare and derived a pricing and tolling policy that can steer selfish agents to the social optimum. Overall, research on demand management has mainly focused on the effectiveness of policy tools, such as dynamic pricing or congestion charges, whereas limited attention has been devoted to examining the effects of request rejection strategies on AMoD system performance and to identifying the trip types that such operational strategies should target. However, existing studies often examine one operational strategy at a time, and few compare several supply- and demand-side strategies using the same framework. This makes it difficult to understand which strategies work better and what trade-offs they involve. Studies on AMoD geofencing have also mainly restricted services to a specific area, such as the city centre, rather than dividing the whole city into service zones. Demand-side research has focused mostly on pricing, with less attention to rejecting requests based on where and when trips occur. This study therefore uses a consistent simulation framework to examine different operational strategies, including vehicle repositioning, fleet size, citywide geofencing, and selective request rejection.

In summary, the major research gaps are as follows. First, there have been very few studies that directly compare the performance of AMoD and street-hailing on the same urban road network under identical passenger trips and fleet-size constraints, even though this is crucial for evaluating the effects of AMoD services. Second, most existing studies assess geofencing impacts on AMoD services by restricting the operational region to the city centre, while demand management research mainly focuses on designing pricing tools for peak shaving. Previous studies have not evaluated citywide geofencing nor examined which types of demand impose greater operational burdens on AMoD services. And few studies have evaluated supply- and demand-side operational strategies within the same framework. To fill these research gaps, referring to the research~\cite{Vazifeh2018}, we employ a simulation-based model to evaluate the performance of AMoD services and benchmark it against observed street-hailing performance derived from the same taxi operations dataset. We also evaluate the individual impacts of vehicle repositioning, fleet size adjustments, geofencing, and request rejection strategies on the operational performance of AMoD services.

\section{Methodology and data}\label{sec:method}
In this section, we propose an AMoD simulation framework in which trips and autonomous vehicles are matched in every time period (e.g., each minute), using a graph theory-based mechanism. Based on this framework, we conduct a case study to compare the performance of AMoD and street-hailing services on three key performance indicators (KPIs) using a real street-hailing taxi operations dataset from the road network of Chengdu city, China.

\subsection{The AMoD simulation framework}
\begin{figure}[!htbp]
\centering
\includegraphics[width=1.0\textwidth]{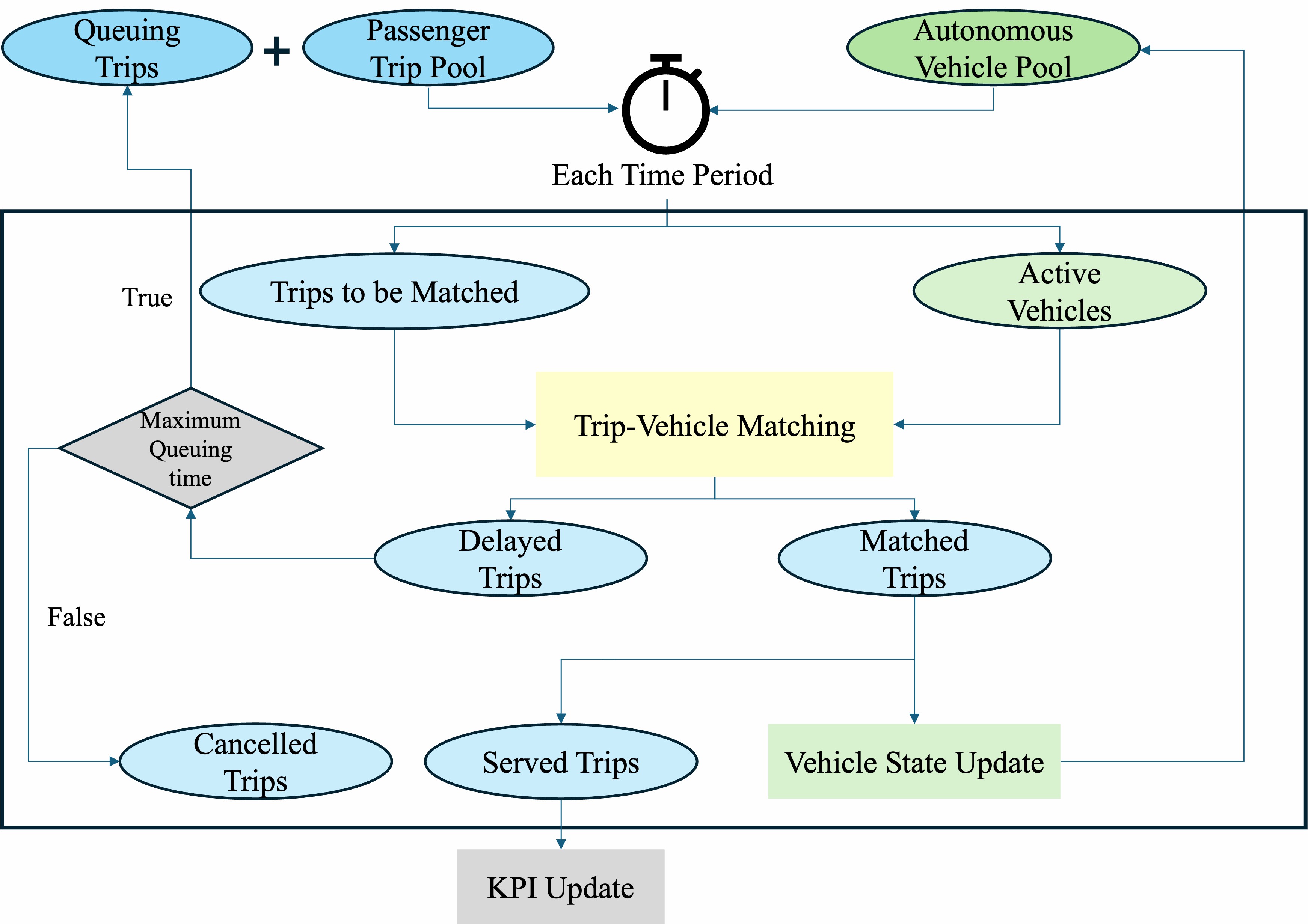}
\caption{Process of the AMoD simulation framework.}
\label{fig:simulation_flow}
\end{figure}

Figure~\ref{fig:simulation_flow} shows the simulation framework used to evaluate the performance of AMoD services based on an urban road network. We consider two types of entities: autonomous vehicles and passenger trips. The trip pool comprises all trips occurring over a single day, according to the historical dataset or generated from a probability distribution. The autonomous vehicle pool of AMoD contains all vehicles operating on a day. For each period, the set of trips to be matched is composed of the newly occurred trips together with the unmatched trips from the previous period, and these trips are matched against the active vehicles according to the trip-vehicle matching mechanism. If a trip is successfully matched in the current period, it is then classified as a matched trip, waiting for a vehicle to pick it up, and finally regarded as a successfully served trip. If a trip is not matched, it is set as a delayed trip. A delayed trip is cancelled if it has queued for longer than the maximum passenger queuing time; otherwise, it is regarded as a backlogged trip and matched in the next period. The vehicle state is updated according to the actual situation (e.g., picking-up, carrying) at the end of the period. Finally, the KPI values are calculated and updated.

Entity attributes of autonomous vehicles and trips are shown in Table~\ref{tab:agent_attributes} and the parameters of the AMoD system that define platform constraints are summarized in Table~\ref{tab:system_parameters}. These parameters are either fixed or used as control inputs in simulation experiments. We assume that passengers will cancel their trips if the queuing delay $\delta^{queue}_u$ exceeds the maximum threshold \( t_{\text{queue}}^{\max} \). A trip-vehicle pair $(u,v)$ is considered feasible if it satisfies the following criteria: (1) Queuing delay \( \delta^{queue}_u \leq t_{\text{queue}}^{\max} \) and (2) Pickup delay \( \delta^{pu}_u \leq t_{\text{pu}}^{\max}\).

\begin{table}[!htbp]
\centering
\caption{Attributes of autonomous vehicles and trips}
\label{tab:agent_attributes}
\resizebox{\textwidth}{!}{
\begin{tabular}{>{\centering\arraybackslash}m{2.5cm} l c m{9cm}}
\toprule
\textbf{Entity} & \textbf{Attribute} & \textbf{Symbol} & \textbf{Description} \\
\midrule
\multirow{3}{*}{Vehicle}
  & Identification & $v$& Unique ID of the vehicle. \\
  & State & $s_t^v$& State of vehicle $v$ at period $t$: \{0: Resting, 1: Vacant, 2: Picking-up, 3: Carrying\}.\\
  & Location (node) & $n_t^v$& Location (current node) of vehicle $v$ on the road network before period $t$.\\
\midrule
\multirow{7}{*}{Trip}
  & Identification & $u$& Unique ID of the trip.\\
  & State & $s_t^u$& State of trip $u$ at period $t$: \{0: In queue, 1: Waiting, 2: In transit\}.\\
  & Request time & $t^{u}_{req}$& Time at which trip $u$ enters the matching system.\\
  & Origin & $o_u$ & Origin node of the trip $u$.\\
  & Destination & $d_u$ & Destination node of the trip $u$.\\
  & Queuing delay & $\delta^{queue}_u$ & Time spent in the queue of trip $u$.\\
  & Pickup delay & $\delta^{pu}_u$ & Time taken for trip $u$ from being matched to being picked up.\\
\bottomrule
\end{tabular}
}
\end{table}

\begin{table}[!htbp]
\centering
\caption{Parameters of the AMoD system}
\label{tab:system_parameters}
\resizebox{\textwidth}{!}{
\begin{tabular}{l c p{10cm}}
\toprule
\textbf{Parameter Name} & \textbf{Symbol} & \textbf{Description} \\
\midrule
Fleet size & $F_0$ & Total number of autonomous vehicles in the system. \\
Vehicles by state & $N^{\text{rest}}, N^{\text{vac}}, N^{\text{pu}}, N^{\text{car}}$ & Number of vehicles in resting, vacant, picking-up, and passenger-carrying states. \\
Trip set & \(O\)& Set of all trips occurring on a single day.\\
Matching interval & $period$& Time interval between successive batches. \\
Max pickup time & $t_{\text{pu}}^{\max}$ & Maximum time limit for a vehicle to pick up a trip.\\
Max queue time & $t_{\text{queue}}^{\max}$ & Maximum time a passenger is willing to wait in the queue for matching.\\
\bottomrule
\end{tabular}
}
\end{table}

We use the same historical trips for both AMoD and street-hailing services. The urban transportation network is represented as a directed graph \( G = (N, E) \), where \( N \) denotes the set of road intersection nodes. The origin and destination of each trip are matched to the nearest road intersection node. Meanwhile, \( E \subseteq N \times N \) represents the set of directed road edges connecting nodes in \( N \). Each edge is associated with a travel cost (e.g., travel time, length or energy consumption), typically derived from real road and travel data, or directly from historical speed or time data. During the baseline simulation, vehicles in resting or vacant states remain at the destination nodes of their recently completed trips until they are activated and assigned to a new trip. This represents a service under a no-repositioning (or drop-and-remain) setting, which is used to isolate the effect of centralized trip–vehicle matching of AMoD services before introducing additional repositioning policies.

\begin{figure}[!htbp]
\centering
\includegraphics[width=1.0\textwidth]{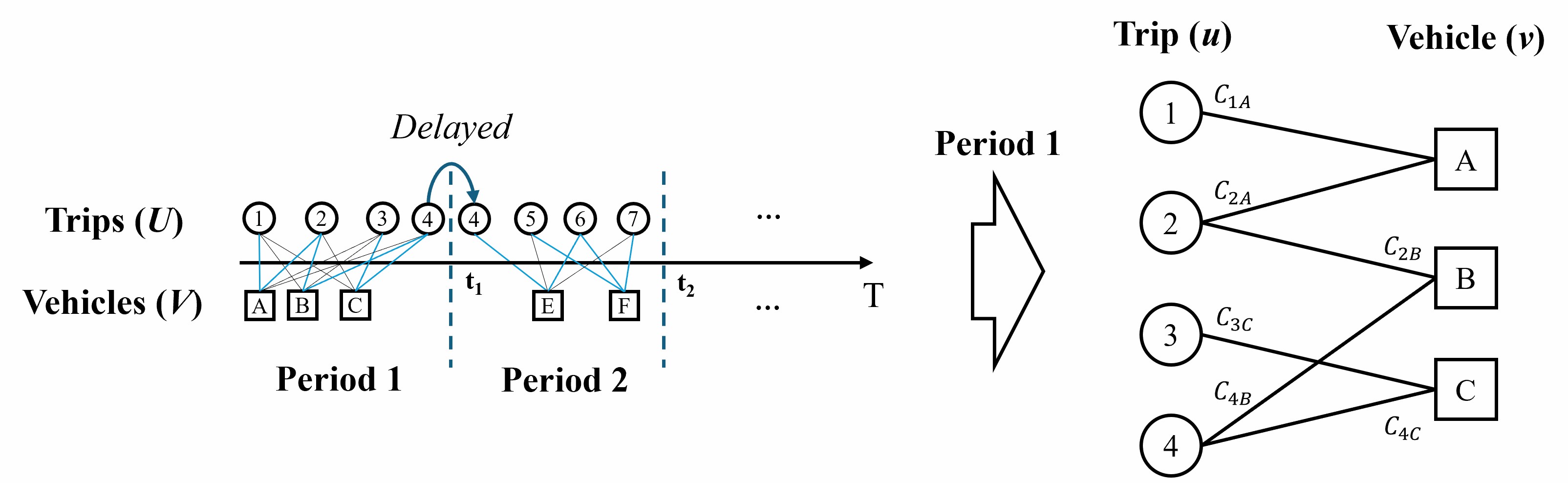}
\caption{Illustration of the batch matching process. Trips and vehicles are collected within a period. The candidate pairs make up a bipartite graph representing an assignment problem.}
\label{fig:matching_flow}
\end{figure}

\subsection{Trip-Vehicle matching}
We use a batch matching method in a non-anticipatory rolling-horizon setting, which aggregates trips within each one-minute period and matches trips and vehicles in each period sequentially without information about the future. Through this, we transform the trip-vehicle matching problem into a bipartite matching and solve it with an algorithm based on graph theory. The objective of batch matching is to minimize the total pickup delay. Let \( U \) denote the set of trips to be matched in the current period and \( V \) denote the set of available autonomous vehicles. Figure~\ref{fig:matching_flow} illustrates the batch matching process. In each period, trips accumulated within the current batch and available vehicles are used to construct a candidate trip--vehicle assignment problem at the end of the batch. All candidate trip--vehicle pairs are represented in a cost matrix \(C\), where each entry \(C_{uv}\) denotes the assignment cost of assigning vehicle \(v \in V\) to trip \(u \in U\). For feasible pairs that satisfy the time constraints, the assignment cost is set to the pickup delay \(\delta_{uv}^{\mathrm{pu}}\). For infeasible pairs, the assignment cost is set to a sufficiently large penalty value. In this way, the trip--vehicle matching problem is formulated as a bipartite assignment problem over the trip set \(U\) and vehicle set \(V\). The assignment cost is defined as:

\[
C_{uv} =
\begin{cases}
\delta_{uv}^{\mathrm{pu}},
& \text{if } \delta_u^{\mathrm{queue}} \leq t_{\mathrm{queue}}^{\max}
\text{ and } \delta_{uv}^{\mathrm{pu}} \leq t_{\mathrm{pu}}^{\max},\\
\Lambda,
& \text{otherwise},
\end{cases}
\]

\noindent where \(\delta_{uv}^{\mathrm{pu}}\) is the pickup time from vehicle \(v\) to trip \(u\), \(\delta_u^{\mathrm{queue}}\) is the accumulated queuing time of trip \(u\), and \(\Lambda\) is set sufficiently larger than the maximum feasible pickup delay so that the assignment solver prioritizes all feasible matches before selecting infeasible assignments.

The assignment problem is formulated as:
\[
\min_{x} \sum_{u\in U}\sum_{v\in V} C_{uv} x_{uv}
\]
subject to:
\[
\sum_{v \in V} x_{uv} \leq 1, \qquad \forall u \in U,
\]
\[
\sum_{u \in U} x_{uv} \leq 1, \qquad \forall v \in V,
\]
\[
\sum_{u\in U}\sum_{v\in V} x_{uv} = \min(|U|,|V|),
\]
\[
x_{uv} \in \{0,1\}, \qquad \forall u\in U,\ v\in V.
\]

\noindent The first two constraints ensure that each trip and each vehicle are assigned at most once. The third constraint requires the solver to return the maximum possible number of tentative assignments given the numbers of trips and available vehicles. After solving the assignment problem, any tentative assignment with \(C_{uv}=\Lambda\) is treated as infeasible. The corresponding trip remains unmatched and is reconsidered in the next matching period unless its accumulated queuing time exceeds \(t_{\mathrm{queue}}^{\max}\). We employed the Jonker–Volgenant (JV) algorithm \citep{JonkerVolgenant1986_hungarian_improvement} to solve the resulting linear assignment problem by iteratively refining a matching between trips and vehicles.

\subsection{Experimental data and settings} \label{sec:data_settings}

All experiments are carried out based on taxi operations data from street-hailing services in Chengdu, China and its road network within the Ring Expressway. The network covers the operating area of cruising taxis \citep{Guo2019_Chengdu_GPS}, which is modelled as a directed graph \(G=(N,E)\) with \(|N| = 1\,902\) nodes and \(|E| = 5\,943\) edges. The primary analysis uses taxi operations data collected in April 2015 (described in Table~\ref{tab:trip_data}). To examine the temporal robustness of the findings, we additionally conduct an out-of-period analysis using taxi GPS trajectory data from the same road network collected in April 2017. This data is processed using the same method as used in the experiments based on 2015 data. We repeat the comparison between street-hailing and simulated AMoD services and examine the effects of different operational strategies. The complete experimental settings and results are reported in~\ref{app:robustness_2017}.

\begin{table}[!htbp]
\centering
\caption{Fields and descriptions of street-hailing trip data}
\label{tab:trip_data}
\begin{tabular}{ll}
\toprule
Field & Description \\
\midrule
Vehicle ID & Unique taxi identifier \\
Pickup location & Latitude and longitude of the trip's start location\\
Drop-off location & Latitude and longitude of the trip's end location\\
Pickup time & Start time of a trip\\
Drop-off time & End time of a trip\\
Deadheading mileage& Distance without passenger between consecutive trips\\
\bottomrule
\end{tabular}
\end{table}

\paragraph{Data cleaning of street-hailing operations data}
First, records with missing values in any key field listed in Table~\ref{tab:trip_data} were removed. Next, duplicate order records were identified based on the same key fields and removed. Then, we deleted records with trip start times outside the operational day, records with unreasonable trip mileage, and records whose latitude or longitude coordinates were more than 500 meters away from the road network links. These abnormal records were generally attributable to logging errors in the taxi terminals. For example, the cleaned data of April 1, 2015 contains more than 405,000 valid trips from more than 11,200 active vehicles.

\paragraph{Passenger arrival time and simulation duration}
Because the historical street-hailing taxi data record the observed pickup time but not the time a passenger starts waiting, we use the observed pickup time as the entry timestamp in all simulation experiments.

Based on this setting, for each simulation day, trip arrivals are generated only within the 24-hour observation window. Dispatching is allowed to continue for an additional period equal to $t_{\text{queue}}^{\max}$ after the end of the observation window. The simulation then continues only to allow vehicles already serving passengers to complete their trips.

\paragraph{Travel speed and time}
Vehicles follow shortest-travel-time routes based on time-dependent travel times estimated based on historical trips. The corresponding road network speed is estimated using a gradient descent algorithm~\citep{ruder2016overview} based on the historical taxi trip data. Then we compute travel time of each road segment in each hour by dividing the length of each road segment by the speed in the corresponding time period. Finally, the Floyd–Warshall algorithm~\citep{floyd1962algorithm} is employed to compute the all-pairs shortest travel time matrices across the entire network over 24 hours. All algorithms and processing details related to the estimation of driving speed and time can be found in \ref{app:prediction}. In the data preprocessing stage, certain trips are excluded to improve the reliability of the analysis. Specifically, trips with loops (where the origin and destination coincide) are rejected, together with trips that are unrealistically short or excessively long. Route computation is then performed by assigning an overall initial speed $v_{\text{int}}=40$ km/h to all road segments, following the procedure described above. Roughly 95\% of the original trips remain available for further estimation. We further set the upper and lower bounds of travel speeds for all road segments to $[5, 80]$ km/h based on historical speed data on the road network of Chengdu. We obtain reasonable travel times according to the following setting: Average road network speeds are highest during off-peak hours (midnight to early morning) with around 38 km/h, gradually decrease during the morning commute, and reach the lowest values during the morning (around 8:00) and evening (17:00–19:00) peak periods with around 28 km/h.

\paragraph{Dynamic fleet-size adjustment and vehicle initialization}
To determine the active fleet size of the simulated AMoD system, we take each hour as an update period. The hourly active fleet size of the simulated AMoD system is set to the number of street-hailing vehicles with at least one valid service trip starting within that period (hour). We assume that the activity status of simulated AMoD vehicles is updated at the beginning of each hour. Let \( N_{h}^{\text{active}} \) denote the number of active vehicles of the simulated AMoD services at the start of hour \( h \), defined as:
\begin{equation}
N_{h}^{\text{active}} = N_{h}^{\text{vac}} + N_{h}^{\text{pu}} + N_{h}^{\text{car}}
\end{equation}
where \( N_{h}^{\text{vac}} \), \( N_{h}^{\text{pu}} \), and \( N_{h}^{\text{car}} \) represent the numbers of vacant, picking-up, and passenger-carrying vehicles in hour \( h \), respectively. These states are updated according to the following fleet activation rule. For hour \( h+1\), vehicle states are updated according to the following rules: (i) If \( N_{h}^{\text{active}} < N^{\text{active}}_{h+1} \), a subset of resting vehicles (\( s_v = 0 \)) is randomly activated and assigned an initial vacant state (\( s_v = 1 \)) until the active fleet size reaches \(N^{\text{active}}_{h+1}\). (ii) If \( N_{h}^{\text{active}} > N^{\text{active}}_{h+1} \), vehicles are deactivated until the active fleet size is reduced to \( N^{\text{active}}_{h+1} \). Vacant vehicles are deactivated first, followed by picking-up vehicles and passenger-carrying vehicles after they complete their current trips. We refer to this rule as the baseline state-priority rule: when the active fleet size decreases, vacant vehicles are selected first, followed by picking-up and passenger-carrying vehicles; vehicles with the same state are selected randomly. This rule is used as the baseline fleet activation rule. For the initialization of vehicle positions, they are initialized by sampling with replacement from the set of distinct historical origin nodes, with equal probability assigned to each unique origin node.

All simulations and algorithm implementations are conducted using Python 3.11 and executed on a workstation with an AMD Ryzen 9 5900X CPU (3.7 GHz) and 32 GB RAM. Due to the uncertainties introduced by random initial vehicle locations each day and stochastic vehicle state transitions every hour, we performed 30 independent Monte Carlo simulations for each $t_{\text{pu}}^{\max}$ and $t_{\text{queue}}^{\max}$ to enhance the robustness of the results.

\subsection{Key performance indicators}
We adopt three key performance indicators, average passenger waiting time (APWT), average deadheading mileage (ADM), and average deadheading energy consumption (ADEC), to evaluate the performance of both AMoD and street-hailing services. APWT is closely related to the objective function of AMoD dispatching in this study and is also a key indicator for evaluating passenger satisfaction and system service capability~\citep{DengChenShaoTangPiGupta2022JMSE}. ADM and ADEC are critical indicators for evaluating the impacts of different mobility service modes on the transportation system and environment: ADM measures the average vehicle mileage travelled without passengers per served trip, which captures the deadheading part of vehicle operations~\citep{KontouGarikapatiHou2020TRC}; ADEC focuses on energy consumption to evaluate environmental performance~\citep{WardMichalekAzevedoSamarasFerreira2019TRC}. Lower ADM and ADEC indicate less traffic and lower energy use associated with deadheading trips. These three indicators capture the main passenger-side and vehicle-operation dimensions considered in the comparative analysis of this study.

The calculations of the three key performance indicators for AMoD services are introduced as follows.

(i) APWT is defined as the average passenger waiting time for trips served. The waiting time for a trip $u$ is the sum of the queuing time $\delta_{u}^{\text{queue}}$ before this trip is accepted and the pickup time $\delta_{u}^{\text{pu}}$ taken for a vehicle to arrive at the origin node of this trip, which is formulated as:
    \begin{equation}
    \text{APWT} = \frac{1}{|U_{\text{served}}|} \sum_{u \in U_{\text{served}}} \left( \delta_{u}^{\text{queue}} + \delta_{u}^{\text{pu}} \right),
    \end{equation}
where \(U_{\text{served}}\) refers to the set of served trips.

(ii) ADM refers to the average mileage of deadheading segments of each trip. A deadheading segment is defined as a segment travelled by a vehicle from completing one trip to picking up the passenger for the next trip, including pickup trips and, where applicable, proactive repositioning trips. Let \mbox{\(R\)} denote the set of proactive repositioning movements conducted during the simulation. The ADM is calculated as

\begin{equation}
\text{\mbox{\(\displaystyle \text{ADM}
=
\frac{1}{|U_{\text{served}}|}
\left(
\sum_{u \in U_{\text{served}}} md^{\text{pu}}_u
+
\sum_{r \in R} md^{\text{rep}}_r
\right).
\)}}\end{equation}

\noindent where \mbox{\(md^{\text{pu}}_u\)} denotes the pickup mileage associated with served trip \mbox{\(u\)}, and \mbox{\(md^{\text{rep}}_r\)} denotes the mileage of proactive repositioning movement \mbox{\(r\)}. When vehicle repositioning is not considered, \mbox{\(R=\varnothing\)}, and ADM includes only pickup mileage.

(iii) ADEC is the average energy consumption of deadheading segments. We use the energy consumption formula of electric vehicles to calculate the ADEC of all deadheading segments. Let \( td \) represent the travel time (minutes). The energy consumption of each deadheading segment DEC is estimated using a macro-level model referring to the work of~\cite{DeCauwer2015}:
    \begin{equation}
    \label{eq:dec}
    DEC(md, td) = B_1 \cdot md + B_2 \cdot \left( \frac{60 \cdot md}{td} \right)^2 \cdot md,
    \end{equation}
where $B_1$ and $B_2$ are model coefficients adopted from~\cite{DeCauwer2015}, set as 0.132 and $5 \times 10^{-6}$, respectively. The ADEC (in kWh) is then formulated as:

    \begin{equation}
    \label{eq:adec}
    \mathrm{ADEC}
    =
    \frac{1}{\left|U_{\mathrm{served}}\right|}
    \left[
    \sum_{u \in U_{\mathrm{served}}}
    DEC\left(md_{u}^{\mathrm{pu}}, td_{u}^{\mathrm{pu}}\right)
    +
    \sum_{r \in R}
    DEC\left(md_{r}^{\mathrm{rep}}, td_{r}^{\mathrm{rep}}\right)
    \right].
    \end{equation}

The three key performance indicators of street-hailing services are obtained as follows. 

APWT of street-hailing services is not directly observable in historical trajectory data. In conventional street-hailing systems, passengers appear along the roadside without registering a digital request. Consequently, the time at which a passenger begins waiting is not recorded in either the taxi GPS trajectories or the trip records. Thus, we adopt an APWT benchmark for street-hailing based on a survey conducted in Chengdu in April 2015~\citep{RolandBerger2015}, shown in Table~\ref{tab:street_hailing_waiting_distribution}. Because the survey reports the distribution of waiting time ranges, we estimate the street-hailing APWT using the lower and upper endpoints of each range. After weighting by the share of each range, the two estimates are 8.5 and 13.5 minutes, respectively. We adopt a simplified cruising-taxi simulation to examine whether the survey-based benchmark is reasonable. The simulation model and experiment are detailed in~\ref{app:street-hailing_simulation}. The experiment shows that, for the service rates between 93.6\% and 96.6\%, the resulting APWT varies from 3.4 to 13.7 minutes, a range that encompasses our survey-based APWT. 

Then, the ADM and ADEC are strictly calculated based on the historical dataset. Specifically, ADM is calculated using historical service trips of all street-hailing taxis in Chengdu during the study period (see \ref{app:ADM_street}); ADEC is then calculated using Equations \ref{eq:dec} and \ref{eq:adec} based on the deadheading mileage and travel time of street-hailing services. Numerically, the ADM and ADEC of street-hailing services are 2.50 km and 0.33 kWh, respectively.

\begin{table}[!htbp]
\centering
\begin{threeparttable}
\caption{Distribution of street-hailing passenger waiting time in a survey and the estimated APWT benchmark~\citep{RolandBerger2015}}
\label{tab:street_hailing_waiting_distribution}

\begin{tabular}{ccc}
\toprule
Passenger waiting time (min)
& Share
& Weighted APWT range (min) \\
\midrule
$[0,5]$   & 10\% & \multirow{4}{*}{$[8.5,13.5]$} \\
$[5,10]$  & 31\% & \\
$[10,15]$ & 39\% & \\
$[15,20]$ & 20\% & \\
\bottomrule
\end{tabular}

\begin{tablenotes}[flushleft]
\footnotesize
\item \textit{Note:} The interval $[8.5,13.5]$ min is used as the
2015 survey-based street-hailing APWT benchmark, following
\citet{RolandBerger2015}.
\end{tablenotes}
\end{threeparttable}
\end{table}

In addition, we assume that both street-hailing and AMoD operate vehicles with identical technical specifications and therefore use the same parameter values for both service modes. This ensures that differences in ADEC reflect operational differences, such as variations in deadheading mileage and time, rather than differences in vehicle technical characteristics.

 \section{Results}\label{sec:results}
\begin{figure}[!htbp]
    \centering
    \begin{subfigure}[b]{0.45\textwidth}
        \centering
        \includegraphics[width=\textwidth]{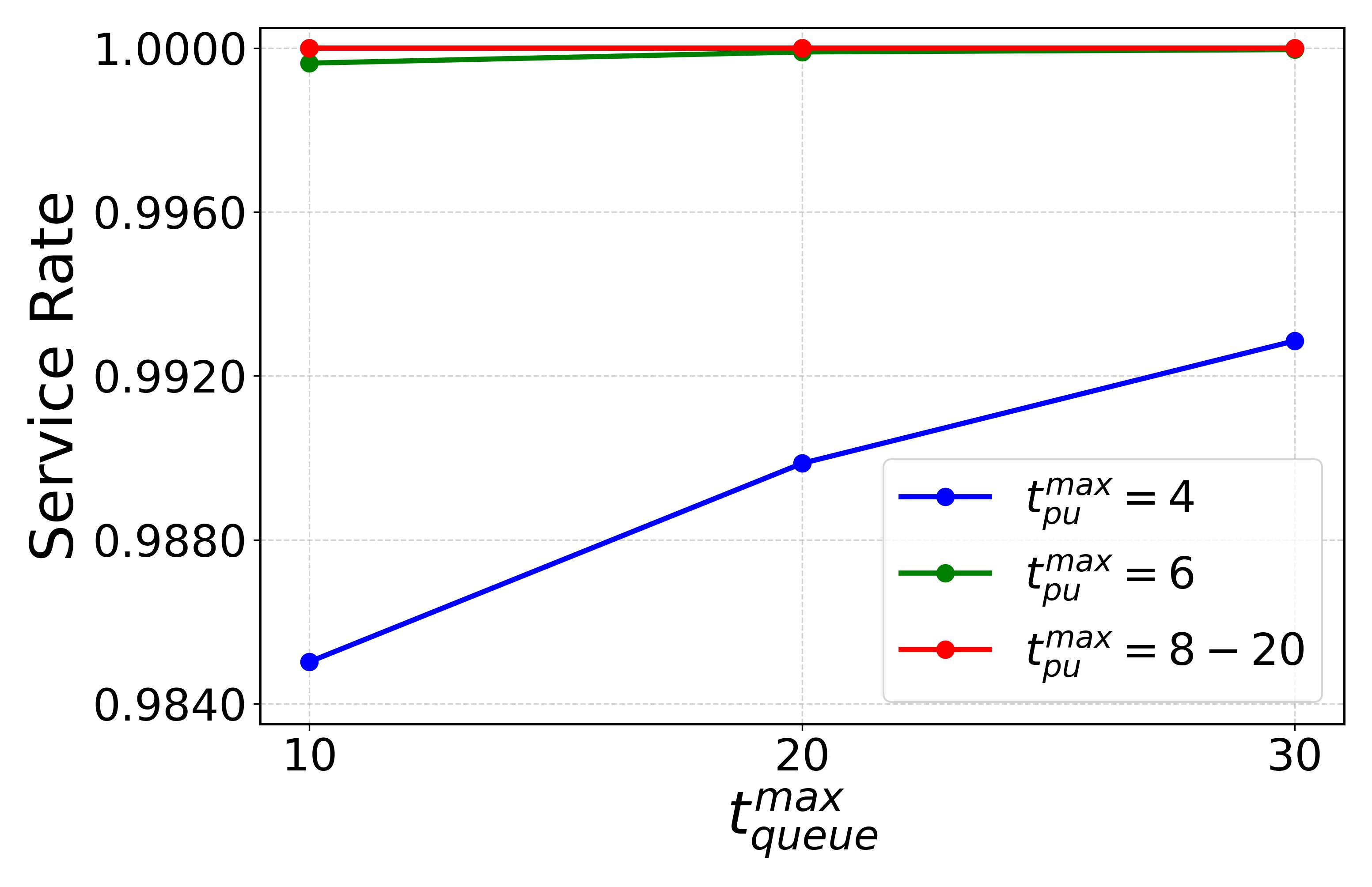}
        \caption{Service rate}
        \label{fig:patience_response}
    \end{subfigure}
    \hfill
    \begin{subfigure}[b]{0.45\textwidth}
        \centering
        \includegraphics[width=\textwidth]{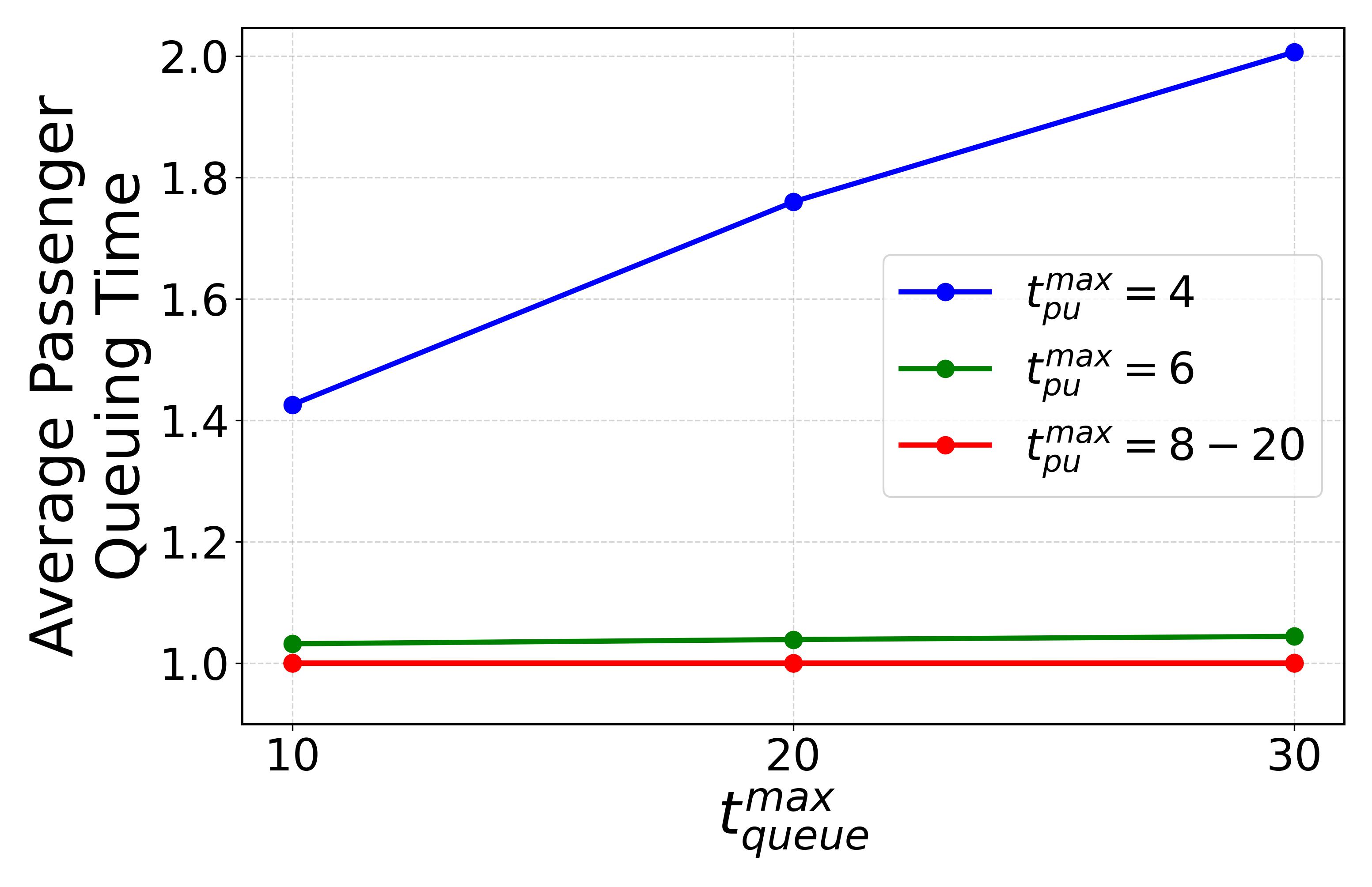}
        \caption{Average passenger queuing time}
        \label{fig:patience_delay}
    \end{subfigure}
    \caption{Service rate and average passenger queuing time of AMoD services influenced by $t_{\text{queue}}^{\max}$ and $t_{\text{pu}}^{\max}$.}
    \label{fig:patience}
\end{figure}

In our experimental setting, all trips in the street-hailing service dataset are served in the simulated AMoD system within reasonable passenger queuing times. Thus, we first explore parameter settings under which all trips can be served. Figure~\ref{fig:patience} shows the changes in service rates and average passenger queuing times over different maximum queuing times $t_{\text{queue}}^{\max}$ and maximum pickup times $t_{\text{pu}}^{\max}$. The results illustrate that when \( t_{\text{pu}}^{\max} < 8 \), not all trips can be served, with longer queuing and pickup time for the served trips. However, when \( t_{\text{pu}}^{\max} \geq 8 \), all trips are successfully served, and the average passenger queuing time of the AMoD system drops rapidly to a very low level (\(\approx 1\) minute). We therefore set $t_{\text{pu}}^{\max}\geq8$ minutes to guarantee a service rate of 100\%. Following this setting, the average passenger queuing time remains nearly constant for $t_{\text{queue}}^{\max}\geq 10$ and we thus simply set $t_{\text{queue}}^{\max}=20$ minutes. We also report the results of sensitivity tests of parameters and key assumptions in~\ref{app:sensitivity_analysis}.

\subsection{Comparison between AMoD and street-hailing services} \label{sec:comparison_rh_cruise}
We compare three key performance indicators of street-hailing services with those of the simulated AMoD service under different settings of maximum passenger pickup time \( t_{\text{pu}}^{\max}\) representing the service radius of vehicles.

\begin{figure}[!htbp]
\centering
\includegraphics[width=1.0\textwidth]{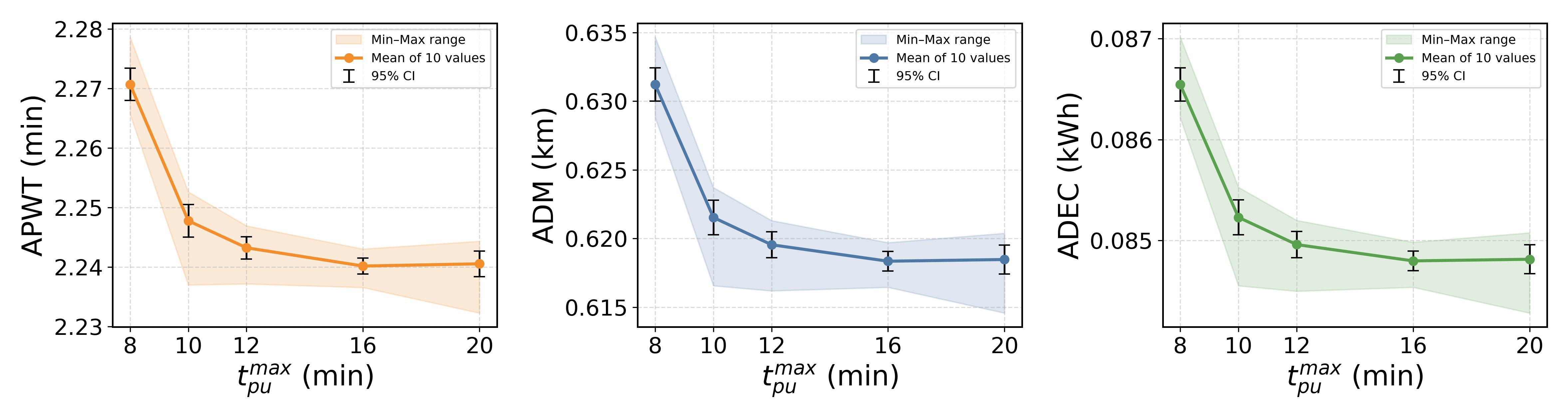}
\caption{Changes of AMoD performance over the maximum pickup time $t_{\text{pu}}^{\max}$. Average passenger waiting time (APWT), average deadheading mileage (ADM), and average deadheading energy consumption (ADEC) generally decrease as $t_{\text{pu}}^{\max}$ increases, while these improvements gradually diminish. The light-coloured bands represent the range of results from all 30 Monte Carlo replications. Each data point represents the mean value derived from the corresponding experiments, with vertical error bars indicating the 95\% confidence interval.}
\label{fig:ride_hailing_errorbands}
\end{figure}

First, we examine the effect of \( t_{\text{pu}}^{\max}\) on the performance of AMoD services for all service trips in the study period. Results in Figure~\ref{fig:ride_hailing_errorbands} show that \( t_{\text{pu}}^{\max}\) is an important operational parameter in centrally dispatched AMoD services and has a significant impact on the overall system efficiency. Increasing \(t_{\text{pu}}^{\max}\) leads to significant improvements in all three performance indicators. Specifically, APWT decreases from 2.27~min at \(t_{\text{pu}}^{\max}=8\)~min to 2.24~min at \(t_{\text{pu}}^{\max}=20\)~min, ADM decreases from 0.63~km to 0.62~km, and ADEC decreases from 0.087~kWh to 0.085~kWh. Notably, although these correlations between \(t_{\text{pu}}^{\max}\) and three key performance indicators are statistically significant and consistent when \(t_{\text{pu}}^{\max}\leq16\) minutes, further relaxation of this constraint yields only marginal improvements in APWT, ADM, and ADEC, indicating diminishing returns, as shown in Figure~\ref{fig:ride_hailing_errorbands}.

\begin{figure}[!htbp]
\centering
\includegraphics[width=0.8\textwidth]{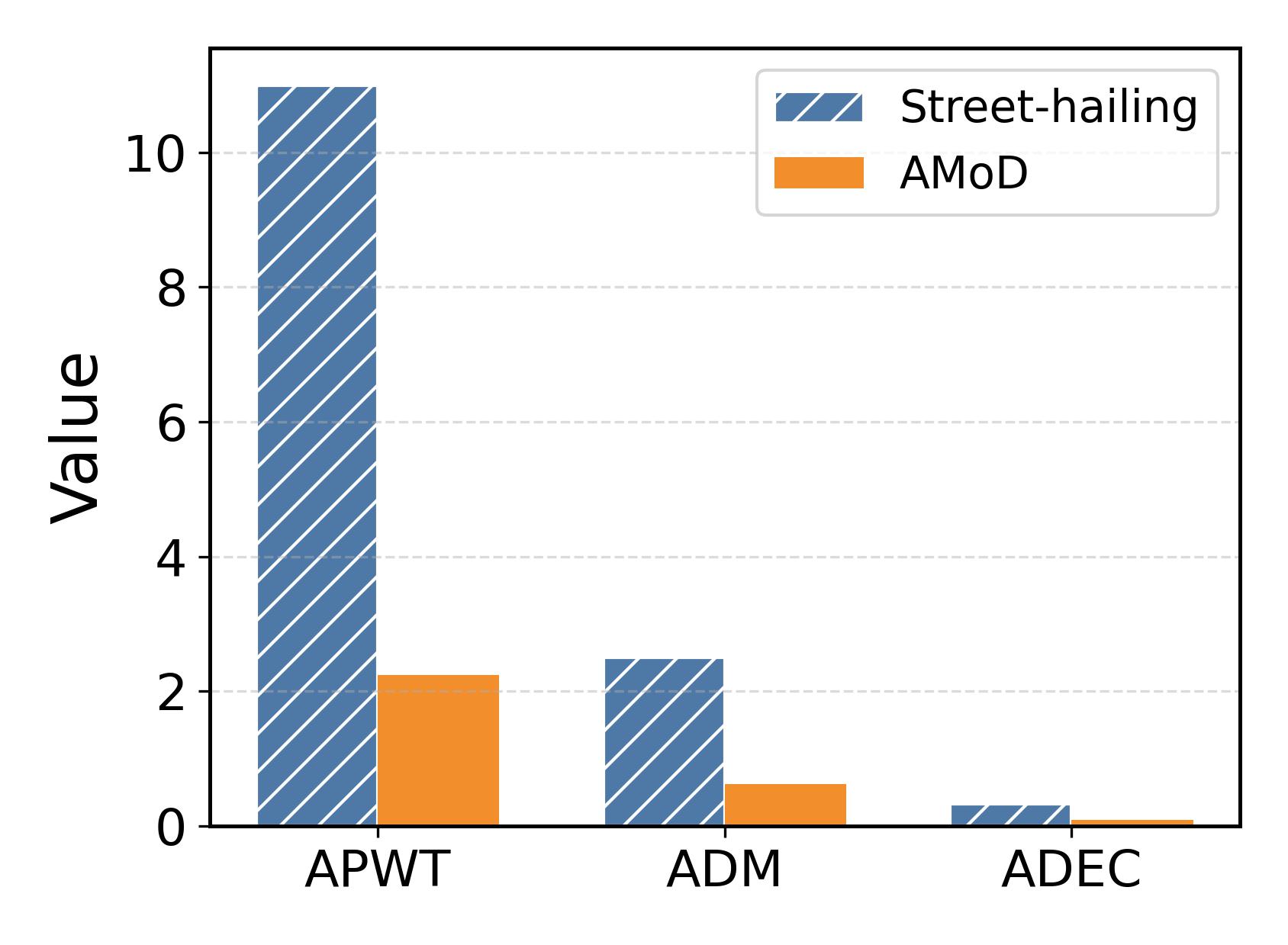}
\caption{Comparison between simulated AMoD services (at $t_{\text{pu}}^{\max} \in [8,20]$ minutes) and street-hailing services. The variation in the AMoD indicators over \mbox{$t_{\mathrm{pu}}^{\max}\in[8,20]$} is small relative to the differences between the two services.}
\label{fig:ride_vs_street}
\end{figure}

Figure~\ref{fig:ride_vs_street} presents the numerical differences between the simulated AMoD service at $t_{\text{pu}}^{\max} \in [8,20]$ and the corresponding street-hailing benchmarks. Under the simulation settings, the APWT of the simulated AMoD system is approximately 73.3\% -- 83.4\% lower than the survey-based street-hailing APWT benchmark when \( t_{\text{pu}}^{\max} \geq 8 \) minutes. The lower end of this percentage range corresponds to the street-hailing APWT of 8.5 minutes, whereas the upper end corresponds to the 13.5-minute scenario. The resulting range includes the 76.5\% difference reported by \cite{Brown2021Transportation}. The difference may be attributed to variations in driver behaviour, data sources, and urban environments.

\begin{figure}[!htbp]
\centering
\includegraphics[width=0.9\textwidth]{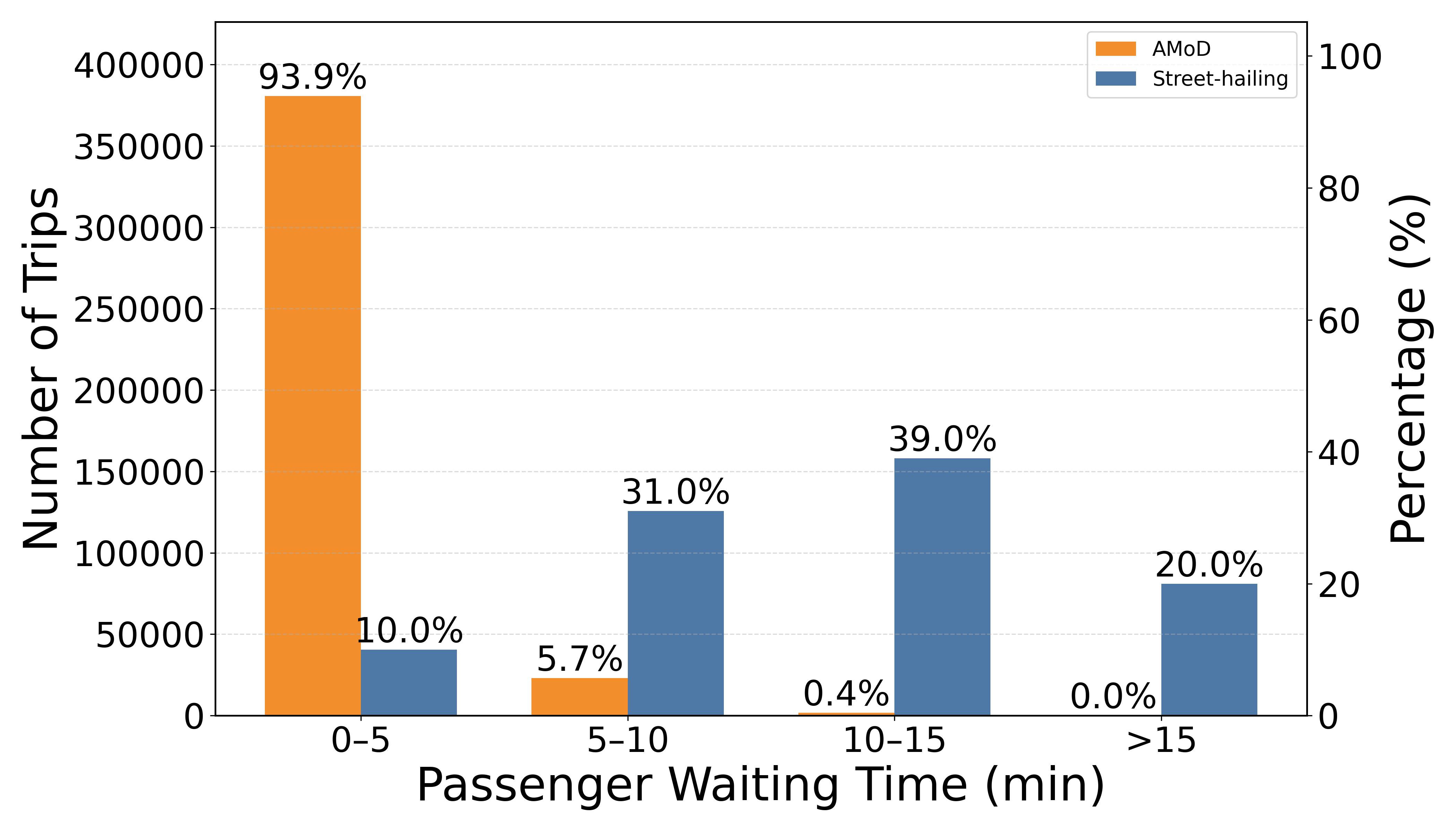}
\caption{Passenger waiting-time distributions for the simulated AMoD service and the survey-based street-hailing benchmark.}
\label{fig:waiting_time_distribution}
\end{figure}

Besides, it shows that over 99.9\% of trips can be matched by the platform within 1 minute, indicating high operational efficiency of the simulated AMoD service with batch matching under the simulation settings. Figure~\ref{fig:waiting_time_distribution} shows that the APWT of street-hailing services according to the survey details is more dispersed and a large proportion of trips require longer waiting times, compared with that of simulated AMoD services. A majority of AMoD trips are responded to and picked up within 5 minutes, while 39\% of street-hailing trips need a 10–15 minute waiting time and 20\% take longer than 15 minutes. This indicates that the simulated AMoD service tends to require less time to serve a trip request than the survey-based street-hailing benchmark, suggesting potential operational gains in service efficiency and stability associated with centralized dispatch. The ADM and ADEC values of the simulated AMoD service are approximately 75.0\% and 74.0\% lower, respectively, than the street-hailing benchmarks. These improvements lie between the values reported by~\cite{Ruch2021JAT} (58\%) and~\cite{Zhan2016TITS} (87\%). This variation may be related to differences in modelling assumptions. For example,~\cite{Ruch2021JAT} approximated street-hailing inter-trip distances using shortest-path estimates rather than real cruising trajectories, which likely understates the actual ADM of street-hailing vehicles and thus reduces the measured improvement. In contrast,~\cite{Zhan2016TITS} assumed a perfect-information, system-wide optimal dispatching scheme, which may overestimate the achievable performance of AMoD systems compared to a rolling-horizon framework. These methodological differences can lead to variations in the reported improvement magnitudes.

We also test whether the performance differences between AMoD and street-hailing services remain similar when the matching interval, maximum pickup time, maximum queuing time, fleet activation rule, and travel-time estimates are changed. The full results are reported in~\ref{app:sensitivity_analysis}. The results show the same overall pattern across all simulation settings.

\subsection{Spatiotemporal analysis of AMoD services} \label{sec:sp_analysis}
We further conduct a spatiotemporal analysis with \( t_{\text{pu}}^{\max} =16\) and \( t_{\text{queue}}^{\max}=20\) to understand how the AMoD performance varies throughout the day and across city areas. We set \( t_{\text{pu}}^{\max} \) to 16 minutes because further increasing \( t_{\text{pu}}^{\max} \) does not bring obvious changes in the three key performance indicators (as shown in Figure~\ref{fig:ride_hailing_errorbands}).

\begin{figure}[!htbp]
    \centering
    \includegraphics[width=\textwidth]{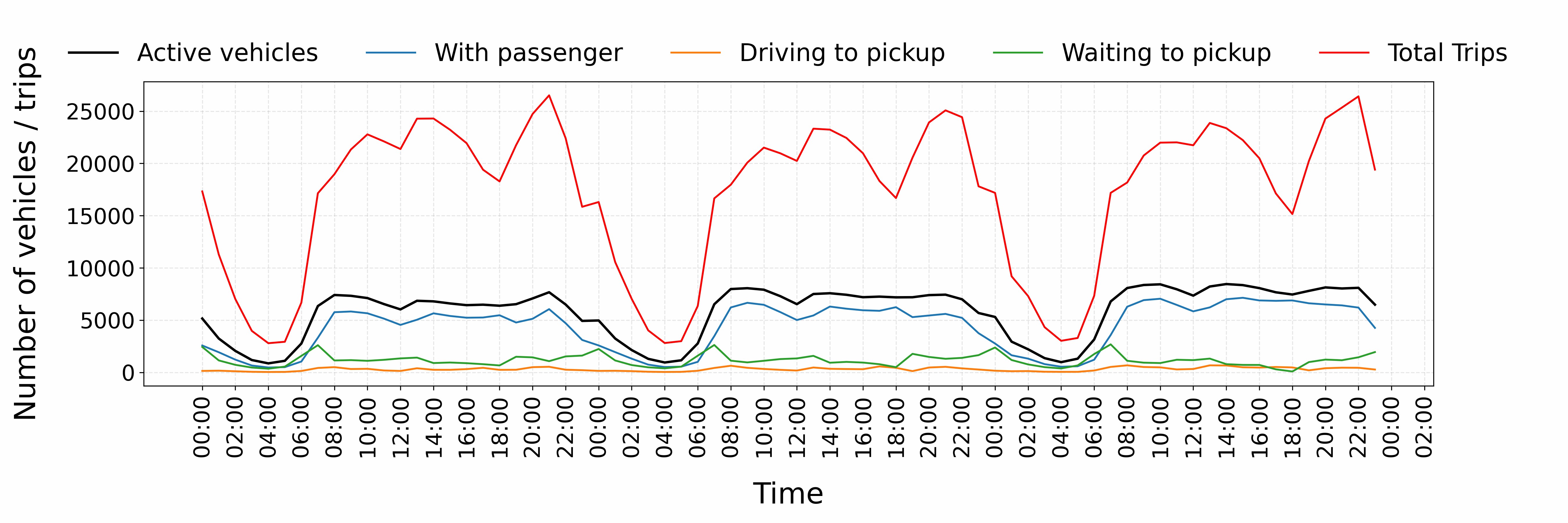}
    \caption{Changes of the number of vehicles in different states over time.}
    \label{fig:vehicle_status}
\end{figure}

\paragraph{Temporal analysis} Figure~\ref{fig:vehicle_status} shows the hourly variation of the AMoD fleet size in different states over a 72-hour operational period (3 simulation days). The utilization of the vehicles remains relatively stable throughout the day, with most of the active vehicles (the black line) being with passengers (the blue line), despite significant fluctuations in the total number of trips (the red line).

\begin{figure}[!htbp]
    \centering
    \includegraphics[width=0.95\linewidth]{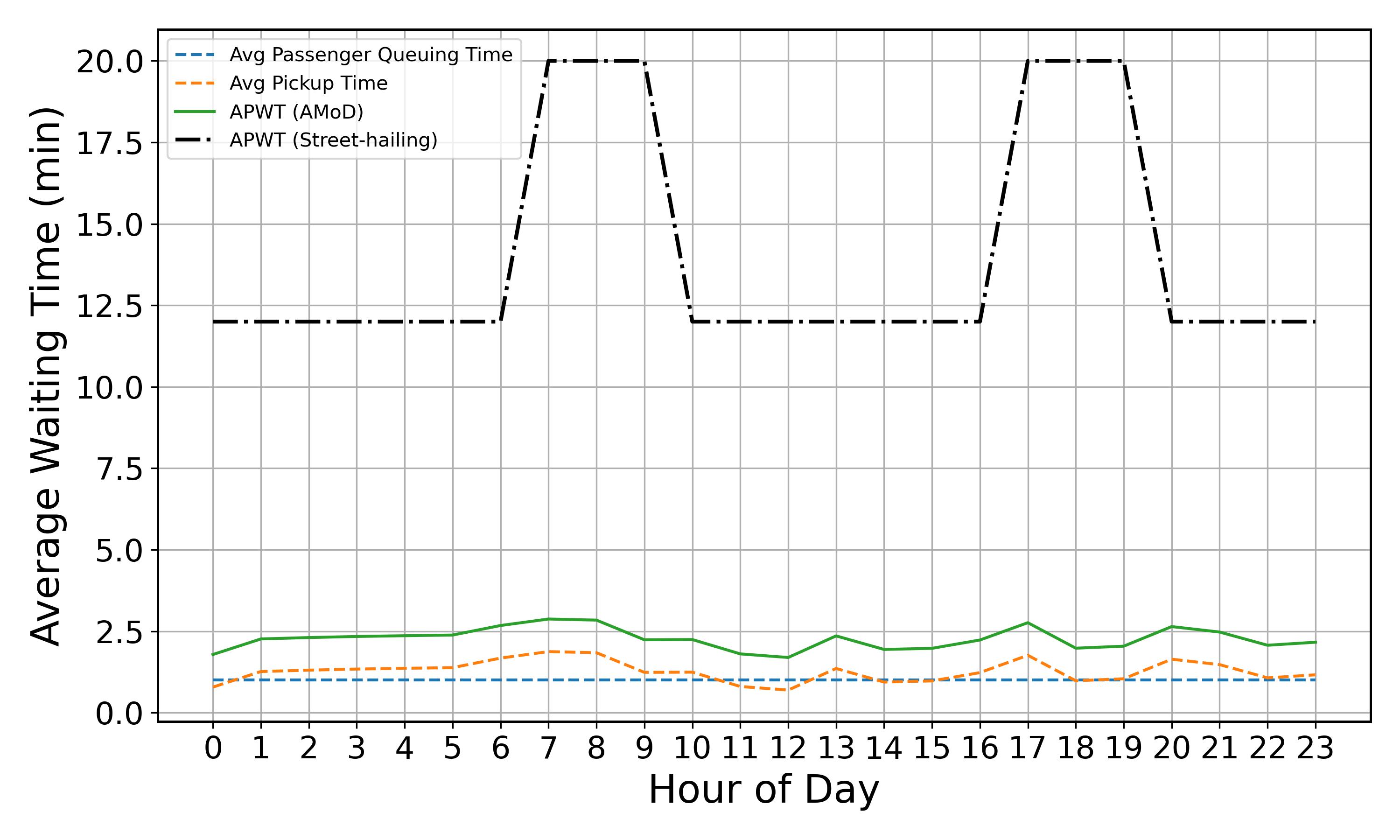}
    \caption{Comparison of average waiting times of the simulated AMoD service and the survey-based benchmark for street-hailing services~\citep{RolandBerger2015}.}
    \label{fig:hourly_waiting_time_realdata}
\end{figure}

Figure~\ref{fig:hourly_waiting_time_realdata} shows the hourly waiting time for AMoD and street-hailing services. In terms of APWT, the two services exhibit similar trends. The simulated AMoD service has lower hourly APWT values than the street-hailing benchmark throughout the day. Specifically, the APWT exceeds 2.8 minutes during the morning rush hour (i.e., from 7:00 to 9:00). This peak is caused by an increase in pickup time, which can be caused by a surge in demand that increases the burden on the available service capacity within a short time window. Considering the trend throughout the entire day, as the queuing time stays at a low level, the variation in APWT is primarily driven by fluctuations in the pickup time. Overall, under our settings, the simulated AMoD service produces lower hourly APWT values than the street-hailing benchmark.

\begin{figure}[!htbp]
    \centering
    \includegraphics[width=0.95\linewidth]{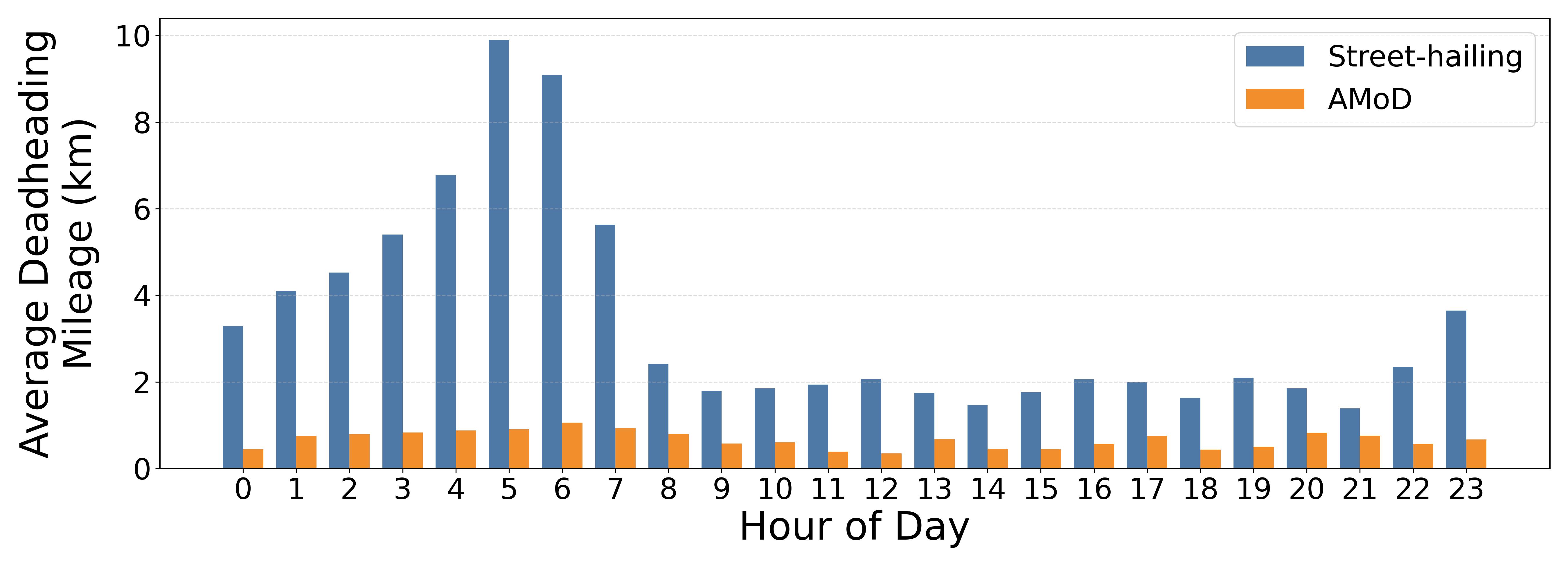}
    \caption{Comparison of ADM in the AMoD and street-hailing systems.}
    \label{fig:spatial_ADM}
\end{figure}

Figure~\ref{fig:spatial_ADM} presents the ADM comparison of the two services. The ADM in the AMoD system remains relatively stable over time, with vehicles travelling approximately 0.35--1.06~km on average to pick up passengers in each hour. This value is substantially lower than the ADM observed in street-hailing services throughout the entire day. The ADM of street-hailing services varies significantly over time. During high-demand periods, such as in daytime and early evening, street-hailing drivers only travel about 1.5–3~km to pick up their next passenger. It indicates that the numerical gap between the simulated AMoD service and the observed street-hailing benchmark is largest under low-demand conditions, because AMoD services can directly assign vehicles to passengers, eliminating the large amount of random cruising in street-hailing services when trip requests are scarce. In addition, although there is a large difference in ADM values between street-hailing and AMoD services, the variation of ADM in both hailing systems shares a similar pattern: periods of high ADM for street-hailing also correspond to high ADM for AMoD; street-hailing exhibits a pronounced peak in ADM between 4:00 and 6:00 a.m., while AMoD shows an obvious increase during the same period; both systems maintain relatively low ADM values during daytime hours, followed by a moderate rise late at night. This indicates that the substantial differences in ADM values are mainly caused by different hailing service modes, while the resemblance in their variation patterns may result from shared characteristics, such as demand patterns, road network topology, and real-time traffic speeds.

\begin{figure}[!htbp]
    \centering
    \includegraphics[width=0.85\linewidth]{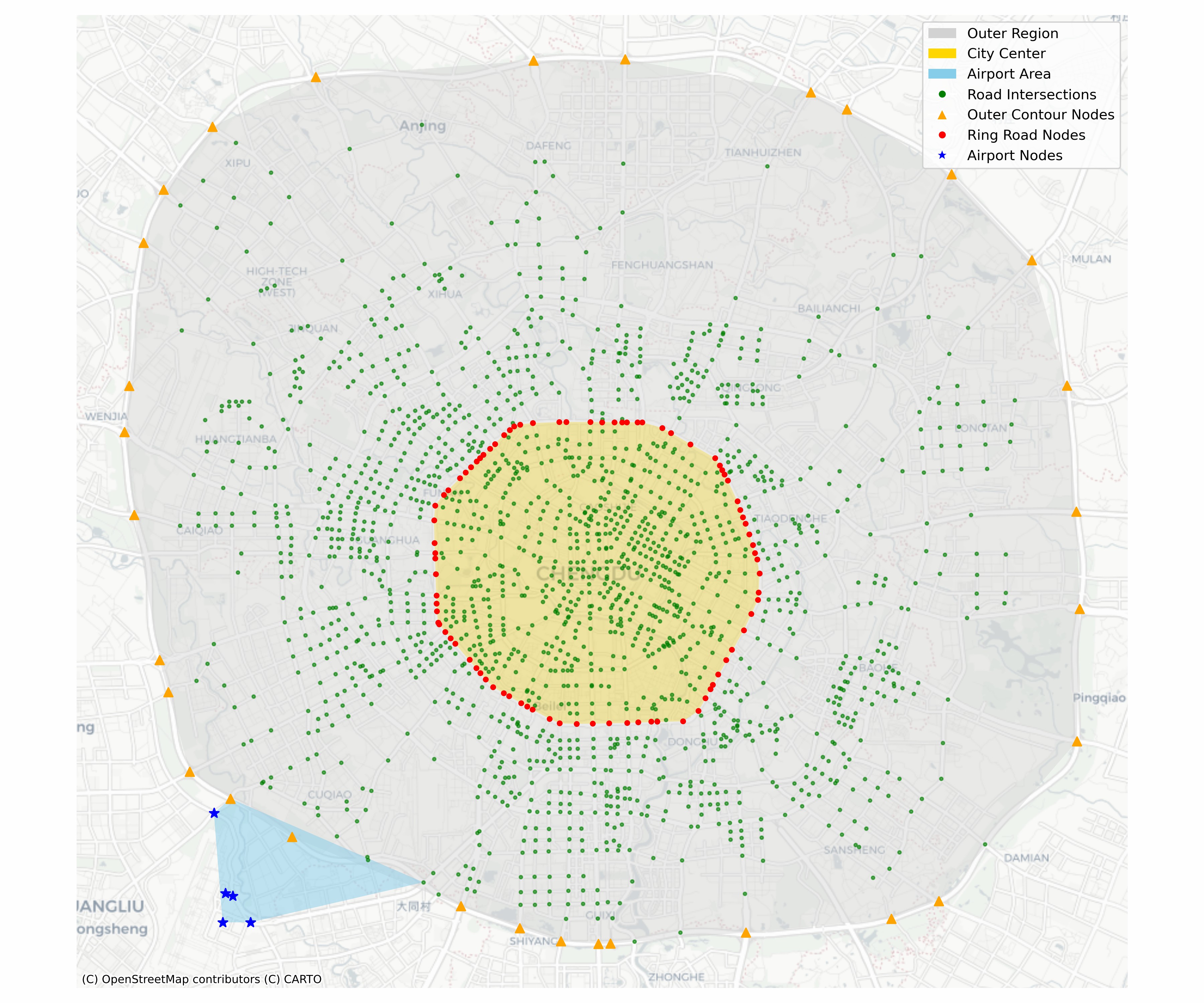}
    \caption{An illustration of the Chengdu map divided into three functional zones: Outer Region, City Centre, and Airport Area. Road intersections and boundary nodes are marked.}
    \label{fig:chengdu_region_zoning}
\end{figure}

\paragraph{Spatial analysis} To further explore the spatial heterogeneity of service performance, we divide the experimental region of Chengdu into three zones, as shown in Figure~\ref{fig:chengdu_region_zoning}. The area within and along the Second Ring Road is defined as the city centre; the area outside the Second Ring Road, excluding the neighbourhood of Shuangliu International Airport, is defined as the outer region; and the area surrounding the airport is treated as a distinct zone (i.e., airport area). This division is also consistent with the actual situation in Chengdu. We analyse all trips originating from these zones between 07:00 and 23:00, which covers the primary operational hours of on-demand travel operations.

\begin{figure}[!htbp]
    \centering
    \includegraphics[width=0.85\linewidth]{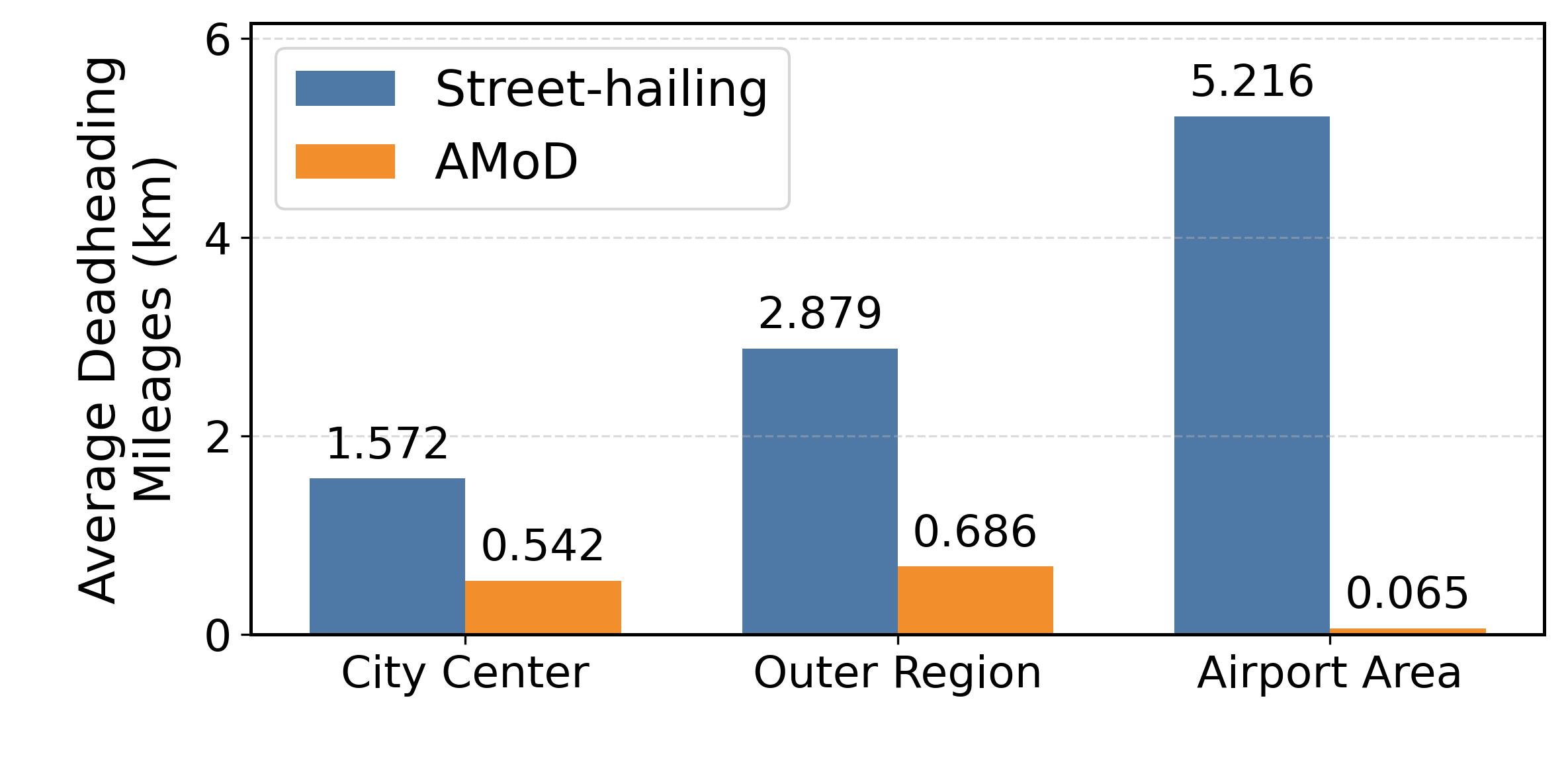}
    \caption{Comparison of ADM in different zones for AMoD and street-hailing services.}
    \label{fig:spatial_ADM_zone}
\end{figure}

Figure~\ref{fig:spatial_ADM_zone} shows the comparison of ADM in different functional zones for AMoD and street-hailing services. The simulated AMoD system yields lower values for deadheading KPIs than the corresponding street-hailing benchmarks in each functional zone. Specifically, street-hailing taxis cruise for 1.57 km to pick up a passenger in the city centre on average, while in AMoD services, autonomous vehicles travel 0.54 km. In the outer region, the ADM values of street-hailing and AMoD services are 2.88 km and 0.69 km, respectively. Notably, in the airport area, the ADM is 5.22 km for street-hailing services but only 0.07 km for AMoD, showing the largest difference among the three areas. This is because street-hailing taxi drivers often head to airports to pick up passengers intentionally, as such trips are typically longer in distance and yield higher fares. In contrast, autonomous vehicles of AMoD services are centrally controlled and cannot select trips by themselves. They usually arrive in the airport area because they receive trips directed there. Owing to the relative scarcity of trips in the vicinity, autonomous vehicles in AMoD will stay in the airport area for extended periods. Consequently, when a new demand arises, they are able to pick up the passenger with a short deadheading mileage.

In summary, under our simulation settings, including the specific demand, active-fleet-size constraints, road network, dispatch mechanism, and parameter settings, the simulated AMoD service maintains a high service rate while achieving lower APWT, ADM, and ADEC values than the corresponding street-hailing benchmarks. Based on the spatiotemporal analysis, we found that the simulated AMoD service yields better deadheading performance than the corresponding street-hailing benchmarks across the spatiotemporal conditions considered, with the largest difference occurring during late-night hours or near the airport.

\section{Discussion of AMoD operational strategies}\label{sec:discussion}
A range of operational strategies have been examined to improve the performance of AMoD services. We evaluate four representative strategies: vehicle repositioning, fleet size management, geofencing and request rejection. These analyses compare system performance under a set of pre-specified operational settings rather than designing or evaluating regulatory policy instruments. Thus, the endogenous economic or behavioural mechanisms are not considered.

\subsection{Vehicle repositioning} \label{sec:vehicle_reposition}
To examine the additional operational effect of vehicle repositioning, we further design an experiment to compare AMoD services with and without a vehicle repositioning strategy. The repositioning strategy is formulated as a minimum-cost flow problem under a non-anticipatory setting. At each decision interval, the strategy identifies vehicle-deficient nodes based on current order backlogs, short-term historical demand patterns, and the spatial distribution of vacant vehicles. A subset of vacant vehicles is then reallocated from supply-abundant nodes to deficit nodes. The details of the algorithm and parameter settings of repositioning strategies are provided in~\ref{app:repositioning}.
For the AMoD service with repositioning, ADM and ADEC account for both proactive repositioning trips and vehicle movements for passenger pickup. The experimental results, as shown in Table~\ref{tab:mode_comparison}, indicate that vehicle repositioning reduces APWT compared with the AMoD baseline without repositioning, but at the cost of increased ADM and ADEC. RD denotes the relative difference with respect to the no-repositioning scenario. Increasing repositioning intensity generally leads to more repositioning movements, reflected in higher ADM and ADEC, and thus produces a more pronounced trade-off between reduced APWT and increased ADM and ADEC. However, at high repositioning intensity, additional repositioning provides limited further reduction in APWT while substantially increasing ADM and ADEC.

\begin{table}[!htbp]
\centering
\caption{Comparison of AMoD services under different repositioning intensities}
\label{tab:mode_comparison}
\begin{threeparttable}
\footnotesize
\setlength{\tabcolsep}{4pt}
\renewcommand{\arraystretch}{1.15}

\begin{tabular}{lcccccc}
\toprule
\multirow{2}{*}{\textbf{Mode}}
& \multicolumn{2}{c}{\textbf{APWT (min)}}
& \multicolumn{2}{c}{\textbf{ADM (km)}}
& \multicolumn{2}{c}{\textbf{ADEC (kWh)}} \\
\cmidrule(lr){2-3}
\cmidrule(lr){4-5}
\cmidrule(lr){6-7}
& \textbf{Mean}& \textbf{RD (\%)}& \textbf{Mean}& \textbf{RD (\%)}& \textbf{Mean}& \textbf{RD (\%)}\\
\midrule
No repositioning& 2.2402& 0\%& 0.6184& 0\%& 0.0848& 0\%\\

Weak repositioning& 2.1805& \mbox{$-2.66\%$}& 0.6211& \mbox{$0.44\%$}& 0.0852& \mbox{$0.46\%$}\\

Moderate repositioning& 2.1310& \mbox{$-4.87\%$}& 0.6685& \mbox{$8.10\%$}& 0.0917& \mbox{$8.18\%$}\\

Strong repositioning& 2.1009& \mbox{$-6.22\%$}& 0.7102& \mbox{$14.86\%$}& 0.0975& \mbox{$14.97\%$}\\

\bottomrule
\end{tabular}

\end{threeparttable}
\end{table}

\subsection{Fleet size management} \label{sec:fleet_policy}
To evaluate how the AMoD fleet size influences its service performance and traffic conditions, we conduct a controlled experiment in which the baseline fleet size is multiplied by a scaling factor $ F $, ranging from 0.8 to 2.0 with a step size of 0.1, which simulates scenarios of fleet size reduction and expansion. Given that autonomous vehicles in AMoD constitute a significant portion of urban traffic, changes in fleet size may affect road congestion and travel speed. Specifically, a larger fleet size exacerbates urban congestion due to additional travel and curbside parking pressure, leading to a lower network speed. We follow the method introduced by~\cite{Beojone2021_inefficiency_ride_sourcing}, and adopt a simplified macroscopic fundamental diagram (MFD) to simulate road network speeds under different fleet sizes (see \ref{app:mfd}).

\begin{figure}[htbp]
    \centering
    \includegraphics[width=0.8\textwidth]{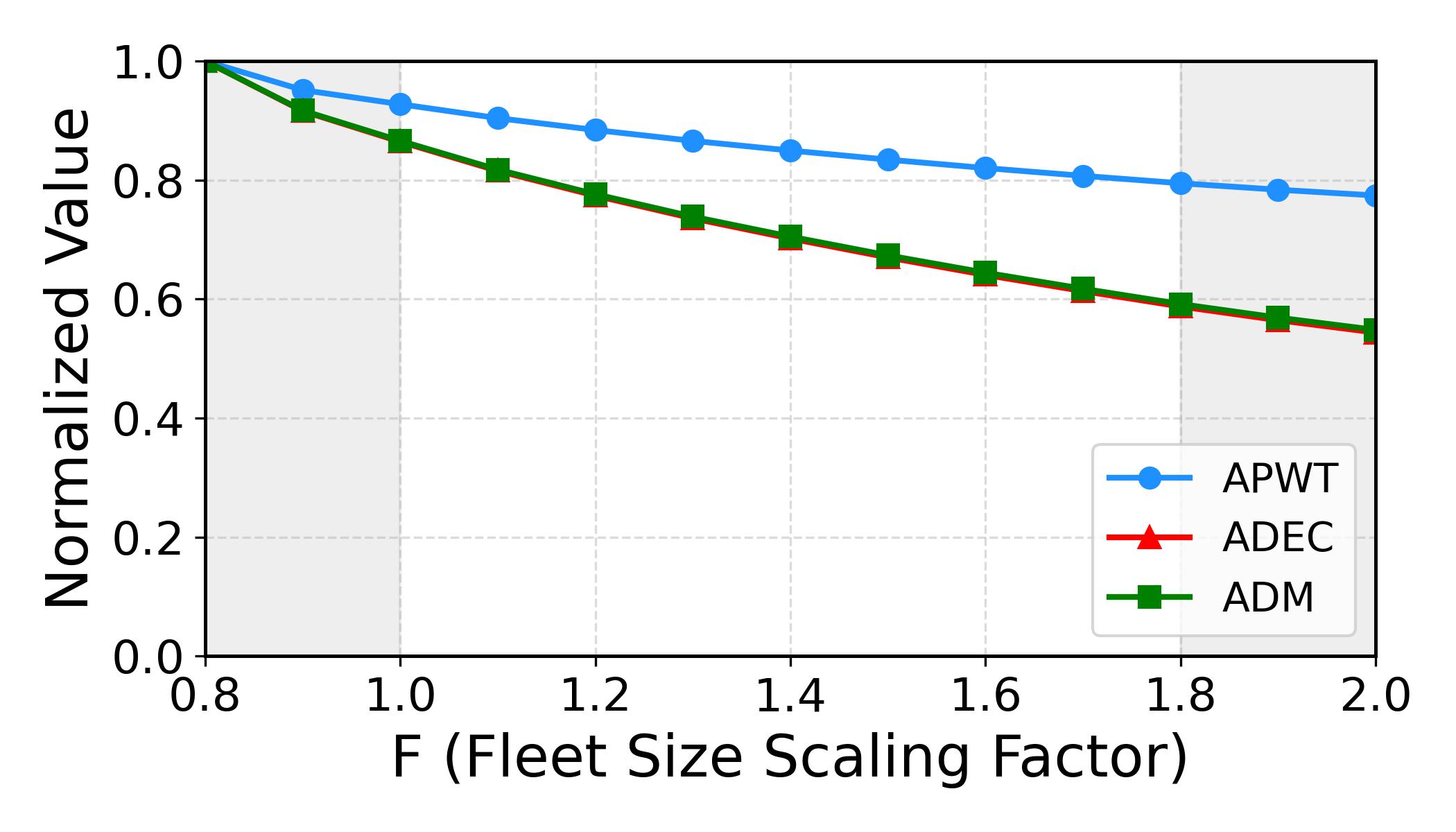}
    \caption{
    Normalized trends of APWT, ADM, and ADEC under different scaling factors of fleet size. The fleet size range $F=1.0$ to $F=1.8$ is highlighted for comparison. Overall, three indicators (APWT, ADM, and ADEC) show diminishing marginal returns as the fleet size increases.
    }
    \label{fig:fleet_tradeoff_normalized}
\end{figure}

Simulation results (Figure~\ref{fig:fleet_tradeoff_normalized}) show that reducing fleet size leads to larger APWT, ADM, and ADEC, where $F \in [0.8, 1.0)$. On average, a 0.1 increase in the fleet-size scaling factor when $F \in [0.8, 2.0]$ yields a reduction of around 2\% in APWT and reductions of around 5\% in both ADM and ADEC. This indicates that when vehicles are allowed to park on roadsides after completing their trips, the bottleneck in AMoD services is not only insufficient supply: merely increasing fleet size brings limited and diminishing benefits.

\subsection{Geofencing}
\label{sec:geofencing_policy}
As shown in Figure~\ref{fig:spatial_ADM_zone}, higher average deadheading occurs in the outer region of the city, whereas vehicles operating in the city centre exhibit higher utilization efficiency. Some cities have restricted AMoD services in specific areas such as airports or city centres, or have required vehicles to operate only within their registered administrative regions, which is also called geofencing~\citep{mot2016policy, Sun2019RidesourcingPolicy}. We thus investigate the effect of geofencing on the performance of AMoD services.

\begin{figure}[H]
    \centering
    \includegraphics[width=0.75\textwidth]{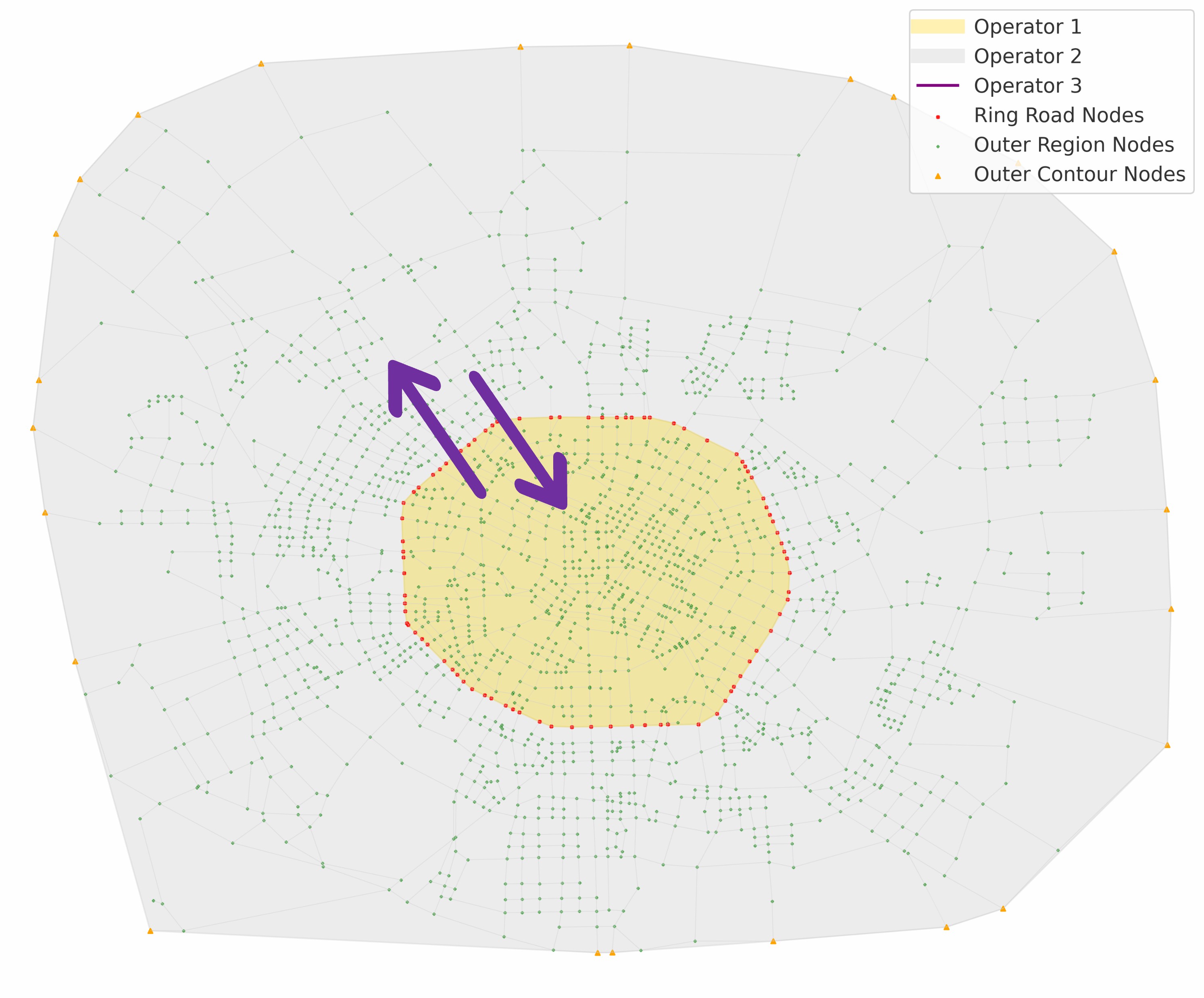}
    \caption{Service area division for Operator 1--3 in Chengdu: Operator 1 covers the city centre (within the Second Ring Road), Operator 2 covers the outer region, and Operator 3 serves cross-zone trips between the two regions.}
    \label{fig:chengdu_operator_zones}
\end{figure}

Following Section~\ref{sec:sp_analysis}, we divide Chengdu into two operational zones using the Second Ring Road as a geographic boundary: the city centre and the outer region (including the airport). This division also aligns with traffic restriction policies implemented in some megacities like Shanghai, China, where certain vehicles are not permitted to enter the city centre. We consider the geofencing with three independent operators (Operator 1--3) operating in different regions respectively, as shown in Figure~\ref{fig:chengdu_operator_zones}. Operator 1 serves trips whose origin and destination are both within the city centre; Operator 2 serves trips within the outer region; Operator 3 serves cross-zone trips that span these two regions. Each operator has a separate fleet and provides AMoD services using its dedicated platform. According to the historical trip data, trips of the three operators account for approximately 37.5\%, 25\%, and 37.5\% of the total number of trips in Chengdu, respectively. We first consider a basic geofencing that allocates the fleet to the three operators proportionally based on the share of their trips in the total number of trips, with the same $t_{\text{pu}}^{\max}$ for all operators.

We find that the geofencing slightly reduces the service capacity of the AMoD system: with the original fleet size, Operator 3 fails to serve all trips, mainly due to the sparser distribution of passenger trips, which results in longer passenger waiting and service times. Overall, we observe a 15.1\% deterioration in APWT, alongside deteriorations of 27.2\% in ADM and 27.3\% in ADEC compared with the performance of the AMoD system without implementing the geofencing. Dividing trips by operators according to their origin-destination locations in AMoD services results in a deterioration of the overall system performance, which is consistent with the findings of previous studies on multi-operator settings~\citep{Vazifeh2018}.

However, this does not imply that the geofencing is completely inefficient. Under the geofencing, AMoD performs significantly better when serving trips within the city centre (accounting for approximately 37.5\% of daily trips). Specifically, APWT, ADM, and ADEC were improved substantially by 23.2\%, 44.6\%, and 44.7\%, respectively. In addition, Operator 1 outperforms the other two operators in terms of APWT, ADM, and ADEC, while also showing relatively low vehicle utilization. In contrast, Operators 2 and 3 exhibit a certain degree of inefficiency due to more dispersed trips and longer average trip mileages.

\begin{figure}[!htbp]
    \centering
    \includegraphics[width=0.75\textwidth]{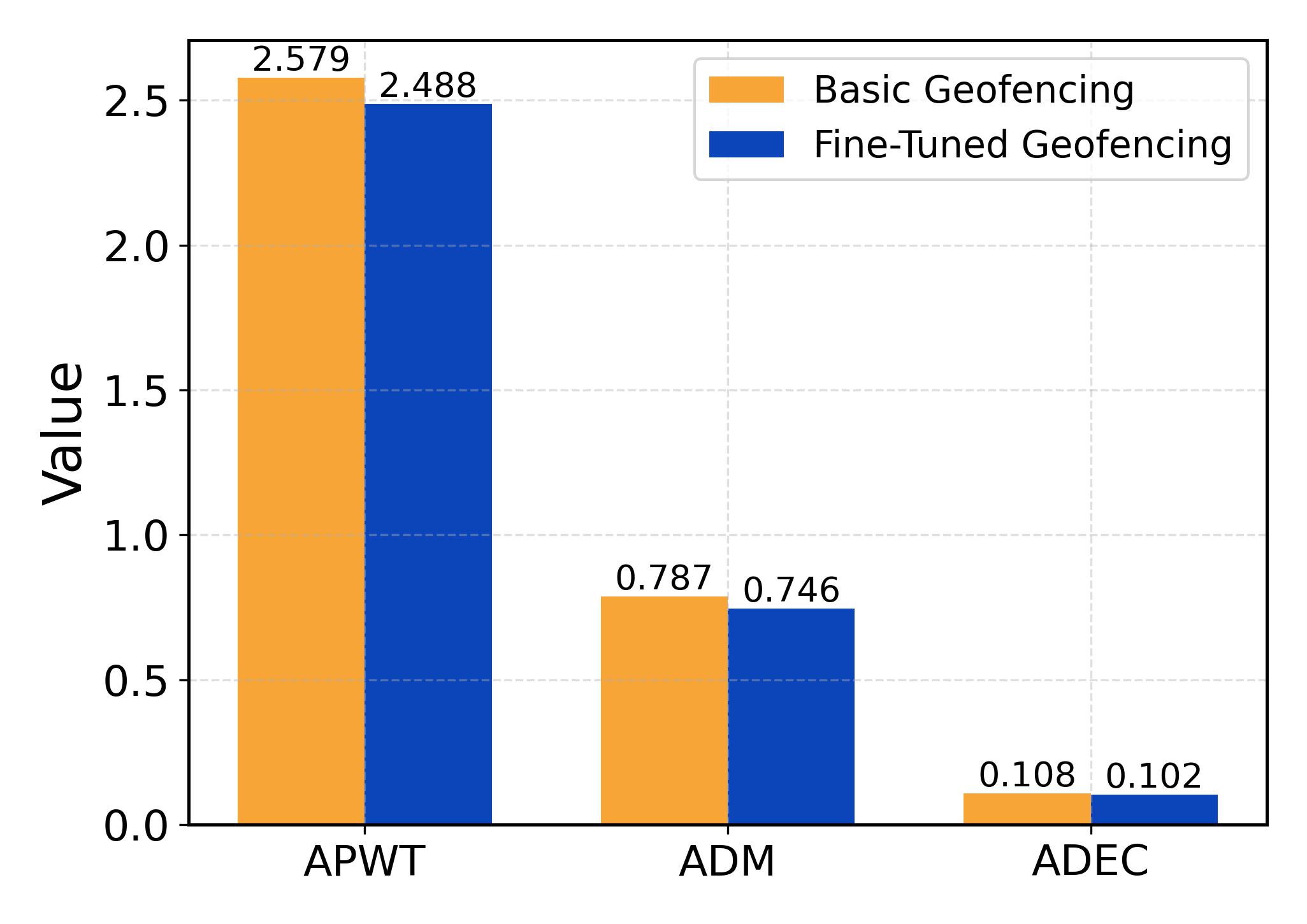}
    \caption{Performance comparison between basic and fine-tuned geofencing settings in terms of the three key performance indicators.}
    \label{fig:geofencing_comparison}
\end{figure}

In light of these results, we conduct an extended experiment in which we reduce the \( t_{\text{pu}}^{\max} \) for Operator 1 from 16 minutes to 12 minutes and increase the fleet sizes of Operators 2 and 3 by 500 vehicles each, reallocated from the fleet of Operator 1. We compare this fine-tuned geofencing with the basic one. The results show that the fine-tuned geofencing system achieves significant improvements in APWT (3.5\%), ADM (5.2\%), and ADEC (5.3\%) compared with the basic geofencing system, as shown in Figure~\ref{fig:geofencing_comparison}.

Overall, our results show that geofencing alone leads to worse performance on APWT, ADM, and ADEC, which is consistent with the result found by~\cite{Vazifeh2018} (i.e., a multi-operator system performs worse than a single-operator system). Meanwhile, the performance in serving trips within the city centre improves significantly when the geofencing treats the city centre and other areas as distinct and independent operating zones. However, this improvement is accompanied by service exclusion effects in other zones, including an increase in unserved orders and worse performance in APWT, ADM, and ADEC. Moreover, compared with the basic geofencing setting, the fine-tuned geofencing setting with more appropriately allocated fleet sizes and clearer geographic service responsibilities can achieve better performance on all three key performance indicators. This fine-tuning problem can be formulated as a fleet/resource allocation problem to achieve better overall performance.

\subsection{Request rejection}\label{sec:order_abandonment}
We further evaluate whether and how a request rejection strategy affects the performance of AMoD services. Some AMoD platforms implement demand management by rejecting or assigning certain requests with lower priority, or by indirectly stimulating or suppressing demand through such mechanisms as dynamic pricing. In this section, we specifically discuss the impacts of rejecting requests with different spatial and temporal characteristics. The amount of system resources required by a trip is determined by its spatiotemporal characteristics, including the origin and destination locations and the departure time~\citep{LiangLuo2024}. As the effect of rejecting a single trip is too small to be evaluated, we design a scenario-based experiment that removes a set of trips with particular spatiotemporal features, and then evaluate the resulting changes in system performance through controlled simulations.

For temporal characteristics, each day is divided into peak and off-peak periods, representing hours with higher or lower travel demand, respectively. For spatial characteristics, we use two labels (i.e., $x$ and $y$) to represent the spatial characteristics of each node in the graph $G$, which measure whether the node is high-demand and high-attraction. If the number of trips starting from a node is large, this node is called a high-demand node and we set $x = 1$; otherwise, it is a low-demand node and we set $x = 0$. Similarly, if the number of trips ending at a node is large, this node is called a high-attraction node and we set $y = 1$; otherwise, it is low-attraction and we set $y = 0$. Let $O_{xy}D_{xy}$ denote the spatial characteristics of a trip from its origin $O$ to the destination $D$. For example, $O_{11}D_{00}$ means that the trip's origin is high-demand and high-attraction, while the destination is low-demand and low-attraction. High-demand and high-attraction nodes are identified using a statistical outlier detection method based on Z-scores of trip origins and destinations (see \ref{app:high_nodes}).

We then conduct a sensitivity analysis based on the spatial and temporal characteristics. A total of 18 scenarios are developed. Each scenario is formed by randomly deleting 2,000 trips with a specific label (e.g., high-demand, low-attraction). Two of these scenarios delete trips started during peak or off-peak hours. The other sixteen ($2^4$) of these scenarios delete trips according to spatial characteristics (i.e., whether the origins and destinations of trips are high-demand/high-attraction or not). However, six scenarios are not analysed because the number of trips in each of these scenarios is less than 2,000. We also consider a baseline random scenario, in which 2,000 trips are randomly deleted. Because each scenario removes the same number of requests, the comparison isolates which types of requests impose greater operational burdens on the system. The resulting APWT, ADM, and ADEC are therefore calculated for the remaining served trips. To ensure the robustness of the results, each scenario is independently simulated 30 times by randomly deleting trips.

\begin{figure}[!htbp]
    \centering

    \begin{subfigure}[b]{0.75\textwidth}
        \centering
        \includegraphics[width=\textwidth]{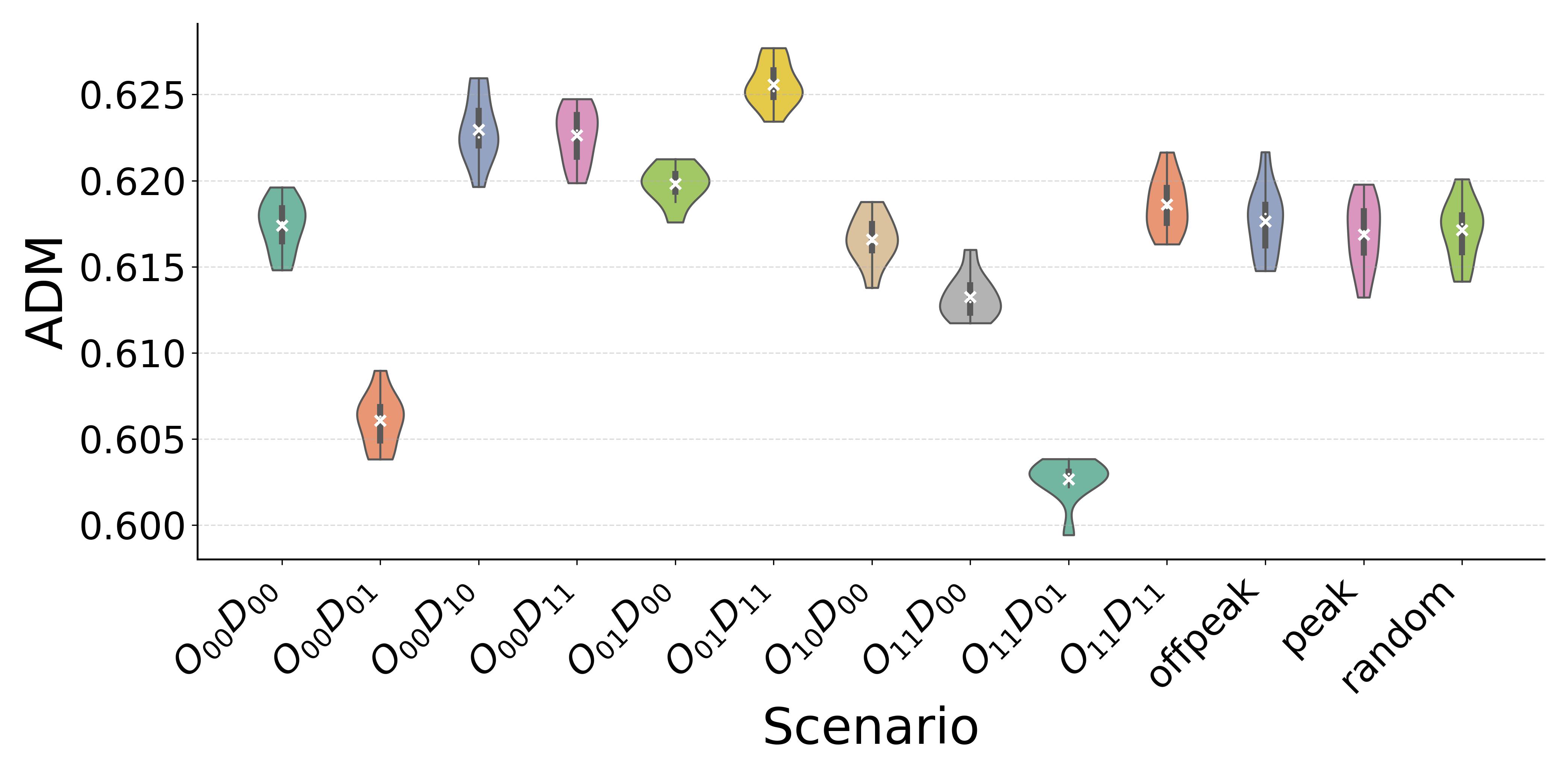}
        \caption{ADM}
        \label{fig:adm_violin}
    \end{subfigure}

    \vspace{0.4cm}

    \begin{subfigure}[b]{0.75\textwidth}
        \centering
        \includegraphics[width=\textwidth]{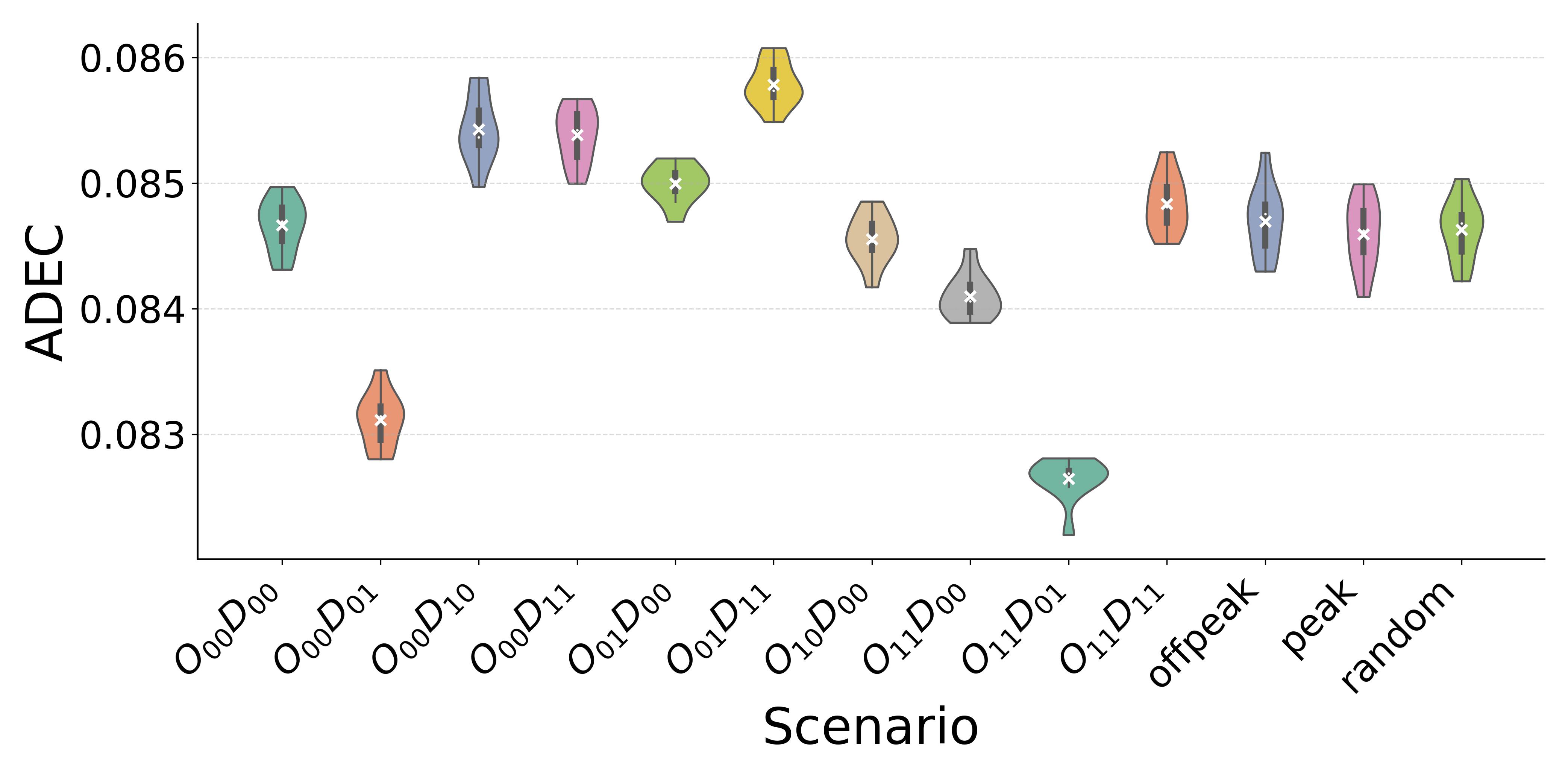}
        \caption{ADEC}
        \label{fig:adec_violin}
    \end{subfigure}

    \vspace{0.4cm}

    \begin{subfigure}[b]{0.75\textwidth}
        \centering
        \includegraphics[width=\textwidth]{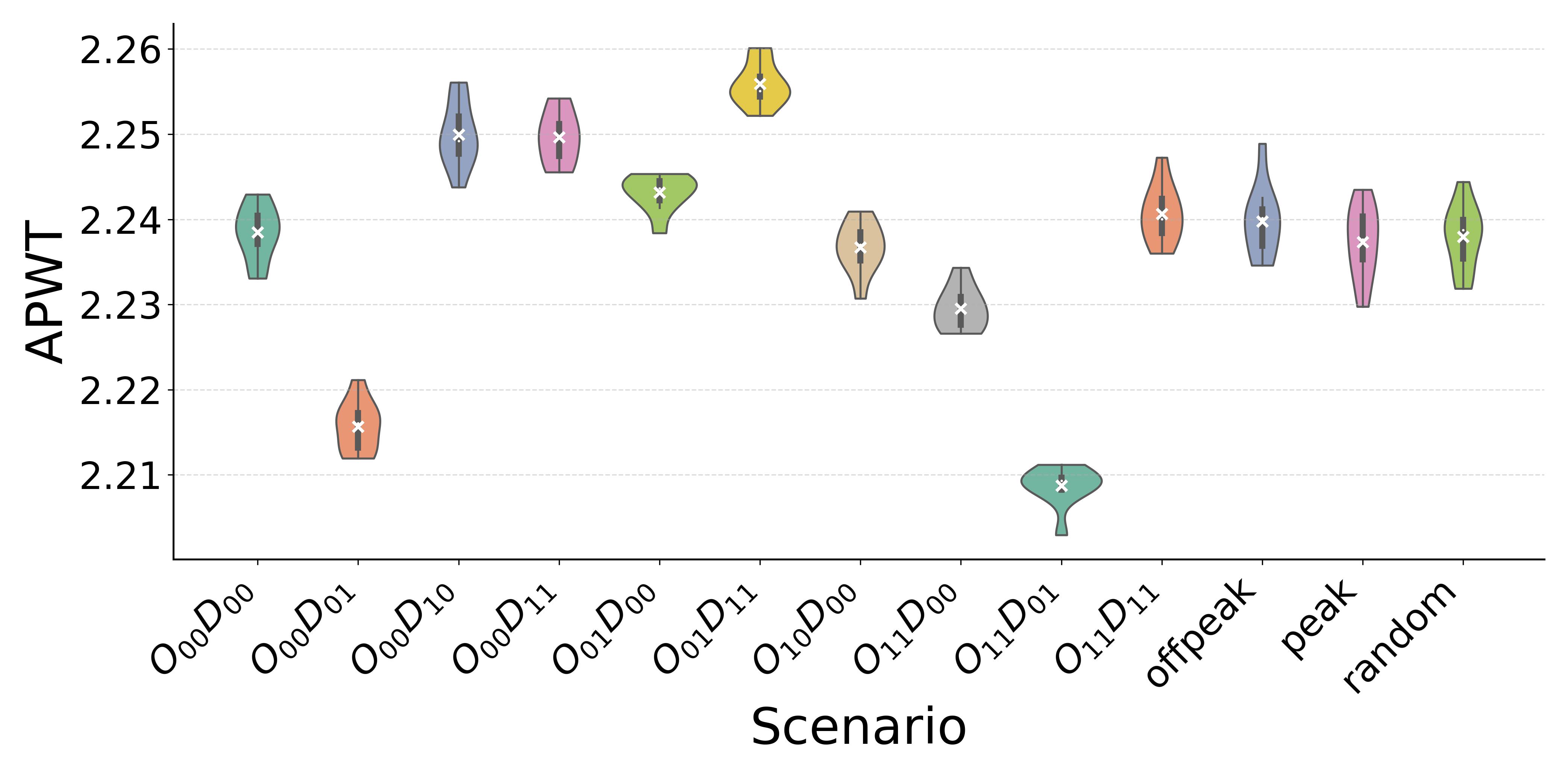}
        \caption{APWT}
        \label{fig:apwt_violin}
    \end{subfigure}

    \caption{Distributions of ADM, ADEC, and APWT under different scenarios. The horizontal line represents the median, and the white cross represents the mean.}
    \label{fig:scenario_violin_all}
\end{figure}

Overall, performance differs across scenarios that remove trips with specific spatial or temporal characteristics, as shown in Figure~\ref{fig:scenario_violin_all}. The scenarios excluding trips ending at high-attraction but not high-demand nodes, i.e., $O_{00}D_{01}$ and $O_{11}D_{01}$, contribute to statistically significant improvements (1.0\% -- 1.3\% on APWT and 1.8\% -- 2.3\% on ADM and ADEC ), compared with the baseline scenario \texttt{random}. The reason is that trips ending at high-attraction but low-demand nodes will cause a number of vehicles to wait at the low-demand node, which will further lead to a decrease in vehicle utilization. Rejecting these trips helps balance supply and demand across trips, which in turn improves all metrics. No significant performance difference is observed in scenarios rejecting trips during peak or off-peak hours. The reason may be that many nodes remain low-demand areas during peak periods, while some areas still have high demand during off-peak hours.

In summary, rejecting trips ending at high-attraction but not high-demand nodes can yield the largest reduction in passenger waiting time without increasing deadheading mileage. In contrast, the origin's characteristics appear to have minimal influence. Additionally, rejecting peak-hour demand yields negligible benefits for all performance indicators. As a result, demand management based solely on time of day is unlikely to have strong effects. AMoD operators may consider spatially differentiated request-management strategies, subject to service-equity requirements, minimum service obligations, and applicable regulatory constraints.

\section{Conclusion}\label{sec:conclusion}
This study presented a simulation framework to evaluate the potential operational performance of AMoD services and benchmark it against traditional street-hailing services under the same historical trips, fleet and road network. Three key performance indicators, including average passenger waiting time (APWT), average deadheading mileage (ADM) and average deadheading energy consumption (ADEC), were evaluated based on a real street-hailing taxi operations dataset. We considered key real-world characteristics of on-demand mobility systems, including real-time trip-vehicle matching, maximum passenger queuing time, and maximum pickup time. Our results revealed that under the simulation settings used in this study, compared with the street-hailing service, the simulated centralized dispatching AMoD service using the same historical trip demand, active fleet size, and road network is estimated to reduce APWT by 73.3\%--83.4\% relative to the survey-based street-hailing APWT scenarios. The corresponding ADM and ADEC are 75.0\% and 74.0\% lower than the GPS-derived street-hailing benchmarks, respectively. These estimated improvements are more pronounced during early-morning periods and in remote areas.

We evaluated the impacts of four operational strategies on the performance of simulated AMoD services. From the supply side, vehicle repositioning creates a trade-off between passenger waiting time and deadheading performance: it consistently reduces passenger waiting time at the cost of poorer deadheading performance; increasing fleet size exhibits gradually converging improvements, suggesting that the marginal benefits decrease as the fleet expands; geofencing improves all three performance indicators for trips within the city centre, but worsens overall system performance in the tested three-operator geofencing system. Meanwhile, tuning the fleet size and operational parameters (i.e., $ t_{\text{pu}}^{\max}$) for each operator improves the overall performance of the AMoD service. On the demand side, the results suggest that rejecting trips destined for high-attraction but not high-demand nodes may improve performance indicators among served trips under the simulated operating conditions.

This study has several limitations that point to future research directions. The AMoD simulation abstracts from some micro-level behavioural and operational frictions, such as temporary parking constraints, passenger locating delays, and minor route deviations. Addressing these issues remains an interesting direction for future work. Future work could also replace the non-anticipatory vehicle repositioning strategy with predictive or learning-based strategies that account for stochastic demand and real-time traffic conditions. The effect of ride-sharing strategies can also be evaluated through the framework in the future. In addition, the framework could be extended to incorporate economic and behavioural mechanisms, such as pricing schemes and incentive structures, to influence demand. Finally, adopting a multi-objective optimization framework and incorporating broader evaluation indicators (e.g., congestion, total VMT, accessibility) would enable a more comprehensive assessment of urban mobility system performance.

\section{Data availability}\label{sec:data_available}
Data and code can be accessed at https://github.com/youkaiwujasper/ride-hailing-dispatch.

\appendix
\renewcommand\thesection{Appendix~\Alph{section}}
\renewcommand\thefigure{\Alph{section}.\arabic{figure}}
\renewcommand\thetable{\Alph{section}.\arabic{table}}
\renewcommand\theequation{\Alph{section}.\arabic{equation}}

\makeatletter
\@addtoreset{figure}{section}
\@addtoreset{table}{section}
\@addtoreset{equation}{section}
\makeatother

\section{Estimation of the travel time on the road network}
\label{app:prediction}
We utilize a gradient descent algorithm leveraging taxi trip records to estimate the time-dependent average travel time of each road segment. Based on these edge-level travel times, the Floyd–Warshall algorithm is then employed to compute the all-pairs shortest travel time matrices across the entire network, thereby capturing the spatiotemporal accessibility dynamics over 24 hours. The procedure for training the model is shown in Algorithm~\ref{alg:travel_time_prediction}.

\begin{table}[!htbp]
\centering
\caption{Notation used in Algorithm~\ref{alg:travel_time_prediction}}
\label{tab:notation}
\begin{tabular}{ll}
\hline
\textbf{Notation} & \textbf{Definition} \\
\hline
$G=(N,E)$ & Directed road network with node set $N$ and edge set $E$ \\
$|N|$ & Number of nodes in the network \\
$|E|$ & Number of edges in the network \\
$e=(i,j)\in E$ & Directed edge from node $i$ to node $j$ \\
$\ell_e$ & Physical length of edge $e$ \\
\hline
$\mathcal{T}$ & Set of all observed trips \\
$(n_o,n_d,H,T^{obs})$ & A trip with origin $n_o$, destination $n_d$, hour $H$, observed time $T^{obs}$ \\
$\mathcal{T}_H$ & Subset of trips starting in hour $H$ \\
\hline
$t_e^H$ & Estimated travel time of edge $e$ during hour $H$ \\
$\hat{T}$ & Predicted travel time of a trip (sum of edge times) \\
$g_e$ & Gradient of loss function w.r.t. edge $e$ \\
$\eta$ & Learning rate in gradient descent \\
$\epsilon$ & Lower bound to ensure nonnegative edge times \\
\hline
$D_{i,j}^H$ & Shortest travel time from node $i$ to $j$ during hour $H$ \\
\hline
\end{tabular}
\end{table}

\begin{algorithm}[!htbp]
\caption{From Trip Records to Hourly All-Pairs Travel Time Matrices}
\label{alg:travel_time_prediction}
\begin{algorithmic}[1]
  \State \textbf{Input:} Directed road network $G=(N,E)$ with edge lengths $\ell_e$; trip dataset
    $\mathcal{T}=\{(n_o, n_d, H, T^{obs})\}$
  \State \textbf{Output:} Hourly all-pairs shortest travel time matrices $\{D^H\}_{H=0}^{23}$

  \For{each hour $H \in \{0,\dots,23\}$}
    \State \textbf{Step 1: Estimate edge travel times via gradient descent}
    \State Extract trips $\mathcal{T}_H = \{(n_o, n_d, T^{obs}) \mid (n_o, n_d, H, T^{obs}) \in \mathcal{T}\}$
    \State Initialize edge times $t_e^H \gets \dfrac{\ell_e}{v_0},\quad v_0=40~\text{km/h},\ \forall e \in E$
    \Repeat
      \State Set gradients $g_e \gets 0,\ \forall e \in E$
      \For{each trip $(n_o, n_d, T^{obs}) \in \mathcal{T}_H$}
        \State Compute shortest-time path $P(n_o,n_d)$ in $G$
        \State Predict travel time $\hat{T} = \sum_{e \in P(n_o,n_d)} t_e^H$
        \For{each $e \in P(n_o,n_d)$}
          \State $g_e \gets g_e + 2(\hat{T} - T^{obs})$
        \EndFor
      \EndFor
      \State Update edge times: $t_e^H \gets \max\{\,t_e^H - \eta g_e,\ \epsilon\,\},\ \forall e$
    \Until{loss convergence}

    \State \textbf{Step 2: Compute all-pairs travel times via Floyd–Warshall}
    \State Initialize $D_{i,j}^H \gets
      \begin{cases}
        0 & i=j \\
        t_{(i,j)}^H & (i,j) \in E \\
        \infty & \text{otherwise}
      \end{cases}$
    \For{$k \in N$}
      \For{$i \in N$}
        \For{$j \in N$}
          \State $D_{i,j}^H \gets \min\{D_{i,j}^H,\ D_{i,k}^H + D_{k,j}^H\}$
        \EndFor
      \EndFor
    \EndFor
  \EndFor
\end{algorithmic}
\end{algorithm}

\begin{figure}[!htbp]
    \centering
    \includegraphics[width=0.75\textwidth]{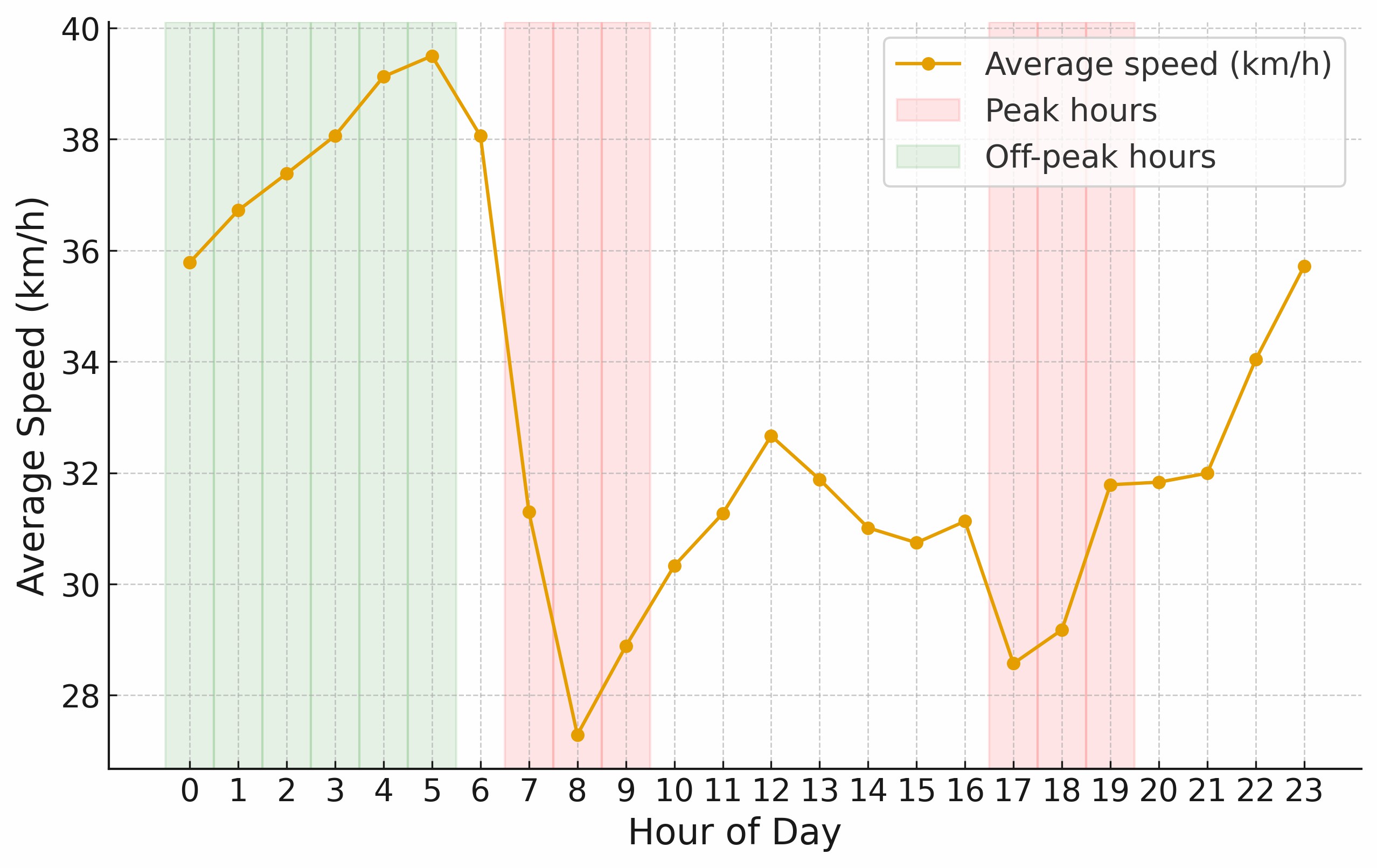}
    \caption{Hourly variation of average travel speed across the road network.}
    \label{fig:hourly_speed_curve}
\end{figure}

Figure~\ref{fig:hourly_speed_curve} illustrates the daily variation of average travel speed across the road network. The results clearly capture typical urban traffic dynamics: Speeds are highest during off-peak hours (midnight to early morning), gradually decrease during the morning commute, and reach the lowest values during the morning (around 8:00) and evening (17:00–19:00) peak periods. The recovery of speeds in late evening further confirms the method’s ability to reflect temporal fluctuations in network-wide traffic conditions.

\section{APWT analysis using a simplified street-hailing simulation}
\label{app:street-hailing_simulation}
We develop a simplified simulation model for street-hailing operations to assess the plausibility of APWT values (i.e., 8.5 -- 13.5 minutes) from the survey conducted on the same road network.

The simulation uses the same historical trips employed in the AMoD experiments. Each trip order enters the system according to its observed boarding time and joins the waiting queue at the corresponding boarding node. The initial positions of vehicles are randomly selected from the boarding nodes of the historical orders. Each vehicle can be in one of three states: cruising, serving passengers, or resting. To ensure consistency of the supply scale, the number of active vehicles is updated hourly according to the observed fleet size in the historical data.

The simulation workflow is as follows: when a vacant vehicle arrives at a node, it first checks for waiting passengers. If there is a waiting passenger at that node, the vehicle serves the first order in the waiting queue at that node; if no passengers are waiting, the vehicle randomly selects one of the downstream neighbouring nodes (including the current node), and continues cruising or remains in place.

\begin{table}[htbp]
\centering
\caption{Sensitivity analysis of the maximum passenger waiting-time threshold in the simplified cruising-taxi simulation}
\label{tab:cruising_waiting_threshold}
\begin{tabular}{ccc}
\hline
Service rate & Waiting-time threshold & Average waiting time \\
\hline
93.6\% & 20 min & 3.4 min \\
95.2\% & 60 min & 8.3 min \\
96.6\% & 120 min & 13.7 min \\
\hline
\end{tabular}
\end{table}

We conduct experiments under three maximum passenger waiting-time thresholds (i.e., 20, 60, and 120 minutes). The results (Table~\ref{tab:cruising_waiting_threshold}) show that when the service rate increases from 93.6\% to 96.6\%, the APWT increases from 3.4 to 13.7 minutes. It shows a trade-off between service rate and APWT among served trips in the street-hailing system. As for trips that cannot be served, about 5\% of the orders could not be served even if they waited for more than 60 minutes. The survey-based APWT value (8.5 -- 13.5 minutes) lies within a plausible range that reflects observed passenger tolerance in real-world operations. Extremely short waiting-time thresholds may artificially suppress the estimated waiting time by prematurely classifying requests as unserved, whereas excessively long thresholds imply unrealistic passenger patience and inflate the average waiting time.
In contrast, the survey-based APWT is based on passengers’ reported waiting times and provides a plausible benchmark for comparison. It is therefore adopted as the benchmark APWT for the street-hailing system in our comparative analysis.

\section{ADM estimation of street-hailing services}
\label{app:ADM_street}
The average deadheading mileage (ADM) of street-hailing services is estimated based on real GPS trajectory data. We use real-world GPS trajectory data from street-hailing services in Chengdu, China. Each record corresponds to a single GPS point and includes vehicle operational and positional information. Key fields of the GPS data are shown in Table~\ref{tab:gps_data_fields}.

\begin{table}[!htbp]
\centering
\caption{Fields of street-hailing GPS data}
\label{tab:gps_data_fields}
\begin{tabular}{l l}
\toprule
\textbf{Field} & \textbf{Description} \\
\midrule
\texttt{vehicle\_id} & Unique identifier for the vehicle \\
\texttt{timestamp} & GPS recorded time \\
\texttt{longitude} & Longitude of the vehicle location \\
\texttt{latitude} & Latitude of the vehicle location \\
\texttt{status} & Operational status of the vehicle: \texttt{deadheading} or \texttt{occupied} \\
\bottomrule
\end{tabular}
\end{table}

To estimate the ADM of street-hailing services, we need to identify valid deadheading trajectory segments for each vehicle and compute their travel time and distance, while excluding segments with missing or discontinuous GPS data. The notation and process are described in Table~\ref{tab:idle_notation_simple} and Algorithm~\ref{alg:idle_extraction}:

\begin{table}[!htbp]
\centering
\caption{Notation in estimating the deadheading travel time and distance}
\label{tab:idle_notation_simple}
\begin{tabular}{c l}
\toprule
\textbf{Symbol} & \textbf{Description} \\
\midrule
$t$ & Recorded time of the GPS point in the trajectory.\\
$\phi$ & Latitude of the GPS point (in radians). \\
$\lambda$ & Longitude of the GPS point (in radians). \\
$s$ & Vehicle status at the GPS point (e.g., \texttt{deadheading}). \\
\midrule
$t_{\text{start}}$ & Recorded time of the first GPS point in a deadheading segment.\\
$t_{\text{end}}$ & Recorded time of the last GPS point in a deadheading segment.\\
$\Delta t$ & Deadheading travel time of a segment.\\
$d_{12}$ & Haversine distance between two consecutive GPS points. \\
$d$ & Sum of all $d_{12}$ values within a vacant segment (in kilometres). \\
\bottomrule
\end{tabular}
\end{table}

\begin{algorithm}[!htbp]
\caption{Deadheading Travel Time and Distance Estimation from Continuous GPS Trajectories}
\label{alg:idle_extraction}
\begin{algorithmic}[1]
  \State \textbf{Input:} GPS records with $(\text{vehicle\_id},\ \text{timestamp},\ \text{latitude},\ \text{longitude},\ \text{status})$
  \State \textbf{Output:} Daily deadheading travel time and mileage of each vehicle

  \State Sort all GPS records by \texttt{vehicle\_id} and \texttt{timestamp}
  \For{each vehicle}
    \State Initialize list of valid deadheading segments as $\mathcal{S} \gets \emptyset$
    \State Initialize empty buffer $\texttt{segment} \gets []$
    \For{each GPS record $(t, \phi, \lambda, s)$}
      \If{$s = \texttt{deadheading}$}
        \State Append point to \texttt{segment}
      \Else
        \If{\texttt{segment} is not empty}
          \State Append \texttt{segment} to $\mathcal{S}$ and reset \texttt{segment} $\gets []$
        \EndIf
      \EndIf
    \EndFor
    \If{\texttt{segment} is not empty}
      \State Append final \texttt{segment} to $\mathcal{S}$
    \EndIf

    \State Initialize $\texttt{deadheading\_time\_min} \gets 0$,\quad $\texttt{deadheading\_distance\_km} \gets 0$
    \For{each segment $s \in \mathcal{S}$}
      \If{any time gap in $s$ exceeds 300 seconds}
    \State \textbf{continue}
  \EndIf
      \State Let $t_{\text{start}}$, $t_{\text{end}}$ be the first and last timestamps in $s$
      \State $\Delta t \gets (t_{\text{end}} - t_{\text{start}})$/60
      \State Initialize $d \gets 0$
      \For{each consecutive point pair $(\phi_1, \lambda_1)$ and $(\phi_2, \lambda_2)$ in $s$}
        \State Compute Haversine distance $d_{12}$ between the points
        \State $d \gets d + d_{12}$
      \EndFor
      \State $\texttt{deadheading\_time\_min} \gets \texttt{deadheading\_time\_min} + \Delta t$
      \State $\texttt{deadheading\_distance\_km} \gets \texttt{deadheading\_distance\_km} + d$
    \EndFor
    \State Record $(\texttt{vehicle\_id},\ \texttt{date},\ \texttt{deadheading\_time\_min},\ \texttt{deadheading\_distance\_km})$
  \EndFor
\end{algorithmic}
\end{algorithm}

Based on the above algorithm, we obtain the daily deadheading travel time and mileage for each vehicle. The ADM of street-hailing services is then computed by dividing the total deadheading mileage of all vehicles by the total number of trips. Thus, ADM measures the average deadheading mileage per trip and does not include passenger-carrying mileage.

\section{Sensitivity tests of modelling parameters and assumptions}
\label{app:sensitivity_analysis}

This appendix examines the sensitivity of the simulation results to three operational parameters—the batch-matching interval, maximum pickup travel time, and maximum passenger queuing time—and two operational assumptions—the fleet activation/deactivation rule and travel-time estimation factor. All settings considered are summarized in Table~\ref{tab:sensitivity_settings}. The baseline settings are a batch-matching interval of 1 minute, a maximum pickup travel time of 16 minutes, and a maximum passenger queuing time of 20 minutes. Each sensitivity scenario is evaluated using 30 independent Monte Carlo replications. Vehicle positions are independently initialized in each replication by sampling from the origin nodes of historical trips. The resulting distributions therefore provide a robustness check with respect to stochastic initial vehicle positioning.

\begin{table}[!htbp]
\centering
\caption{Parameter and assumption settings for the sensitivity analysis}
\label{tab:sensitivity_settings}
\footnotesize
\setlength{\tabcolsep}{6pt}
\renewcommand{\arraystretch}{1.18}
\begin{tabular}{@{}
>{\centering\arraybackslash}p{3.8cm}
>{\centering\arraybackslash}p{4.5cm}
>{\centering\arraybackslash}p{5.2cm}
@{}}
\toprule
\textbf{Sensitivity group}& \textbf{Parameter / assumption}& \textbf{Tested settings}\\
\midrule

\multirow{3}{3.8cm}{\centering operational parameters}
& Batch-matching interval& 1, 2, and 3 min\\

& Maximum pickup travel time& 6, 8, 10, 12, 16, and 20 min\\

& Maximum passenger queuing time& 10, 20, and 30 min\\

\midrule

\multirow{2}{3.8cm}{\centering Operational assumptions}
& Fleet activation rule& Baseline and random\\

& Travel-time estimation factor& 0.9, 1.0, and 1.1\\

\bottomrule
\end{tabular}
\end{table}

The baseline fleet activation rule refers to the state-priority rule defined in Section~\ref{sec:data_settings}. Under the random fleet activation rule, vehicles are selected without applying state-priority ordering. For the travel-time sensitivity test, all edge travel times in the baseline setting are multiplied by a common factor $\kappa$, such that 
\[
\text{\mbox{\(\displaystyle t_{ij}^{\mathrm{test}}=\kappa t_{ij}^{\mathrm{baseline}}.
\)}}\]
A factor of $\kappa=0.9$ represents travel times that are 10\% shorter than the baseline estimates, whereas $\kappa=1.1$ represents travel times that are 10\% longer.

\begin{figure}[!htbp]
    \centering
    \includegraphics[width=0.8\textwidth]{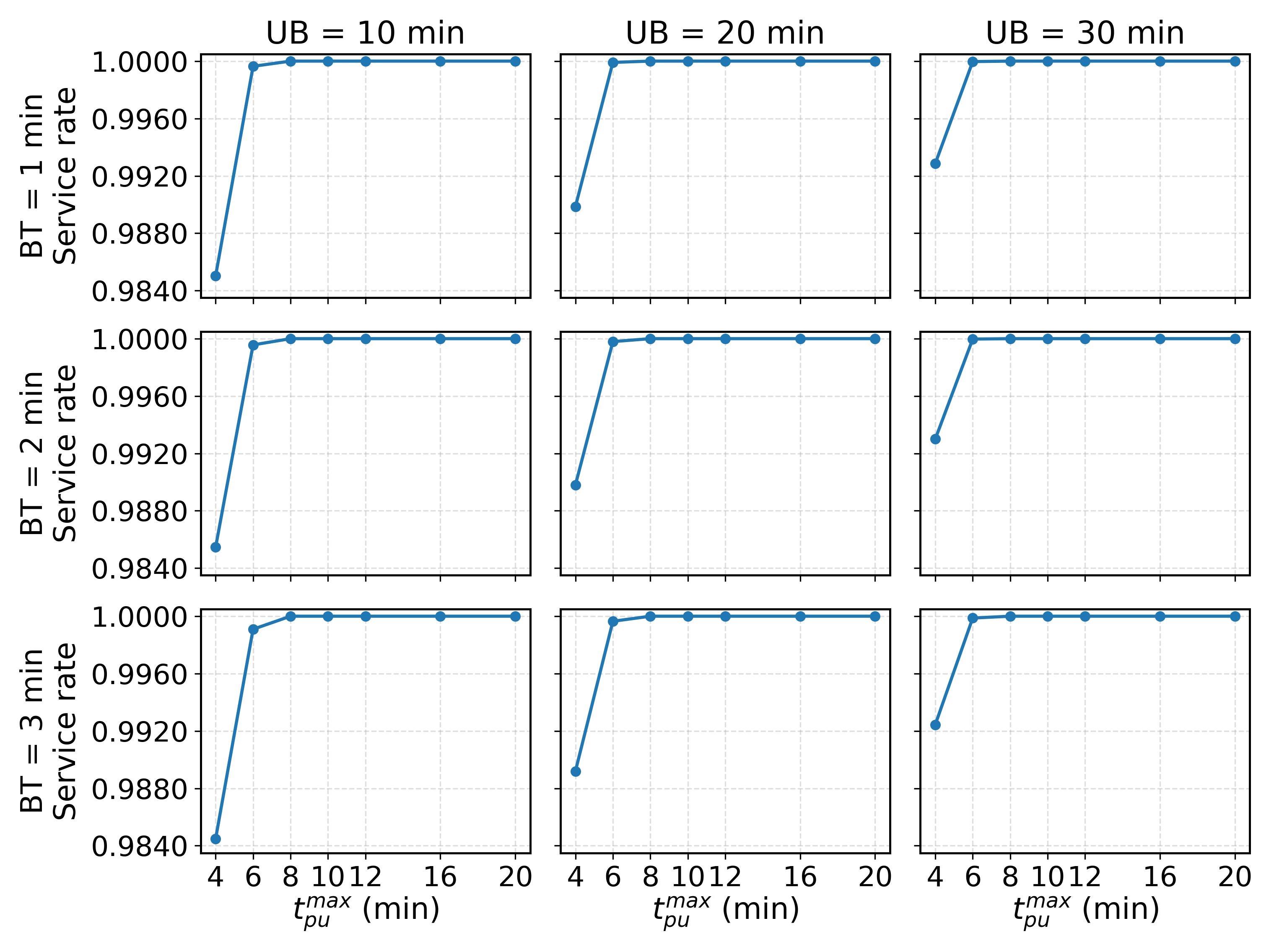}
    \caption{Sensitivity of service rate to maximum pickup travel time under different batch-matching intervals and maximum passenger queuing times.}
    \label{fig:service_rate_sensitivity}
\end{figure}

First, Figure~\ref{fig:service_rate_sensitivity} reports the service rate under different combinations of maximum pickup travel time, batch-matching interval, and maximum passenger queuing time. The results show that the service rate increases as the maximum pickup travel time is relaxed. Longer batch-matching intervals tend to reduce the service rate when the pickup travel time constraint is tight.

\begin{figure}[!htbp]
    \centering
    \includegraphics[width=0.8\textwidth]{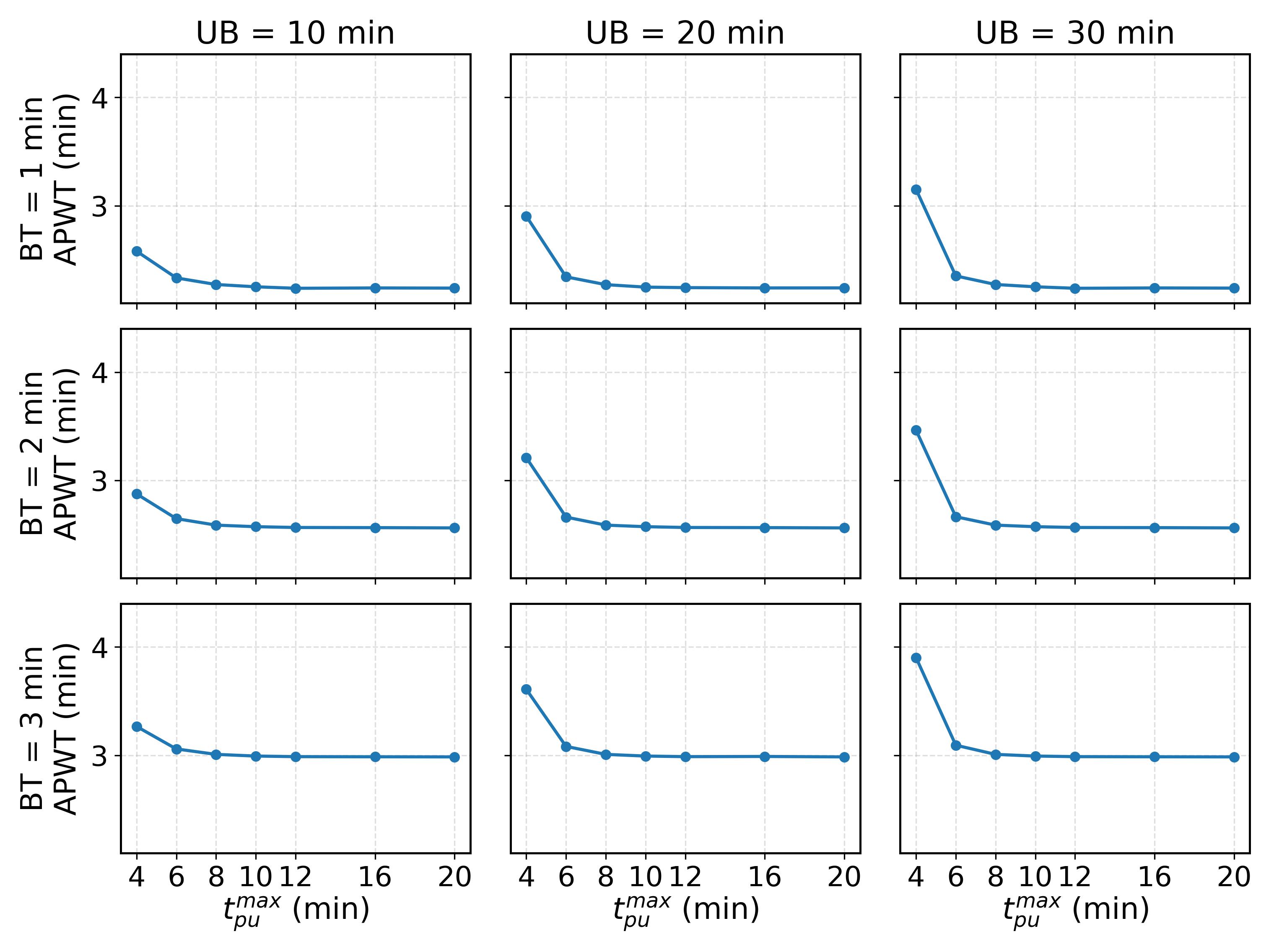}
    \caption{Sensitivity of APWT to maximum pickup travel time under different batch-matching intervals and maximum passenger queuing times.}
    \label{fig:apwt_sensitivity}
\end{figure}

\begin{figure}[!htbp]
    \centering
    \includegraphics[width=0.8\textwidth]{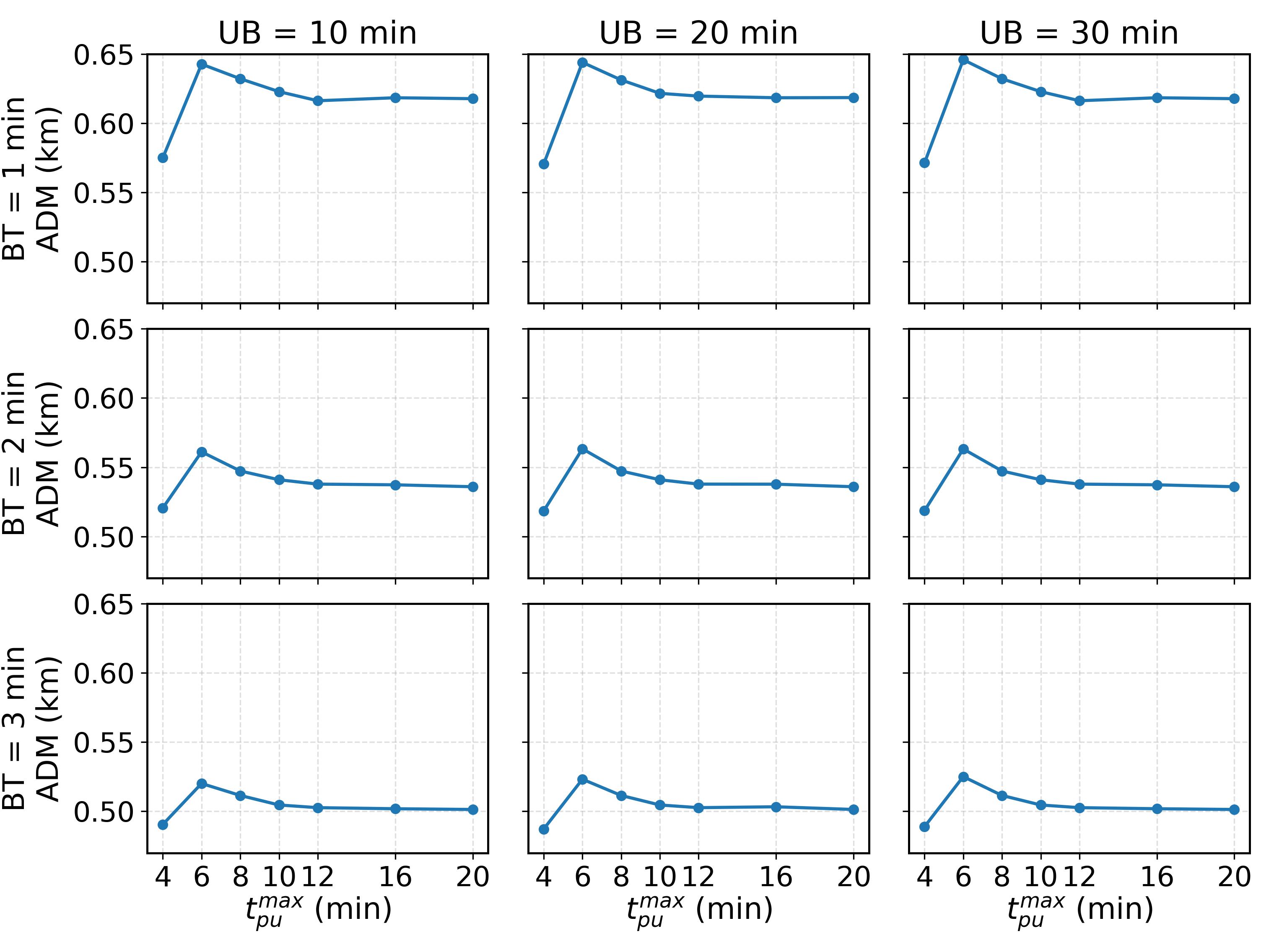}
    \caption{Sensitivity of ADM to maximum pickup travel time under different batch-matching intervals and maximum passenger queuing times.}
    \label{fig:adm_sensitivity}
\end{figure}

\begin{figure}[!htbp]
    \centering
    \includegraphics[width=0.8\textwidth]{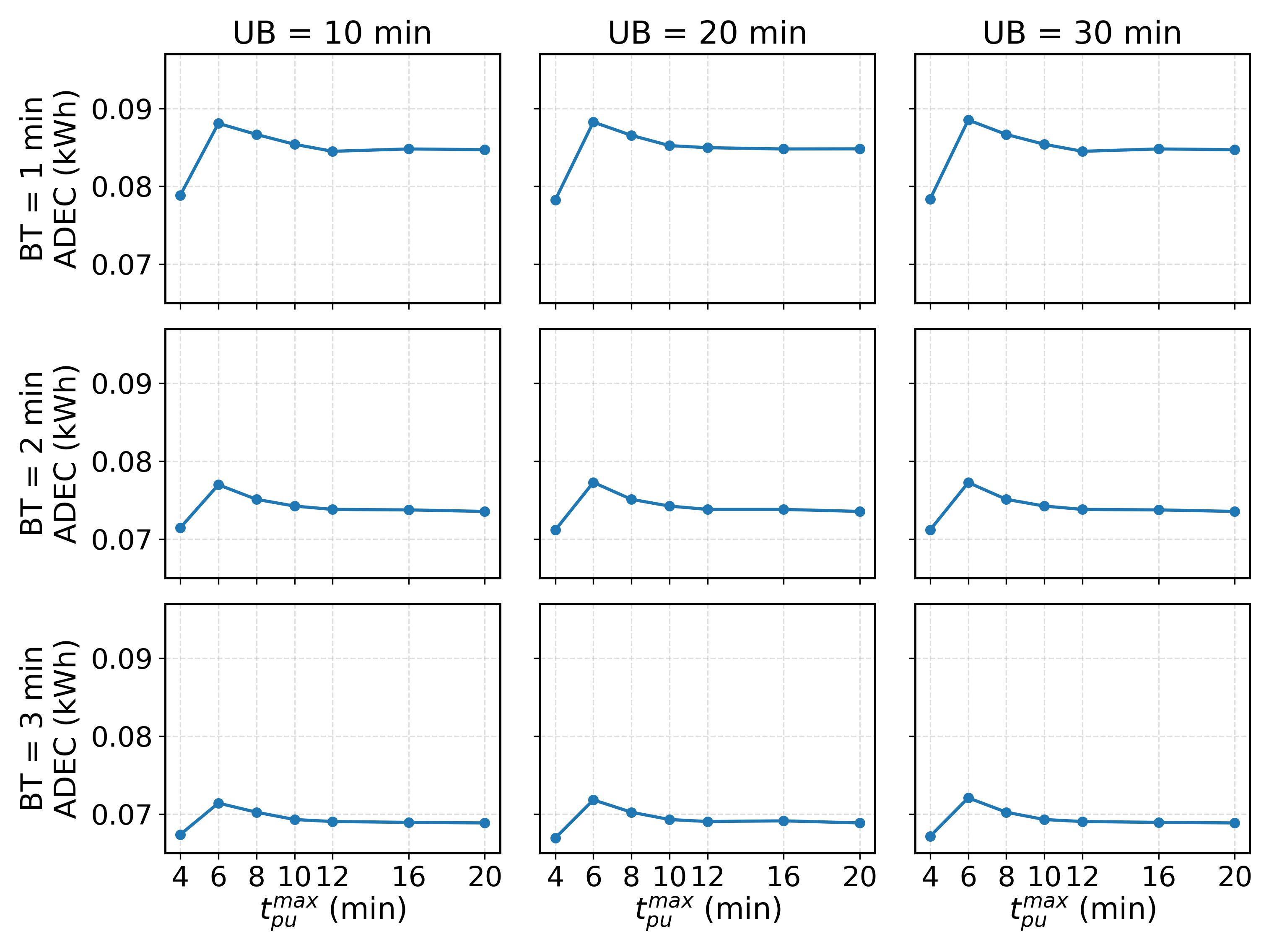}
    \caption{Sensitivity of ADEC to maximum pickup travel time under different batch-matching intervals and maximum passenger queuing times.}
    \label{fig:adec_sensitivity}
\end{figure}

Next, Figure~\ref{fig:apwt_sensitivity},~\ref{fig:adm_sensitivity}, and~\ref{fig:adec_sensitivity} report the sensitivity of APWT, ADM, and ADEC to the dispatch parameters. The results show that longer batch-matching intervals increase APWT but reduce ADM and ADEC, indicating that larger matching batches improve system operating efficiency at the expense of passenger waiting time. The maximum pickup travel time has a relatively weak effect on APWT but slightly reduces ADM and ADEC as the constraint is relaxed.

\begin{figure}[!htbp]
    \centering
    \includegraphics[width=0.95\textwidth]{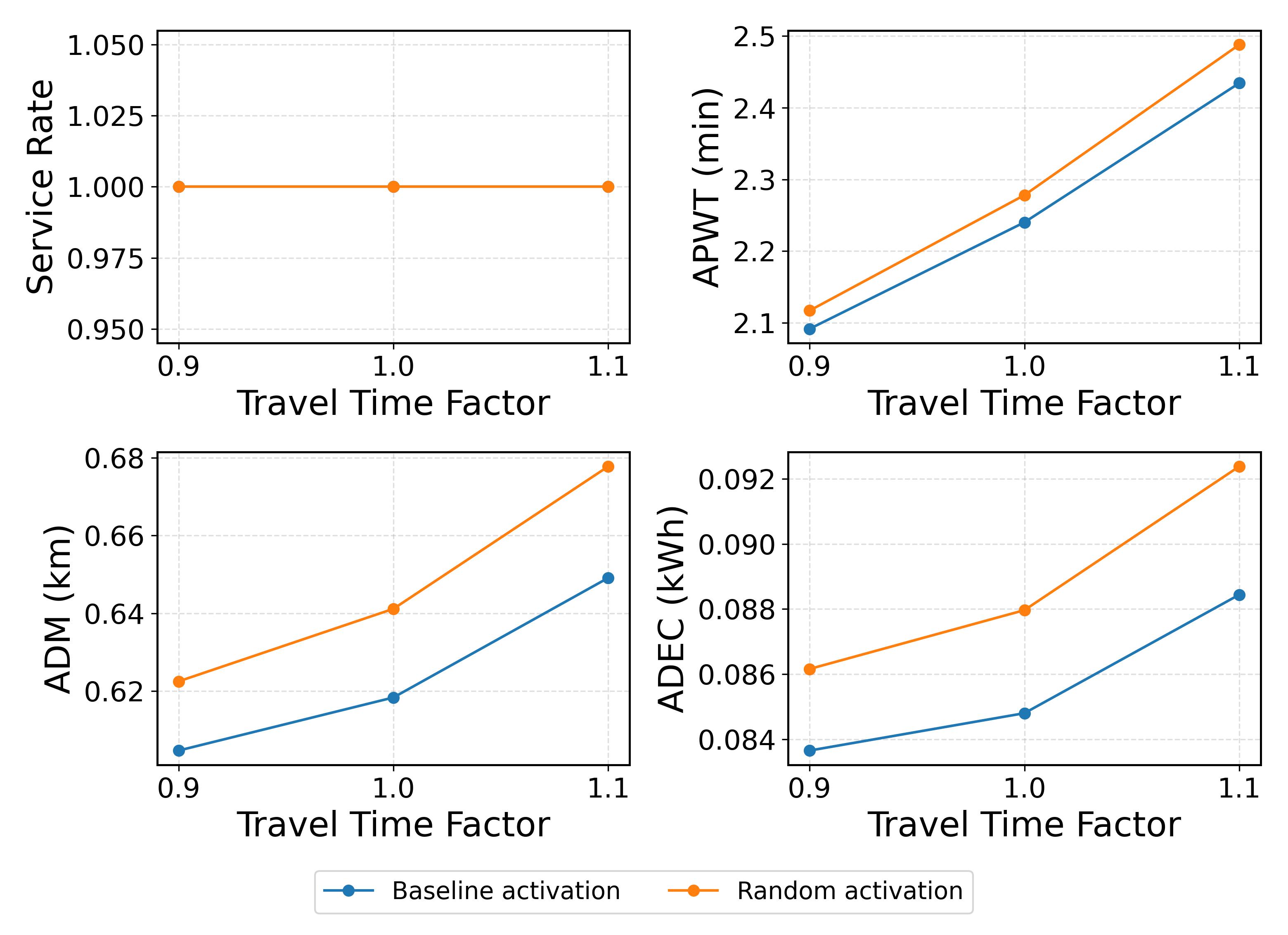}
    \caption{Sensitivity of APWT, ADM, and ADEC to fleet activation rules and travel-time estimation factors.}
    \label{fig:operational_sensitivity}
\end{figure}

Finally, Figure~\ref{fig:operational_sensitivity} reports the sensitivity of the results to the operational assumptions, including the fleet activation rule and travel-time estimation factor. In this test, the operational parameters are fixed at their baseline values: a batch-matching interval of 1 minute, a maximum pickup travel time of 16 minutes, and a maximum passenger queuing time of 20 minutes. The baseline and random fleet activation rules produce similar results. As the travel-time estimation factor increases, representing more congested traffic conditions, APWT, ADM, and ADEC generally increase. It means that slower vehicle speeds lead to longer passenger waiting times and worse deadheading performance. Meanwhile, the baseline state-priority vehicle activation setting performs better than activating/deactivating vehicles randomly. Overall, our findings remain consistent across activation rules and travel-time factors.

\section{Experimental results based on 2017 data}
\label{app:robustness_2017}

We use street-hailing service data from Chengdu, China, in 2017 to examine the robustness of the main conclusions. Specifically, we use GPS trajectory data from April 2017. For the simulated AMoD service, we apply the same modelling framework as in the main analysis, with a maximum queuing time of 30 minutes and a maximum pickup time of 24 minutes to ensure that all trips can be served. Because no survey-based APWT benchmark is available for 2017, we use the 2015 survey-based APWT range only as a reference for interpreting passenger waiting performance.

The overall comparison results between the street-hailing service and the simulated AMoD benchmark, as well as the effects of potential operational strategies, are shown in Table~\ref{tab:street_hailing_amod_baseline} and Table~\ref{tab:amod_scenario_comparison}, respectively. RD denotes the relative difference with respect to the street-hailing baseline in Table~\ref{tab:street_hailing_amod_baseline} and the AMoD baseline without repositioning or geofencing in Table~\ref{tab:amod_scenario_comparison}.

\begin{table}[!htbp]
\centering
\renewcommand{\thetable}{E.1}
\caption{Comparison between street-hailing and the AMoD baseline using 2017 data}
\label{tab:street_hailing_amod_baseline}

\begin{threeparttable}
\footnotesize
\setlength{\tabcolsep}{4pt}
\renewcommand{\arraystretch}{1.15}

\begin{tabular}{lcccccc}
\toprule
\multirow{2}{*}{\textbf{Scenario}}
& \multicolumn{2}{c}{\textbf{APWT (min)}}
& \multicolumn{2}{c}{\textbf{ADM (km)}}
& \multicolumn{2}{c}{\textbf{ADEC (kWh)}} \\
\cmidrule(lr){2-3}
\cmidrule(lr){4-5}
\cmidrule(lr){6-7}
& \textbf{Mean} & \textbf{RD}
& \textbf{Mean} & \textbf{RD}
& \textbf{Mean} & \textbf{RD} \\
\midrule
Street-hailing
& 8.5--13.5\tnote{a}
& 0\%
& 3.462
& 0\%
& 0.471
& 0\% \\
AMoD baseline
& 2.858
& $[-78.8\%, -66.4\%]$
& 0.944
& $-72.7\%$
& 0.130
& $-72.5\%$ \\
\bottomrule
\end{tabular}

\begin{tablenotes}[flushleft]
\footnotesize
\item[a] The street-hailing APWT is based on the survey-based reference
from 2015~\citep{RolandBerger2015}; all other values are based on the 2017 data.
\end{tablenotes}

\end{threeparttable}
\end{table}

\begin{table}[!htbp]
\centering
\caption{Comparison among AMoD scenarios using 2017 data}
\label{tab:amod_scenario_comparison}
\begin{threeparttable}
\footnotesize
\setlength{\tabcolsep}{4pt}
\renewcommand{\arraystretch}{1.15}

\begin{tabular}{lcccccc}
\toprule
\multirow{2}{*}{\textbf{Scenario}}
& \multicolumn{2}{c}{\textbf{APWT (min)}}
& \multicolumn{2}{c}{\textbf{ADM (km)}}
& \multicolumn{2}{c}{\textbf{ADEC (kWh)}} \\
\cmidrule(lr){2-3}
\cmidrule(lr){4-5}
\cmidrule(lr){6-7}
& \textbf{Mean}& \textbf{RD (\%)}& \textbf{Mean}& \textbf{RD (\%)}& \textbf{Mean}& \textbf{RD (\%)}\\
\midrule
AMoD baseline& 2.858& 0\%& 0.944& 0\%& 0.130& 0\%\\
AMoD with repositioning& 2.808& \mbox{$-1.8\%$}& 0.954& \mbox{$1.0\%$}& 0.131& \mbox{$1.0\%$}\\
AMoD with geofencing& 3.308& \mbox{$15.7\%$}& 1.150& \mbox{$21.8\%$}& 0.158& \mbox{$22.0\%$}\\
\bottomrule
\end{tabular}

\end{threeparttable}
\end{table}

According to Table~\ref{tab:street_hailing_amod_baseline}, in the experiment based on the 2017 data, the magnitude of improvement is comparable to that reported in the main analysis. Relative to the street-hailing reference values, APWT, ADM, and ADEC decrease by 66.4\%--78.8\%, 72.7\%, and 72.5\%, respectively.

\begin{figure}[htbp]
    \centering
    \includegraphics[width=0.8\textwidth]{\detokenize{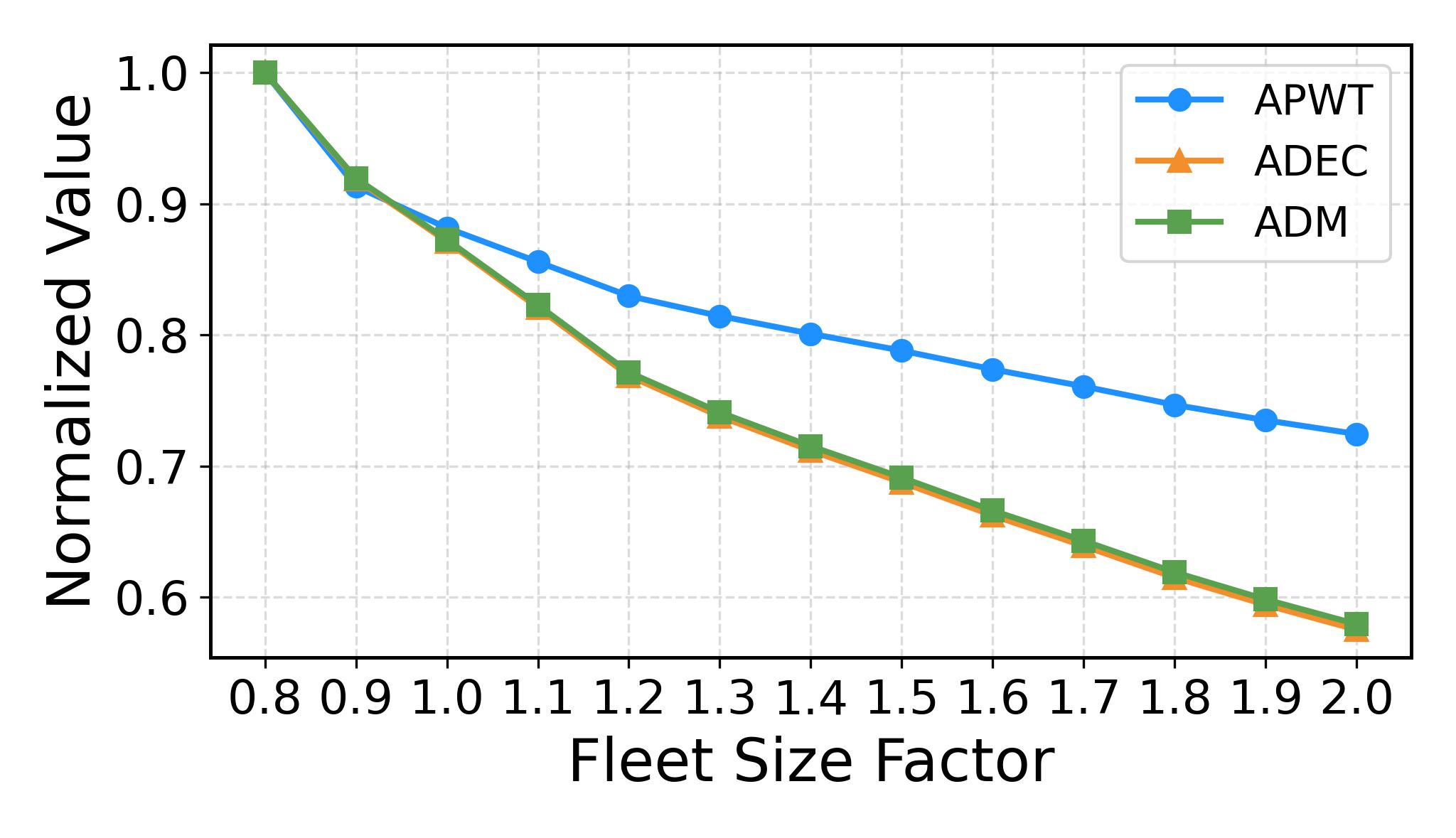}}
    \caption{Normalized values under different fleet size scaling factors using 2017 data.}
    \label{fig:normalized_values}
\end{figure}

Table~\ref{tab:amod_scenario_comparison} further shows that vehicle repositioning with weak intensity slightly reduces APWT by 1.8\%, while increasing ADM and ADEC by 1.0\%. Geofencing increases APWT, ADM, and ADEC by 15.7\%, 21.8\%, and 22.0\%, respectively. These results indicate qualitative patterns consistent with those found in the analysis based on 2015 data, particularly the trade-off between passenger waiting time and deadheading performance. Figure~\ref{fig:normalized_values} shows the effect of fleet-size variation using the 2017 data. The three indicators are normalized by the corresponding values under the baseline fleet size. The results show that increasing fleet size generally improves APWT, ADM, and ADEC, but the marginal improvements diminish as the fleet size scaling factor increases. This pattern is consistent with the findings based on the 2015 dataset.

\begin{figure}[!htbp]
    \centering

    \begin{subfigure}[b]{0.75\textwidth}
        \centering
        \includegraphics[width=\textwidth]{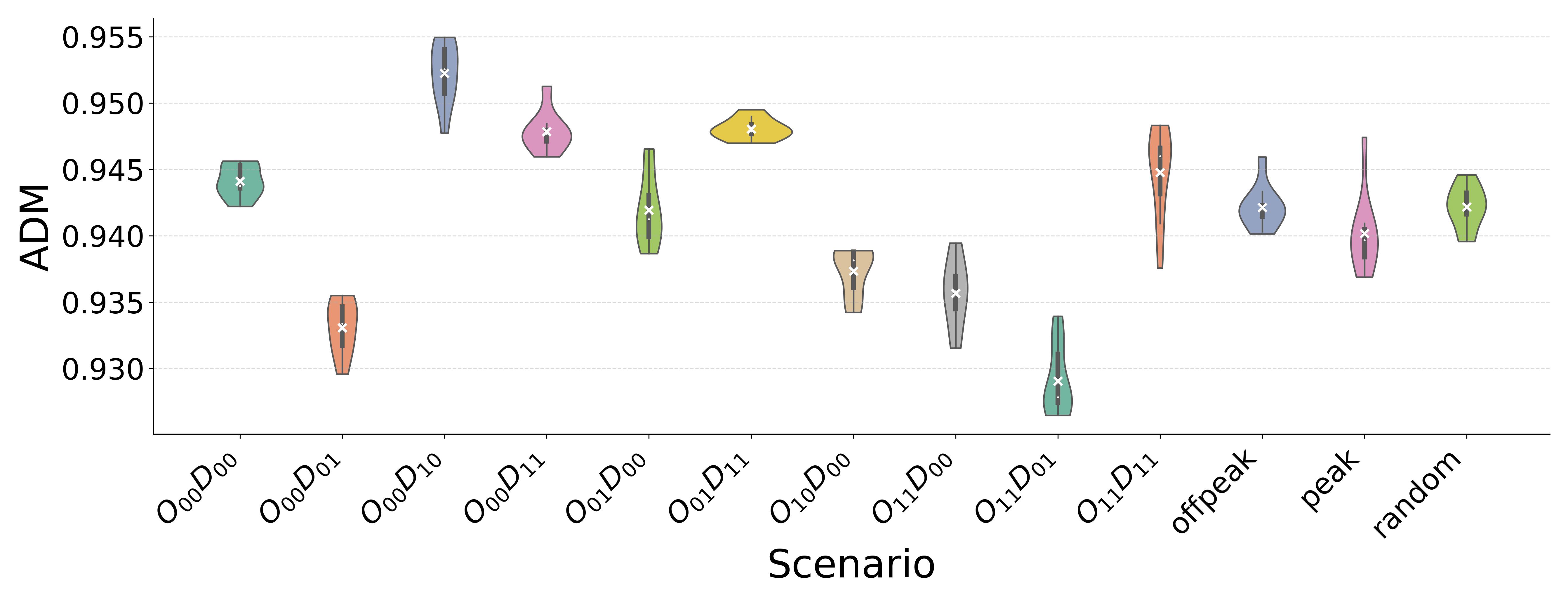}
        \caption{ADM}
        \label{fig:demand_removal_adm_violin}
    \end{subfigure}

    \vspace{0.4cm}

    \begin{subfigure}[b]{0.75\textwidth}
        \centering
        \includegraphics[width=\textwidth]{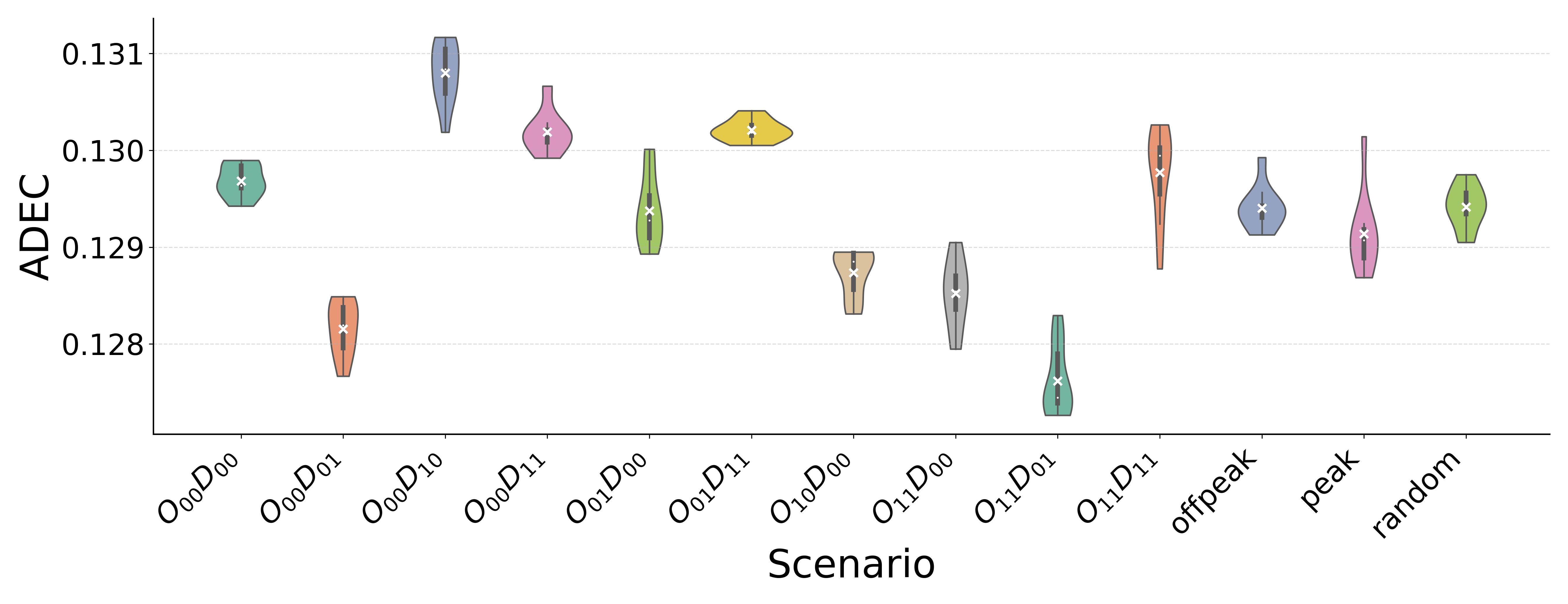}
        \caption{ADEC}
        \label{fig:demand_removal_adec_violin}
    \end{subfigure}
    \vspace{0.4cm}

    \begin{subfigure}[b]{0.75\textwidth}
        \centering
        \includegraphics[width=\textwidth]{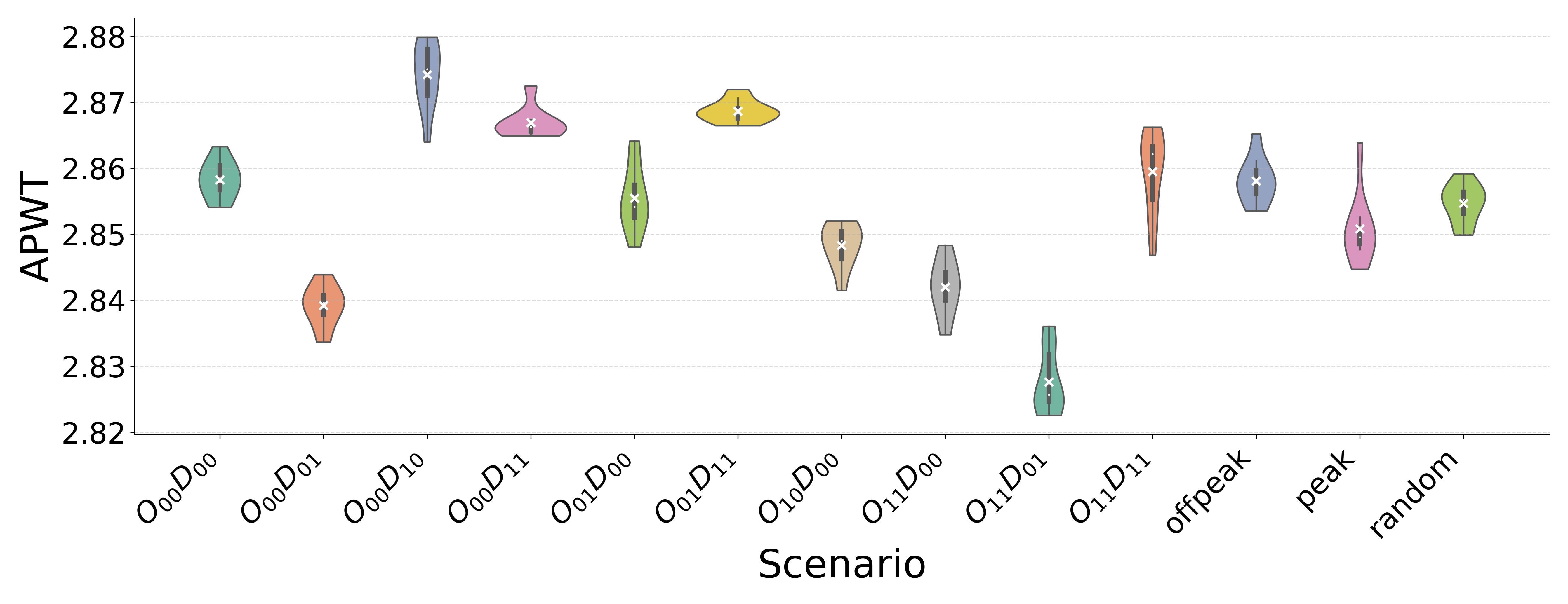}
        \caption{APWT}
        \label{fig:demand_removal_apwt_violin}
    \end{subfigure}

    \caption{Distributions of ADM, ADEC, and APWT under different request-rejection scenarios using the 2017 data. The horizontal line represents the median, and the white cross represents the mean.}
    \label{fig:demand_removal_violin_all}
\end{figure}

We further test the effect of request rejection using the 2017 dataset. We repeat the request-rejection experiment using the same scenario definitions as the analysis based on 2015 data. In each scenario, 2,000 requests with specified spatial or temporal characteristics are rejected, and APWT, ADM, and ADEC are calculated for the remaining served trips. Each scenario is independently simulated 30 times. As shown in Figure~\ref{fig:demand_removal_violin_all}, rejecting requests ending at high-attraction but low-demand nodes results in lower APWT, ADM, and ADEC than the random-removal baseline. This qualitative pattern is consistent with the results obtained using the 2015 dataset.

Overall, the results based on the 2017 dataset are broadly consistent with those based on the 2015 dataset. The AMoD baseline substantially improves APWT, ADM, and ADEC relative to the street-hailing benchmark, and the additional operational-scenario analyses show similar qualitative patterns. These results indicate that the main findings show a degree of temporal robustness across the 2015 and 2017 datasets.

\section{Mathematical details of the minimum cost flow model for vehicle repositioning }
\label{app:repositioning}
We design a vehicle repositioning model based on the minimum cost flow formulation to quantify impacts of the vehicle repositioning strategy on the operational performance of AMoD systems. This model implements a responsive and localized repositioning scheme. The framework operates in a rolling simulation manner without using predictive or future information. It evaluates vehicle shortage pressure at each network node using current unmatched orders, the real-time spatial distribution of vacant and repositioning vehicles, and historical demand weighted by an adjustable parameter.

At each repositioning period, a node is defined as vehicle-deficient if current queued orders and historical demand exceed the supply of currently and historically available vehicles. As a setting, only vehicles in the vacant state are accessible for repositioning. Occupied vehicles, resting vehicles, and those already under repositioning are excluded. For each node with surplus vacant vehicles, a number of vehicles can be reserved at the node, while the remaining vehicles are treated as candidates for repositioning. In each repositioning process, the maximum proportion of vehicles that can be dispatched is controlled by a tunable parameter.

Specifically, the shortage pressure of node $j$ at repositioning epoch $t$ is defined as:
\begin{equation}
\Delta_j^t = Q_j^t + \omega D_j^t - H_j^t - S_j^t - R_j^t,
\end{equation}
where $Q_j^t$ is the number of unmatched orders at node $j$, $D_j^t$ is the number of historical orders at node $j$ within the previous historical window, $H_j^t$ is the average number of vacant vehicles at node $j$ within the previous historical window, $S_j^t$ is the current number of vacant vehicles at node $j$, $R_j^t$ is the number of vehicles currently en route to node $j$ for repositioning, and $\omega$ is the weight of historical demand.

The repositioning demand of node $j$ is then given by:
\begin{equation}
d_j^t = \left\lceil \max(0, \Delta_j^t) \right\rceil.
\end{equation}

For each origin node $i$, the number of movable vacant vehicles is defined as:
\begin{equation}
m_i^t = \max(0, S_i^t - \gamma),
\end{equation}
where $\gamma$ is the minimum number of vacant vehicles reserved at each node.

The total repositioning budget at epoch $t$ is:
\begin{equation}
B_t = \left\lfloor \rho \sum_{i\in N} S_i^t \right\rfloor,
\end{equation}
where $\rho$ is the maximum proportion of currently vacant vehicles that can be repositioned in one repositioning epoch.

Based on these definitions, the repositioning decision is formulated as a minimum cost flow problem. In this model, covering vehicle-deficient nodes is treated as a benefit, while the travel time of vehicle repositioning is regarded as the cost. Vehicles can only be repositioned to adjacent neighbourhood nodes of their current location, with a single repositioning duration not exceeding a threshold. If the net benefit from the origin node to the target node is positive—i.e., the benefit of covering the vehicle-deficient node outweighs the travel time cost—this movement edge is incorporated into the minimum cost flow network. After solving the problem, the selected vehicles enter the repositioning state and cannot accept new orders until they reach the target node. Upon arrival, their positions are updated to the target node, and they are reset to the vacant state. The minimum cost flow model is formulated as:

\begin{equation}
\min_{r} \sum_{(i,j)\in A_t^R} \left( \beta \tau_{ij}^t - \alpha \right) r_{ij}^t
\end{equation}

\noindent Subject to:

\begin{equation}
\sum_{j:(i,j)\in A_t^R} r_{ij}^t \leq m_i^t, \qquad \forall i \in I_t
\end{equation}

\begin{equation}
\sum_{i:(i,j)\in A_t^R} r_{ij}^t \leq d_j^t, \qquad \forall j \in J_t
\end{equation}

\begin{equation}
\sum_{(i,j)\in A_t^R} r_{ij}^t \leq B_t
\end{equation}

\begin{equation}
r_{ij}^t \in \mathbb{Z}_{\geq 0}, \qquad \forall (i,j)\in A_t^R
\end{equation}

Here, $A_t^R$ only includes feasible repositioning arcs satisfying
$(i,j)\in E$, $\tau_{ij}^t \leq \tau_{\max}^R$, and $\alpha-\beta\tau_{ij}^t>0$. $r_{ij}^t$ denotes the number of vacant vehicles repositioned from origin node $i$ to destination node $j$ at repositioning epoch $t$; $A_t^R$ denotes the set of feasible repositioning arcs; $\tau_{ij}^t$ denotes the travel time from node $i$ to node $j$ at time $t$; $\alpha$ denotes the reward for covering one unit of vehicle shortage; $\beta$ denotes the travel-time cost coefficient of repositioning; $I_t$ denotes the set of origin nodes with movable idle vehicles; $J_t$ denotes the set of destination nodes with positive vehicle shortage; $m_i^t$ denotes the number of movable idle vehicles at node $i$; $d_j^t$ denotes the repositioning demand at node $j$; and $B_t$ denotes the maximum number of vehicles allowed to be repositioned at epoch $t$.

The parameters used in this study are shown in Table~\ref{tab:repositioning_parameters} including sensitivity tests of alternative repositioning settings. A stronger repositioning strategy indicates a tendency for more repositioning actions.

\begin{table}[!htbp]
\centering
\caption{Parameter settings for the minimum-cost-flow-based vehicle repositioning sensitivity analysis}
\label{tab:repositioning_parameters}
\footnotesize
\setlength{\tabcolsep}{4pt}
\renewcommand{\arraystretch}{1.15}
\begin{tabular}{@{}>{\centering\arraybackslash}p{2.4cm}
                >{\centering\arraybackslash}p{2.0cm}
                >{\centering\arraybackslash}p{1.8cm}
                >{\centering\arraybackslash}p{1.7cm}
                >{\centering\arraybackslash}p{2.0cm}@{}}
\toprule
\textbf{Parameter}
& \textbf{No RP}
& \textbf{Weak}
& \textbf{Moderate}
& \textbf{Strong} \\
\midrule
$\Delta^R$
& 5 min
& 5 min
& 3 min
& 2 min \\

$W$
& 30 min
& 30 min
& 30 min
& 30 min \\

$\tau_{\max}^R$
& 0 min
& 15 min
& 20 min
& 30 min \\

$\rho$
& 0
& 0.3
& 0.6
& 1.0 \\

$\alpha$
& 0
& 20
& 30
& 100 \\

$\beta$
& 1
& 1
& 1
& 1 \\

$\omega$
& 0
& 0.5
& 0.5
& 0.5 \\

$\gamma$
& 0
& 1
& 0
& 0 \\
\bottomrule
\end{tabular}
\end{table}

\section{Simulating road network speeds under different fleet sizes based on the macroscopic fundamental diagram model}
\label{app:mfd}
The average travel speed on the road network is computed based on the total travel time and mileage using the macroscopic fundamental diagram (MFD) model. We construct the MFD data points by aggregating vehicle operation records every 30 minutes across three consecutive days. For each 30-minute interval, we compute the total number of active vehicles in the street-hailing system and assume that street-hailing vehicles account for 10\% of total vehicles on the road network, which is a normal allocation ratio scenario \citep{Ding2022TRC}. $v$ and $n$ represent the average road network speed and the number of vehicles running on the road, respectively. The aggregated data points $(n, v)$ are then used to fit multiple candidate functional forms for estimating speeds $v(n)$, including linear, exponential decay, power law, and log-shift models (Algorithm~\ref{alg:mfd_regression}).

\begin{algorithm}[!htbp]
\caption{Macroscopic Fundamental Diagram (MFD) Regression Procedure}
\label{alg:mfd_regression}
\begin{algorithmic}[1]
  \State \textbf{Input:} Vehicle trip records from multiple days
  \State \textbf{Output:} Regression parameters for $v(n)$
  \State Assign each record to a 30-minute time interval $t$

  \For{each interval $t$}
    \State Compute active taxi count $n_t^{\text{taxi}}$
    \State Estimate total vehicle count $n_t \gets n_t^{\text{taxi}} / 0.1$
    \State Compute total travel time and mileage
    \State Calculate average speed $v_t$ (km/h)
  \EndFor

  \State Construct MFD points: $\mathcal{P} = \{(n_t, v_t)\}$
  \For{each candidate model $f \in \{\text{linear}, \text{exponential}, \text{power}, \text{log-shift}\}$}
    \State Estimate parameters $\theta_f$ via non-linear regression
  \EndFor
\end{algorithmic}
\end{algorithm}

\begin{figure}[!htbp]
    \centering
    \includegraphics[width=0.75\textwidth]{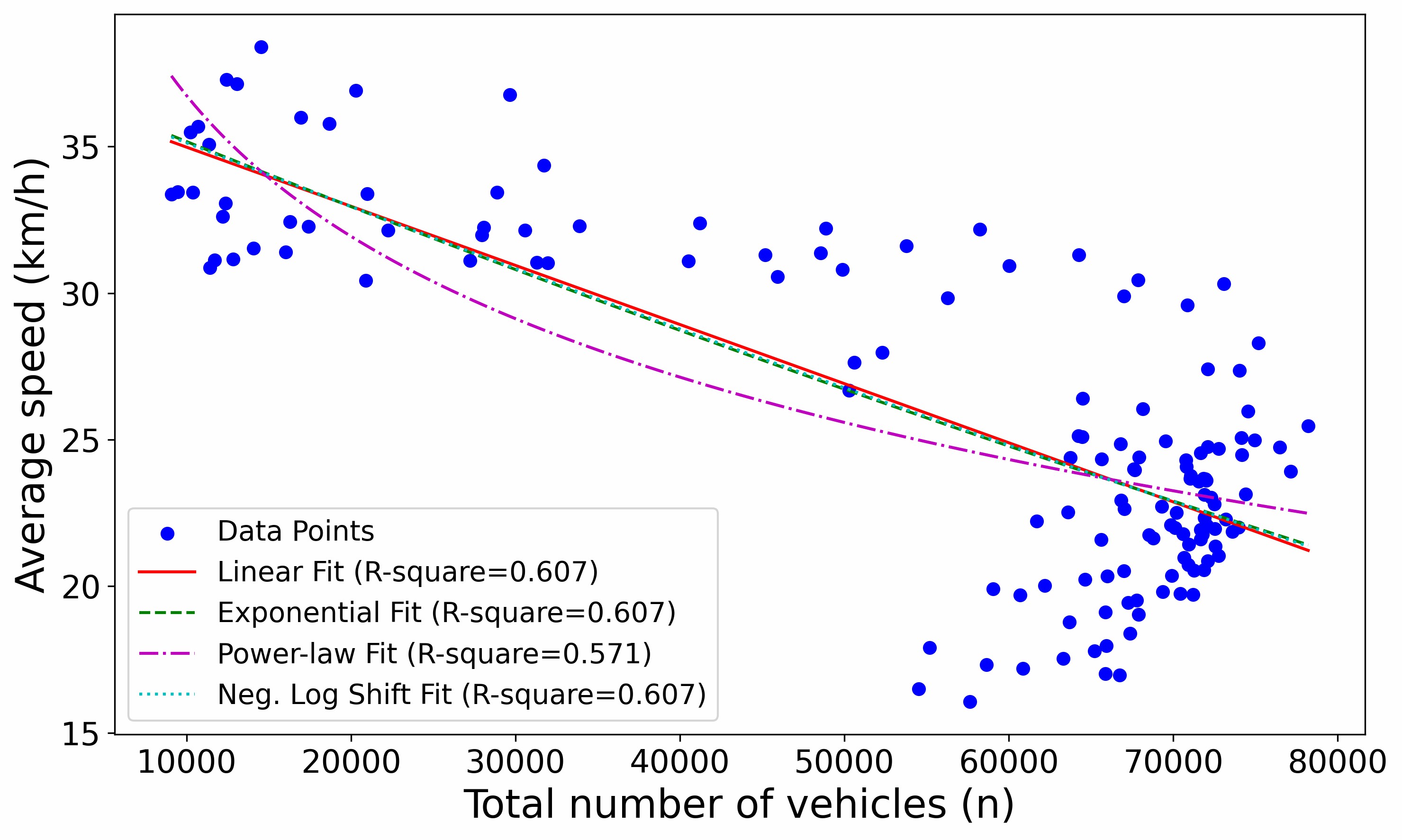}
    \caption{Comparison of model fittings for MFD. The scatter points represent the aggregated data over 30-minute intervals. Four candidate models are tested: linear, exponential, power-law, and shifted logarithmic functions.}
    \label{fig:MFD}
\end{figure}

We evaluate each functional form using the coefficient of determination ($R^2$). Among all candidate functions, the exponential and log-shift models yield fitted curves that are nearly linear. We select the linear model due to its sufficient predictive accuracy ($R^2 = 0.61$) and simplicity, as shown in Table~\ref{tab:mfd_models} and Figure~\ref{fig:MFD}. The linear model is used to estimate the average travel speed.

\begin{table}[!htbp]
\centering
\caption{Performance of candidate MFD models}
\label{tab:mfd_models}
\begin{tabular}{l c c}
\toprule
\textbf{Model} & \textbf{Parameters} & \textbf{$R^2$} \\
\midrule
Linear & $A = 36.998$, $B = -0.000202$ & 0.61 \\
Exponential & $a = 69.868$, $b = 3 \times 10^{-6}$, $c = -32.403$ & 0.61 \\
Power law & $a = 52429.814$, $b = 0.000$, $c = -52329.194$ & 0.5708 \\
Log-shift & $a = 1063.048$, $b = 80.264$, $c = 354669.766$ & 0.61 \\
\bottomrule
\end{tabular}
\end{table}

\section{Identification of high-demand and high-attraction nodes}
\label{app:high_nodes}

To analyse spatial disparities in trip demand, we classify intersections on the road network by two characteristics: demand and attractiveness. If the number of trip origins at an intersection (node) is statistically significantly higher than that at other intersections, this intersection is defined as a high-demand node; otherwise, it is a non-high-demand node. If the number of trip destinations at an intersection node is statistically significantly higher than that at other intersections, this intersection is defined as a high-attraction node; otherwise, it is a non-high-attraction node.

We use a statistical outlier detection method (Algorithm~\ref{alg:high_node_classification}) based on the log-transformed origin and destination counts. A node is considered significantly high-demand/high-attraction if its Z-score exceeds a predefined threshold $z_{\text{threshold}}$. This approach ensures that only nodes with demand/attractiveness far exceeding the average can be classified as high-demand/high-attraction. We set $z_{\text{threshold}} = 1.2816$, which corresponds to the 90\% quantile in a one-tailed Z-test ($\alpha = 0.10$), marking the number of trips at a node as statistically significantly higher than the average origin or destination counts of all nodes.

\vspace{1em}

\begin{algorithm}[H]
\caption{Identification of High-Demand and High-Attraction Nodes}
\label{alg:high_node_classification}
\begin{algorithmic}[1]
  \State \textbf{Input:} Origin (pickup) counts $\texttt{on\_count}(n)$ and destination (drop-off) counts $\texttt{off\_count}(n)$ for all nodes $n \in \mathcal{N}$
  \State \textbf{Output:} Node $\texttt{is\_high\_demand}(n)$ and $\texttt{is\_high\_attractive}(n)$ for all $n$
  \Statex

  \State Compute Z-scores of log-transformed pickup counts:
    $z^{\text{on}}_n \gets \text{Zscore}(\log(\texttt{1 + on\_count}(n)))$
  \State Compute Z-scores of log-transformed drop-off counts:
    $z^{\text{off}}_n \gets \text{Zscore}(\log(\texttt{1 + off\_count}(n)))$

  \For{each node $n \in \mathcal{N}$}
    \If{$z^{\text{on}}_n \geq z_{\text{threshold}}$}
      \State $\texttt{is\_high\_demand}(n) \gets 1$
    \Else
      \State $\texttt{is\_high\_demand}(n) \gets 0$
    \EndIf

    \If{$z^{\text{off}}_n \geq z_{\text{threshold}}$}
      \State $\texttt{is\_high\_attraction}(n) \gets 1$
    \Else
      \State $\texttt{is\_high\_attraction}(n) \gets 0$
    \EndIf
  \EndFor
\end{algorithmic}
\end{algorithm}

\bibliographystyle{elsarticle-harv}
\bibliography{reference_2}

@article{Rayle2016,
  author  = {Lisa Rayle and Danielle Dai and Nelson Chan and Robert Cervero and Susan Shaheen},
  title   = {Just a Better Taxi? A Survey-Based Comparison of Taxis, Transit, and Ridesourcing Services in San Francisco},
  journal = {Transport Policy},
  volume  = {45},
  pages   = {168--178},
  year    = {2016},
  doi     = {10.1016/j.tranpol.2015.10.004}
}

@article{Vazifeh2018,
  author  = {Mohammad M. Vazifeh and Paolo Santi and Giovanni Resta and Steven H. Strogatz and Carlo Ratti},
  title   = {Addressing the minimum fleet problem in on-demand urban mobility},
  journal = {Nature},
  volume  = {557},
  pages   = {534--538},
  year    = {2018},
  doi     = {10.1038/s41586-018-0095-1}
}

@article{FengKongWang2021,
  author  = {Guiyun Feng and Guangwen Kong and Zizhuo Wang},
  title   = {We Are on the Way: Analysis of On‑Demand Ride‑Hailing Systems},
  journal = {Manufacturing \& Service Operations Management},
  volume  = {23},
  number  = {5},
  pages   = {1237--1256},
  year    = {2021},
  doi     = {10.1287/msom.2020.0880}
}

@article{Beojone2021_inefficiency_ride_sourcing,
  title     = {On the inefficiency of ride\-sourcing services towards urban congestion},
  author    = {Beojone, Caio Vitor and Geroliminis, Nikolas},
  journal   = {Transportation Research Part C: Emerging Technologies},
  volume    = {124},
  pages     = {102890},
  year      = {2021},
  doi       = {10.1016/j.trc.2021.102890}
}

@article{JonkerVolgenant1986_hungarian_improvement,
  author    = {Jonker, Roy and Volgenant, Ton},
  title     = {Improving the Hungarian assignment algorithm},
  journal   = {Operations Research Letters},
  volume    = {5},
  number    = {4},
  pages     = {171--175},
  year      = {1986},
  month     = oct,
  doi       = {10.1016/0167-6377(86)90073-8}
}

@article{Guo2019_Chengdu_GPS,
  title     = {Urban link travel speed dataset from a megacity road network},
  author    = {Guo, Feng and Zhang, Dongqing and Dong, Yucheng and Guo, Zhao Xia},
  journal   = {Scientific Data},
  volume    = {6},
  pages     = {61},
  year      = {2019},
  doi       = {10.1038/s41597-019-0060-3}
}

@article{DeCauwer2015,
  author       = {Cedric De Cauwer and Joeri Van Mierlo and Thierry Coosemans},
  title        = {Energy Consumption Prediction for Electric Vehicles Based on Real-World Data},
  journal      = {Energies},
  volume       = {8},
  number       = {8},
  pages        = {8573--8593},
  year         = {2015},
  doi          = {10.3390/en8088573}
}

@report{RolandBerger2015,
  title        = {Think Act: Urban Mobility in the Age of Mobile Internet},
  author       = {{Roland Berger Strategy Consultants}},
  year         = {2015},
  month        = {May},
  institution  = {Roland Berger Strategy Consultants},
  address      = {Beijing, China},
  url          = {https://www.rolandberger.com.cn}
}

@online{mot2016policy,
  author    = {{Ministry of Transport of the People's Republic of China}},
  title     = {Guiding Opinions on Deepening Reform and Advancing the Healthy Development of the Taxi Industry},
  year      = {2016},
  url       = {https://www.gov.cn/xinwen/2016-07/28/content_5097333.htm}
}

@article{LiangLuo2024,
  author  = {Liang, Yile and Luo, Haocheng and Duan, Haining and Li, Donghui and Liao, Hongsen and Feng, Jie and Zhao, Jiuxia and Ren, Hao and Ding, Xuetao and Cha, Ying and Zhou, Qingte and Situ, Chenqi and Hao, Jinghua and Xing, Ke and Yin, Feifan and He, Renqing and Sun, Yang and Zheng, Yueqiang and Feng, Yipeng and Sun, Zhizhao and Chen, Jingfang and Zheng, Jie and Wang, Ling},
  title   = {Meituan's Real-Time Intelligent Dispatching Algorithms Build the World's Largest Minute-Level Delivery Network},
  journal = {INFORMS Journal on Applied Analytics},
  volume  = {54},
  number  = {1},
  pages   = {84--101},
  year    = {2024},
  doi     = {10.1287/inte.2023.0084}
}

@article{Oh2020_AMOD_Singapore,
  title        = {Assessing the impacts of automated mobility‑on‑demand through agent‑based simulation: A study of Singapore},
  author       = {Simon Oh and Ravi Seshadri and Carlos Lima Azevedo and Nishant Kumar and Kakali Basak and Moshe Ben‑Akiva},
  journal      = {Transportation Research Part A: Policy and Practice},
  year         = {2020},
  volume       = {138},
  pages        = {367--388},
  doi          = {10.1016/j.tra.2020.06.004}
}

@article{Sun2019RidesourcingPolicy,
  title     = {Spatiotemporal evolution of ridesourcing markets under the new restriction policy: A case study in Shanghai},
  author    = {Sun, Daniel (Jian) and Ding, Xueqing},
  journal   = {Transportation Research Part A: Policy and Practice},
  volume    = {130},
  pages     = {227--239},
  year      = {2019},
  publisher = {Elsevier},
  doi       = {10.1016/j.tra.2019.09.052}
}

@article{ruder2016overview,
  title={An overview of gradient descent optimization algorithms},
  author={Ruder, Sebastian},
  journal={arXiv preprint arXiv:1609.04747},
  year={2016}
}

@article{floyd1962algorithm,
  title={Algorithm 97: Shortest path},
  author={Floyd, Robert W.},
  journal={Communications of the ACM},
  volume={5},
  number={6},
  pages={345},
  year={1962},
  publisher={ACM}
}

@article{Aguilera-Garcia2022,
  author       = {Aguilera-Garc{\'\i}a, {\'A}lvaro and G{\'o}mez, Juan and Vel{\'a}zquez, Guillermo and Vassallo, Jos{\'e} Manuel},
  title        = {Ridesourcing vs.\ traditional taxi services: Understanding users’ choices and preferences in Spain},
  journal      = {Transportation Research Part A: Policy and Practice},
  volume       = {155},
  pages        = {161--178},
  year         = {2022},
  doi          = {10.1016/j.tra.2021.11.002}
}

@article{siddiq2022ride,
  title={Ride-Hailing Platforms: Competition and Autonomous Vehicles},
  author={Siddiq, Auyon and Taylor, Terry A.},
  journal={Manufacturing \& Service Operations Management},
  volume={24},
  number={3},
  pages={1511--1528},
  year={2022},
  publisher={INFORMS},
  doi={10.1287/msom.2021.1013}
}

@article{NahmiasBiran2019TRR,
  title   = {From Traditional to Automated Mobility on Demand: A Comprehensive Framework for Modeling On-Demand Services in SimMobility},
  author  = {Nahmias-Biran, Bat Hen and Oke, Jimi B. and Kumar, Nishant and Basak, Kakali and Araldo, Andrea and Seshadri, Ravi and Akkinepally, Arun and Lima Azevedo, Carlos M. and Ben-Akiva, Moshe},
  journal = {Transportation Research Record: Journal of the Transportation Research Board},
  year    = {2019},
  volume  = {2673},
  number  = {12},
  pages   = {15--29},
  doi     = {10.1177/0361198119853553}
}

@online{Waymo2025YIR,
  title        = {Delivering more for our riders in a year of incredible growth},
  author       = {{The Waymo Team}},
  organization = {Waymo LLC},
  year         = {2025},
  url          = {https://waymo.com/blog/2025/12/2025-year-in-review}
}

@article{Lokhandwala2018TRC,
  title   = {Dynamic ride sharing using traditional taxis and shared autonomous taxis: A case study of {NYC}},
  author  = {Lokhandwala, Mustafa and Cai, Hua},
  journal = {Transportation Research Part C: Emerging Technologies},
  year    = {2018},
  volume  = {97},
  pages   = {45--60},
  doi     = {10.1016/j.trc.2018.10.007}
}

@article{Pakusch2020TBS,
  title   = {Traditional taxis vs automated taxis -- Does the driver matter for Millennials?},
  author  = {Pakusch, Christina and Meurer, Johanna and Tolmie, Peter and Stevens, Gunnar},
  journal = {Travel Behaviour and Society},
  year    = {2020},
  volume  = {21},
  pages   = {214--225},
  doi     = {10.1016/j.tbs.2020.06.009}
}

@article{Pemberton2025BIT,
  title   = {Choosing a driverless or conventional taxi: a quick or slow decision?},
  author  = {Pemberton, Laura and Sutcliffe, Alistair and Mehandjiev, Nik},
  journal = {Behaviour \& Information Technology},
  year    = {2025},
  volume  = {44},
  number  = {20},
  pages   = {5048--5067},
  doi     = {10.1080/0144929X.2025.2504515}
}

@article{Bosch2018TP,
  title   = {Cost-based analysis of autonomous mobility services},
  author  = {B{\"o}sch, Patrick M. and Becker, Felix and Becker, Henrik and Axhausen, Kay W.},
  journal = {Transport Policy},
  year    = {2018},
  volume  = {64},
  pages   = {76--91},
  doi     = {10.1016/j.tranpol.2017.09.005}
}

@article{Freedman2018IE,
  title   = {Autonomous vehicles are cost-effective when used as taxis},
  author  = {Freedman, Isaac G. and Kim, Ellen and Muennig, Peter A.},
  journal = {Injury Epidemiology},
  year    = {2018},
  volume  = {5},
  number  = {1},
  pages   = {24},
  doi     = {10.1186/s40621-018-0153-z}
}

@article{Wadud2021TRA,
  title   = {Fully automated vehicles: A cost-based analysis of the share of ownership and mobility services, and its socio-economic determinants},
  author  = {Wadud, Zia and Mattioli, Giulio},
  journal = {Transportation Research Part A: Policy and Practice},
  year    = {2021},
  volume  = {151},
  pages   = {228--244},
  doi     = {10.1016/j.tra.2021.06.024}
}

@article{Nahmias-Biran2021Transportation,
  title        = {Evaluating the impacts of shared automated mobility on-demand services: an activity-based accessibility approach},
  author       = {Nahmias-Biran, Bat Hen and Oke, Jimi B. and Kumar, Nishant and Lima Azevedo, Carlos M. and Ben-Akiva, Moshe},
  journal      = {Transportation},
  year         = {2021},
  volume       = {48},
  number       = {4},
  pages        = {1613--1638},
  doi          = {10.1007/s11116-020-10106-y},
}

@article{SimoniKockelmanGurumurthyBischoff2019,
  title        = {Congestion pricing in a world of self-driving vehicles: An analysis of different strategies in alternative future scenarios},
  author       = {Simoni, Michele D. and Kockelman, Kara M. and Gurumurthy, Krishna M. and Bischoff, Joschka},
  journal      = {Transportation Research Part C: Emerging Technologies},
  year         = {2019},
  volume       = {98},
  pages        = {167--185},
  doi          = {10.1016/j.trc.2018.11.002},
}

@article{SalazarLanzettiRossiSchifferPavone2020,
  title        = {Intermodal Autonomous Mobility-on-Demand},
  author       = {Salazar, Mauro and Lanzetti, Nicolas and Rossi, Federico and Schiffer, Maximilian and Pavone, Marco},
  journal      = {IEEE Transactions on Intelligent Transportation Systems},
  year         = {2020},
  volume       = {21},
  number       = {9},
  pages        = {3946--3960},
  doi          = {10.1109/TITS.2019.2950720},
}

@article{DengChenShaoTangPiGupta2022JMSE,
  title        = {Multi-objective vehicle rebalancing for ridehailing system using a reinforcement learning approach},
  author       = {Deng, Yuntian and Chen, Hao and Shao, Shiping and Tang, Jiacheng and Pi, Jianzong and Gupta, Abhishek},
  journal      = {Journal of Management Science and Engineering},
  year         = {2022},
  volume       = {7},
  number       = {2},
  pages        = {346--364},
  doi          = {10.1016/j.jmse.2021.12.004},
}

@article{KontouGarikapatiHou2020TRC,
  title        = {Reducing ridesourcing empty vehicle travel with future travel demand prediction},
  author       = {Kontou, Eleftheria and Garikapati, Venu and Hou, Yi},
  journal      = {Transportation Research Part C: Emerging Technologies},
  year         = {2020},
  volume       = {121},
  pages        = {102826},
  doi          = {10.1016/j.trc.2020.102826},
}

@article{WardMichalekAzevedoSamarasFerreira2019TRC,
  title        = {Effects of on-demand ridesourcing on vehicle ownership, fuel consumption, vehicle miles traveled, and emissions per capita in U.S. states},
  author       = {Ward, Jacob W. and Michalek, Jeremy J. and Azevedo, In{\^e}s L. and Samaras, Constantine and Ferreira, Pedro},
  journal      = {Transportation Research Part C: Emerging Technologies},
  year         = {2019},
  volume       = {108},
  pages        = {289--301},
  doi          = {10.1016/j.trc.2019.07.026},
}

@article{Yang2024Robotaxi,
  author       = {Yang, Jerry and Yeh, Andrew},
  title        = {Chinese Apollo Go commands robotaxi market adoption amid Tesla delay reports},
  journal      = {Digitimes Asia},
  year         = {2024},
  url          = {https://www.digitimes.com/news/a20240717PD217/china-robotaxi-av-market-tesla.html}
}

@misc{Baidu2025HKEXQ2,
  author       = {{Baidu, Inc.}},
  title        = {Announcement of the 2025 Second Quarter Results},
  year         = {2025},
  month        = aug,
  day          = {20},
  url          = {https://staticpacific.blob.core.windows.net/press-releases-attachments./3911462/HKEX-EPS_20250820_11802582_0.PDF},
  note         = {HKEX filing PDF (accessed 2026-02-02)}
}

@article{Zhou2023RTE,
  author    = {Zhou, Yimin and Xu, Meng},
  title     = {Robotaxi service: The transition and governance investigation in {China}},
  journal   = {Research in Transportation Economics},
  volume    = {100},
  pages     = {101321},
  year      = {2023},
  doi       = {10.1016/j.retrec.2023.101321},
}

@article{Pande2023JEPR,
  author    = {Pande, Devyani and Taeihagh, Araz},
  title     = {Navigating the governance challenges of disruptive technologies: insights from regulation of autonomous systems in {Singapore}},
  journal   = {Journal of Economic Policy Reform},
  year      = {2023},
  volume    = {26},
  number    = {4},
  pages     = {523--540},
  doi       = {10.1080/17487870.2023.2197599},
}

@article{Ding2022TRC,
  title   = {Dynamic dispatch of connected taxis for large-scale urban road networks with stochastic demands: An {MFD}-enabled hierarchical and cooperative approach},
  author  = {Ding, Heng and Li, Jiye and Zheng, Nan and Zheng, Xiaoyan and Huang, Wenjuan and Bai, Haijian},
  journal = {Transportation Research Part C: Emerging Technologies},
  volume  = {142},
  pages   = {103792},
  year    = {2022},
  doi     = {10.1016/j.trc.2022.103792}
}

@article{Asgari2018TRR,
  title        = {A Stated Preference Survey Approach to Understanding Mobility Choices in Light of Shared Mobility Services and Automated Vehicle Technologies in the US},
  author       = {Asgari, Hamidreza and Jin, Xia and Corkery, Terrence},
  journal      = {Transportation Research Record},
  year         = {2018},
  volume       = {2672},
  number       = {47},
  pages        = {12--22},
  doi          = {10.1177/0361198118790124},
}

@article{Mo2021TRC,
  title        = {Impacts of subjective evaluations and inertia from existing travel modes on adoption of autonomous mobility-on-demand},
  author       = {Mo, Baichuan and Wang, Qing Yi and Moody, Joanna and Shen, Yu and Zhao, Jinhua},
  journal      = {Transportation Research Part C: Emerging Technologies},
  year         = {2021},
  volume       = {130},
  pages        = {103281},
  doi          = {10.1016/j.trc.2021.103281}
}

@article{Brown2021Transportation,
  title   = {Hailing a change: comparing taxi and ridehail service quality in {Los Angeles}},
  author  = {Brown, Anne and LaValle, Whitney},
  journal = {Transportation},
  volume  = {48},
  number  = {2},
  pages   = {1007--1031},
  year    = {2021},
  doi     = {10.1007/s11116-020-10086-z}
}

@article{Ruch2021JAT,
  author  = {Ruch, Claudio and H{\"o}rl, Sebastian and G{\"a}chter, Joel and Hakenberg, Jan},
  title   = {The Impact of Fleet Coordination on Taxi Operations},
  journal = {Journal of Advanced Transportation},
  volume  = {2021},
  pages   = {1--14},
  year    = {2021},
  doi     = {10.1155/2021/2145716}
}

@article{Zhan2016TITS,
  title   = {A graph-based approach to measuring the efficiency of an urban taxi service system},
  author  = {Zhan, Xianyuan and Qian, Xinwu and Ukkusuri, Satish V.},
  journal = {IEEE Transactions on Intelligent Transportation Systems},
  volume  = {17},
  number  = {9},
  pages   = {2479--2489},
  year    = {2016},
  doi     = {10.1109/TITS.2016.2521862}
}

@article{FagnantKockelman2018Transportation,
  author  = {Fagnant, Daniel J. and Kockelman, Kara M.},
  title   = {Dynamic ride-sharing and fleet sizing for a system of shared autonomous vehicles in Austin, Texas},
  journal = {Transportation},
  year    = {2018},
  volume  = {45},
  number  = {1},
  pages   = {143--158},
  doi     = {10.1007/s11116-016-9729-z}
}

@article{HylandMahmassani2020TRA,
  author  = {Hyland, Michael and Mahmassani, Hani S.},
  title   = {Operational benefits and challenges of shared-ride automated mobility-on-demand services},
  journal = {Transportation Research Part A: Policy and Practice},
  year    = {2020},
  volume  = {134},
  pages   = {251--270},
  doi     = {10.1016/j.tra.2020.02.017}
}

@article{Cai2025TIMOD,
  author  = {Cai, Mingming and Abu Ashour, Lamis and Shen, Qing and Chen, Cynthia},
  title   = {Incorporating Mobility-on-Demand into Public Transit in Suburban Areas: A Comparative Cost-Effectiveness Evaluation},
  journal = {Transportation Research Part D: Transport and Environment},
  volume  = {144},
  pages   = {104775},
  year    = {2025},
  doi      = {10.1016/j.trd.2025.104775},
}

\end{document}